\documentclass[reqno]{amsart}
\usepackage[T1]{fontenc}
\usepackage{amsxtra}
\usepackage{enumerate}
\usepackage{color}
\usepackage{amscd,amsfonts,amsmath,amssymb,
amsthm,amstext,latexsym}
\usepackage{pifont,graphicx}
\usepackage[english]{babel}
\usepackage[all,cmtip]{xy} 
\usepackage{cases} 
\DeclareFontFamily{U}{mathx}{\hyphenchar\font45}
\DeclareFontShape{U}{mathx}{m}{n}{
      <5> <6> <7> <8> <9> <10>
      <10.95> <12> <14.4> <17.28> <20.74> <24.88>
      mathx10
      }{}
\DeclareSymbolFont{mathx}{U}{mathx}{m}{n}
\DeclareFontSubstitution{U}{mathx}{m}{n}
\DeclareMathAccent{\widecheck}{0}{mathx}{"71}
\DeclareMathAccent{\widetilde}{0}{mathx}{"72}
\DeclareMathAccent{\widebar}{0}{mathx}{"73}
\DeclareMathAccent{\widevec}{0}{mathx}{"74}
\DeclareMathAccent{\widehat}{0}{mathx}{"70}
\newcommand{\A}{{\mathbb A}}
\newcommand{\Z}{{\mathbb Z}}
\renewcommand{\S}{{\mathbb S}}
\newcommand{\R}{{\mathbb R}}
\newcommand{\C}{{\mathbb C}}
\newcommand{\D}{{\mathbb D}}
\newcommand{\N}{{\mathbb N}}
\def\boA{{\mathcal A}}
\def\boB{{\mathcal B}}

\def\CC{\overline{\C}}
\def\boE{{\mathcal E}}
\def\boD{{\mathcal D}}
\def\boF{{\mathcal F}}
\def\boG{{\mathcal G}}
\def\boH{{\mathcal H}}
\def\boL{{\mathcal L}}
\def\boM{{\mathcal M}}
\def\boP{{\mathcal P}}
\def\boR{{\mathcal R}}

\def\boW{{\mathcal W}}
\def\cqfd{\hfill$\Box$}
\def\Res{{\,\rm Res}}

\def\Re{{\,\rm Re}}
\def\Im{{\,\rm Im}}
\def\[{[\![}
\def\]{]\!]}
\def\ii{{\rm i}}
\def\wtalpha{\widetilde{\alpha}}

\def\whgamma{{\widehat{\gamma}}}
\def\whPhi{{\widehat{\Phi}}}

\def\whxi{{\widehat{\xi}}}

\def\wtboA{\widetilde{\mathcal{A}}}

\def\whM{{\widehat{M}}}

\def\whP{{\widehat{P}}}
\def\whf{\widehat{f}}
\def\whalpha{{\widehat{\alpha}}}
\def\whbeta{{\widehat{\beta}}}
\def\whgamma{{\widehat{\gamma}}}
\def\wt0{\widetilde{0}}
\def\wta{{\widetilde{a}}}
\def\wtb{{\widetilde{b}}}
\def\wtc{{\widetilde{c}}}
\def\wtz{\widetilde{z}}
\def\wtbeta{\widetilde{\beta}}
\def\wtA{{\widetilde{A}}}
\def\wtB{{\widetilde{B}}}
\def\wtC{{\widetilde{C}}}
\def\wtF{{\widetilde{F}}}
\def\wtg{{\widetilde{g}}}

\def\wtgamma{{\widetilde{\gamma}}}
\def\wtp{\widetilde{p}}
\def\wtP{{\widetilde{P}}}
\def\wtPhi{{\widetilde{\Phi}}}
\def\wtf{{\widetilde{f}}}
\def\wtt{\widetilde{t}}
\def\wtw{\widetilde{w}}
\def\wtg{{\widetilde{g}}}
\def\wtxi{{\widetilde{\xi}}}
\def\wtzeta{\widetilde{\zeta}}

\def\wtM{{\widetilde{M}}}
\def\wtOmega{{\widetilde{\Omega}}}
\def\wtr{\widetilde{r}}
\def\wtq{\widetilde{q}}

\def\whzeta{\widehat{\zeta}}
\def\whTheta{\widehat{\Theta}}
\def\whC{\widehat{C}}
\def\whchi{\widehat{\chi}}

\def\wtSigma{\widetilde{\Sigma}}
\def\wcM{{\widecheck{M}}}

\def\wcxi{{\widecheck{\xi}}}
\def\wczeta{\widecheck{\zeta}}

\def\wcPhi{\widecheck{\Phi}}
\def\wcF{\widecheck{F}}
\def\wcB{\widecheck{B}}
\def\wcf{\widecheck{f}}
\def\wcN{\widecheck{N}}

\def\bft{{\bf t}}

\def\su{{\mathfrak su}}
\def\sl{{\mathfrak sl}}
\def\pos{{\geq 0}}
\def\pospos{{+}}
\def\neg{{\leq 0}}
\def\negneg{{-}}

\def\pij0{p_{ij}^{\circ}}
\def\Ord{\mbox{Ord}}
\def\Uni{\mbox{Uni}}
\def\Pos{\mbox{Pos}}
\def\force{{\bf f}}
\def\xiS{\xi^S}
\def\PhiS{\Phi^S}
\def\FS{F^S}
\def\BS{B^S}
\def\fS{f^S}
\def\NS{N^S}
\def\xiC{\xi^C}
\def\PhiC{\Phi^C}
\def\FC{F^C}
\def\BC{B^C}
\def\fC{f^C}
\def\NC{N^C}

\def\SSigma{\overline{\Sigma}}
\def\u{{\bf u}}
\def\v{{\bf v}}
\def\x{{\bf x}}
\def\y{{\bf y}}
\def\I{{I^*}}
\def\Uni{\mbox{Uni}}
\def\Pos{\mbox{Pos}}
\def\Sym{\mbox{Sym}}
\def\Nor{\mbox{Nor}}
\newcommand{\smallfrac}[2]{\mbox{$\frac{#1}{#2}$}}
\renewcommand{\matrix}[1]{\left(\begin{array}{cc} #1\end{array}\right)}

\newcommand{\minimatrix}[1]{\left(\begin{smallmatrix}#1\end{smallmatrix}\right)}

\newtheorem{theorem}{Theorem}
\newtheorem{lemma}{Lemma}
\newtheorem{proposition}{Proposition}
\newtheorem{remark}{Remark}
\newtheorem{corollary}{Corollary}
\newtheorem{claim}{Claim}

\newtheorem{definition}{Definition}

\newenvironment{myenumerate}{\begin{enumerate}[1.]}{\end{enumerate}}

\begin{document}
\title{Opening nodes and the DPW method}
\author{Martin Traizet}
\maketitle
{\em Abstract: we combine the DPW method and Opening Nodes to construct
embedded surfaces of positive constant mean curvature with Delaunay ends in euclidean
space, with no limitation to the genus or number of ends.}
\medskip

\section{Introduction}
In \cite{dorfmeister-pedit-wu}, Dorfmeister, Pedit and Wu have shown that harmonic maps
from a Riemann surface to a symmetric space admit a Weierstrass-type
representation, which means that they can be represented in terms of
holomorphic data.
In particular, surfaces with constant mean curvature one (CMC-1 for short) in
euclidean space admit such a representation, owing to the fact that the Gauss map
of a CMC-1 surface is a harmonic map to the 2-sphere.
This representation is now called the DPW method and has been widely used
to construct CMC-1 surfaces in $\R^3$ and also constant mean curvature surfaces in homogeneous spaces such as the sphere ${\mathbb S}^3$ or hyperbolic space ${\mathbb  H}^3$:
see for example \cite{dorfmeister-haak,
dorfmeister-wu,
heller-heller-schmitt,
heller1,
heller2,
heller3,
kilian-kobayashi-rossman-schmitt,
kilian-mcintosh-schmitt,
schmitt,
schmitt-kilian-kobayashi-rossman}.
Also the DPW method has been implemented by N. Schmitt to make computer
images of CMC-1 surfaces.
\medskip

The main limitation to the construction of examples is the Monodromy
Problem,
so either the topology of the constructed examples is limited
or symmetries are imposed to the construction, in order to reduce the number of equations to be solved.
In constract, Kapouleas \cite{kapouleas} has constructed embedded CMC-1 surfaces with
no limitation on the genus or number of ends by gluing round spheres and
pieces of Delaunay surfaces, using Partial Differential Equations techniques.
Our goal in this paper is to carry on the construction of Kapouleas using the DPW method.
\medskip

The underlying Riemann surface is defined by Opening Nodes, which is a model for
Riemann surfaces with ``small necks''.
The theory of Opening Nodes has been used by the author to construct minimal surfaces
in euclidean space via the classicial Weierstrass Representation
(see for example \cite{nosym} or \cite{triply})
or CMC-1 surfaces in hyperbolic space via Bryant Representation
\cite{bryant}.
This paper opens up the possibility of opening nodes in the DPW method.
Since the DPW method can be used to construct CMC surfaces in all space forms, or
more generally harmonic maps from a Riemann surface to a  symmetric space,
there should be many applications.
\tableofcontents
\section{Main result}
\label{section-intro}
\subsection{Weighted graphs}
\label{section-intro-graph}
\begin{figure}
\begin{center}
\includegraphics[height=4.5cm]{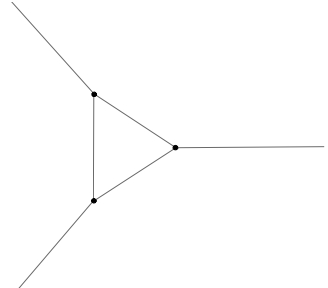}
\hspace{0.5cm}
\includegraphics[height=4.5cm]{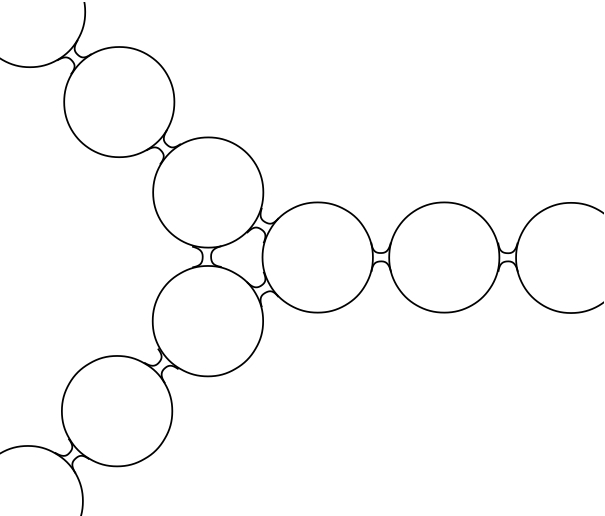}
\end{center}
\caption{A balanced graph with three length-2 edges and three rays (left) and a CMC-1 surface with genus one and three Delaunay-type ends
in the corresponding family (right).}
\label{fig-genus1}
\end{figure}

We want to construct CMC-1 surfaces by gluing spheres and half-Delaunay surfaces.
The layout of these pieces is described by a weighted graph in euclidean space.
\begin{definition}
A weighted graph $\Gamma$ is a triple $(V,E,R)$, where
\begin{myenumerate}
\item $V=\{\v_1,\cdots,\v_N\}$ is the set of vertices: each vertex $\v_i$ is a point in euclidean space.
\item $E=\{e_{ij}\}_{(i,j)\in I}$ is the set of edges: for each $(i,j)\in I$, $i<j$ and $e_{ij}$ is the segment from $\v_i$ to $\v_j$. The edge $e_{ij}$ is assigned a non-zero weight $\tau_{ij}$.
The length of the edge $e_{ij}$ is denoted $\ell_{ij}$.
\item $R=\{r_1,\cdots,r_n\}$ is the set of rays: each ray $r_k$ is a half-line issued from a vertex and is assigned a non-zero weight $\tau_k$.
\end{myenumerate}
\end{definition}
Given a weighted graph $\Gamma$ all whose edges have even length, we can construct a singular surface $M_0$ as follows:
\begin{myenumerate}
\item For $i\in[1,N]$, place a radius-1 sphere centered at the vertex $\v_i$.
\item For each $(i,j)\in I$, place $\frac{1}{2}(\ell_{ij}-2)$
radius-1 spheres centered at the points on the edge $e_{ij}$ which are at even distance
from $\v_i$ (this connects the spheres centered at $\v_i$ and $\v_j$ by a necklace of
spheres).
\item For each ray $r_k$ issued from the vertex $\v_i$, place an infinite number of radius-1 spheres centered at the points on the half-line $r_k$ which are at even distance from $\v_i$.
\end{myenumerate}
Our goal in this paper is to construct a family of CMC-1 surfaces $(M_t)_{0<t<\epsilon}$ by desingularizing $M_0$, replacing
all tangency points between adjacent spheres by small catenoidal necks
(see Figure \ref{fig-genus1}).
The neck-sizes should be approximately $t\tau_{ij}$ in the case of edges and 
$t\tau_k$ in the case of rays.
(This is only a heuristic way to describe the result, and not the way we will construct $M_t$.)

\subsection{Forces}
\label{section-intro-forces}
The weighted graph $\Gamma$ must satisfy a balancing condition
for the construction to succed. For $i\in [1,N]$:
\begin{itemize}
\item Let $E_i$ be the set of indices $j\in[1,N]$ such that
the vertices $\v_i$ and $\v_j$ are connected by an edge, that is $(i,j)\in I$ or $(j,i)\in I$.
For $j\in E_i$, let
$$\u_{ij}=\frac{\v_j-\v_i}{||\v_j-\v_i||}$$
be the unitary vector in the direction of the edge.
Finally, in case $i>j$, define $\tau_{ij}=\tau_{ji}$.
\item Let $R_i$ be the set of indices $k\in[1,n]$ such that $r_k$ is a half-line issued from
$\v_i$. For $k\in R_i$, let $\u_k$ be the unitary vector spanning the half-line $r_k$
\end{itemize}
We define the force $\force_i$ acting on the vertex $\v_i$ by
$$\force_i=\sum_{j\in E_i} \tau_{ij}\u_{ij}+\sum_{k\in R_i}\tau_k \u_k$$
\begin{definition}
A weighted graph $\Gamma$ is balanced if for all $i\in[1,N]$, $\force_i=0$.
\end{definition}
\subsection{Non-degeneracy}
\label{section-intro-nondegeneracy}
To apply the Implicit Function Theorem, we need to perturb our graph $\Gamma$ in order
to prescribe small variations of edge-lengths and forces.
The parameters available to deform $\Gamma$ are: the position of the vertices $\v_i$,
the direction of the rays $r_k$ and the weights of the edges and rays.
The abstract structure of the graph is fixed under this deformation (so the
edges are determined by the positions of the vertices).
\begin{definition}
A weighted graph $\Gamma$ is non-degenerate if the jacobian of
the map
$$\left((\force_i)_{1\leq i\leq N},(\ell_{ij})_{(i,j)\in I}\right)\in\R^{3N}\times\R^I$$ with respect to the parameters (vertices, rays, weights) is onto.
\end{definition}
\begin{remark} There are plenty of interesting examples of non-degenerate graphs in
\cite{kapouleas}.
Our non-degeneracy condition is slightly different from the
flexibility condition of Kapouleas (Definition 1.18 in \cite{kapouleas}).
It is, however, very natural for our implicit function approach.
\end{remark}
\begin{remark} If a weighted graph has no rays, then
$\sum_{i=1}^N\force_i=0$ so it is always degenerate. We will not construct compact CMC-1 surfaces in this paper.
\end{remark}
\subsection{Main result}
\label{section-intro-main}
\begin{theorem}
\label{thm-main}
Assume that $\Gamma$ has even-length edges, is balanced and non-degenerate.
There exists a smooth 1-parameter family of immersed CMC-1 surfaces $(M_t)_{0<t<\epsilon}$
with the following properties:
\begin{myenumerate}
\item $M_t$ converges as $t\to 0$ to $M_0$. The convergence is for the
Hausdorf distance on compact sets of $\R^3$.
\item $M_t$ is homeomorphic to a tubular neighborhood of $\Gamma$.
\item For each $k\in[1,n]$, $M_t$ has a Delaunay end with weight
$\simeq 8\pi t\tau_k$ and whose axis converges as $t\to 0$ to the ray $r_k$.
\item If $\Gamma$ is pre-embedded, then $M_t$ is embedded.
\end{myenumerate}
\end{theorem}
\begin{definition}
\label{def-preembedded}
Following Kapouleas (Definition 2.2 in \cite{kapouleas}),
we say that $\Gamma$ is pre-embedded if:
\begin{myenumerate}
\item All weights are positive.
\item The distance between any two edges or rays which have no common endpoint
is greater than $2$.
\item The angle between any two edges or rays with a common endpoint is greater than $60^\circ$.
\end{myenumerate}
\end{definition}
As already said, Theorem \ref{thm-main} was proved using a completely different method in
\cite{kapouleas}.
In the simplest case $N=1$, there is no need for Opening Nodes and Theorem \ref{thm-main} is proved using the DPW method in \cite{nnoids}. We follow the same strategy
to define the DPW potential and we will use some of the results in \cite{nnoids}.
\subsection{Reduction to length-2 edges}
Let $\Gamma$ be a graph with $N$ vertices and even-length edges. Assume that $\Gamma$ has an edge $e_{ij}$ of length $\ell_{ij}\geq 4$. We can define a new graph $\widetilde{\Gamma}$ as follows: insert a new vertex $\v_{N+1}$
on the edge $e_{ij}$ at distance $2$ from $\v_i$. Replace the edge $e_{ij}$ by
the edges $[\v_i,\v_{N+1}]$ and $[\v_{N+1},\v_j]$, with respective lenghts
$2$ and $\ell_{ij}-2$. Assign to each new edge the weight $\tau_{ij}$.
The new graph $\widetilde{\Gamma}$ is clearly balanced.
\begin{proposition}
\label{prop-reduction}
If $\Gamma$ is non-degenerate then $\widetilde{\Gamma}$ is
non-degenerate. If $\Gamma$ is pre-embedded then $\widetilde{\Gamma}$
is pre-embedded.
\end{proposition}
The proof of Proposition \ref{prop-reduction} is elementary and is omitted.
Thanks to Proposition \ref{prop-reduction}, we can transform by induction the graph $\Gamma$ into
a balanced, non-degenerate graph with length-2 edges. 
Therefore, it suffices to prove Theorem \ref{thm-main} in the case where
all edges have length 2.
\section{Background}
\label{section-background}
\subsection{Opening nodes}
\label{section-background-openingnodes}
In this section, we recall the standard construction of opening nodes.
Consider $n$ copies of the Riemann sphere $\CC=\C\cup\{\infty\}$, labelled
$\CC_1,\cdots,\CC_n$. Consider $2m$ distinct points $p_1,\cdots,p_m$,
$q_1,\cdots,q_m$ in the disjoint union $\CC_1\cup\cdots\cup\CC_n$.
Identify $p_i$ with $q_i$ for all $i\in[1,m]$. This defines a Rieman surface with
nodes which we denote $\Sigma_0$. 
(The nodes refer to the double points $p_i\sim q_i$).
\medskip

To open nodes, consider local complex coordinates $v_i:V_i\stackrel{\sim}{\to} D(0,\varepsilon)$ in a neighborhood
of $p_i$ and $w_i:W_i\stackrel{\sim}{\to} D(0,\varepsilon)$ in a neighborhood of $q_i$, with
$v_i(p_i)=0$ and $w_i(q_i)=0$. We assume that the neighborhoods
$V_1,\cdots,V_m,W_1,\cdots,W_m$ are disjoint in $\CC_1\cup\cdots\cup \CC_n$.
Consider, for each $i\in[1,m]$, a complex parameter
$t_i$ with $|t_i|< \varepsilon^2$.
If $t_i=0$, identify $p_i$ with $q_i$ as above.
If $t_i\neq 0$, remove the disks $|v_i|\leq \frac{|t_i|}{\varepsilon}$ and $|w_i|\leq \frac{|t_i|}{\varepsilon}$. Identify
each point $z$ in the annulus $\frac{|t_i|}{\varepsilon}<|v_i|<\varepsilon$
with the point $z'$ in the annulus $\frac{|t_i|}{\varepsilon}<|w_i|<\varepsilon$
such that $$v_i(z)w_i(z')=t_i.$$
(In particular, the circle $|v_i|=|t_i|^{1/2}$ is identified with the circle $|w_i|=|t_i|^{1/2}$,
with the reverse orientation.)
This creates a Riemann surface (possibly with nodes)
which we denote $\Sigma_{\bft}$, where $\bft=(t_1,\cdots,t_n)$.
If $t_i\neq 0$, we can use $v_i$ and $w_i$ as local coordinates in $\Sigma_{\bft}$,
and the change of coordinate $\psi_i=w_i\circ v_i^{-1}$ is given by
$\psi_i(z)=\frac{t_i}{z}$.
When all $t_i$ are non-zero, $\Sigma_{\bft}$ is a genuine compact Riemann surface.
If $\Sigma_0$ is connected, its genus is $m-n+1$.
\medskip

One can define meromorphic 1-forms on a compact Riemann surface by prescribing
principal parts and periods. In the case of opening nodes, this can be formulated as follows.
\begin{definition}
\label{def-regulardifferential}
Let $r_1,\cdots,r_k$ be points in $\Sigma_0$, distinct from the nodes.
A meromorphic regular differential $\omega_0$ on $\Sigma_0$ with poles at $r_1,\cdots,r_k$
is a meromorphic 1-form
on the disjoint union $\CC_1\cup\cdots\cup \CC_n$ with poles
at $r_1,\cdots,r_k$ and simple poles at $p_i,q_i$ for $i\in[1,m]$ such that
the residues at $p_i$ and $q_i$ are opposite.\end{definition}
In the case $k=0$, $\omega_0$ is called a holomorphic regular differential, or simply
a regular differential.
\begin{theorem}
\label{thm-openingnodes}
Given a meromorphic regular differential $\omega_0$ on $\Sigma_0$, for $\bft$ in
a neighborhood of $0$,
there exists a unique meromorphic regular differential $\omega_{\bft}$ on $\Sigma_{\bft}$,
such that:
\begin{myenumerate}
\item For $i\in[1,m]$,
$$\int_{\gamma_i}\omega_{\bft}=\int_{\gamma_i}\omega_0$$
where $\gamma_i$ denotes the circle $|v_i|=\varepsilon$.
\item For $i\in[1,k]$, $\omega_{\bft}-\omega_0$
extends holomorphically at $r_i$.
\end{myenumerate}
Moreover, $\omega_{\bft}-\omega_0$ depends holomorphically on $\bft$ on compact
subsets of $\Sigma_0$ minus the nodes.
\end{theorem}
As indicated by the notation, when $\bft=0$, $\omega_{\bft}$ is equal to the given
meromorphic regular differential $\omega_0$.
In case all $t_i$ are non-zero, $\omega_{\bf t}$ is of course a genuine meromorphic
1-form on $\Sigma_{\bft}$.
Regarding the last point, a compact subset of $\Sigma_0$ minus the nodes is
included in $\Sigma_{\bft}$ for $\bft$ small enough.
\medskip

Theorem \ref{thm-openingnodes} was first proved by Fay \cite{fay} in the case of holomorphic regular differentials,
using sheaf theory.
This was extended to the case of meromorphic differentials with simple poles
by Masur \cite{masur}.
An elementary proof of the general case is given in \cite{crelle} using an
Implicit Function Theorem argument.
\begin{remark}
When all $t_i$ are non-zero, the existence of the meromorphic 1-form $\omega_{\bft}$
on $\Sigma_{\bft}$ follows from the
standard theory of compact Riemann surfaces. The content of Theorem \ref{thm-openingnodes} is
really that the limit of $\omega_{\bft}$ as $\bft\to 0$ exists and equals $\omega_0$.
\end{remark}
\medskip

An important point for our construction is that $\omega_0$ can be explicitely computed,
since meromorphic 1-forms on the Riemann sphere are rational fractions.
In fact, we can also compute the partial derivatives of $\omega_{\bft}$ with respect to $t_i$
at $\bft=0$, to any order. In this paper, we only need the first order derivative:
\begin{theorem}
\label{thm-openingnodes-derivee}
In compact subsets of $\Sigma_0$ minus the nodes,
the partial derivative $\frac{\partial \omega_{\bft}}{\partial t_i}|_{\bft=0}$
is equal to the unique meromorphic differential on the disjoint union $\CC_1\cup\cdots\cup \CC_n$
which has only two double poles at $p_i$ and $q_i$, with principal parts
$$\frac{-dv_i}{v_i^2}\Res_{q_i}\left(\frac{\omega_0}{w_i}\right)\qquad
\mbox{ at $p_i$}$$
$$\frac{-dw_i}{w_i^2}\Res_{p_i}\left(\frac{\omega_0}{v_i}\right)\qquad
\mbox{ at $q_i$}$$
\end{theorem}
This is proven in \cite{triply}, Lemma 3. See also \cite{crelle}, Remark 5.6.
\subsection{The DPW method}
\label{section-background-DPW}
In this section, we recall standard notations and results used in the
DPW method in the ``untwisted'' setting.
For a comprehensive introduction to the DPW method, we suggest
\cite{fujimori-kobayashi-rossman}.
\subsubsection{Loop groups}
We use blackboard letters for domains in the $\lambda$-plane.
For $\rho>1$, we denote:
\begin{itemize}
\item $\S^1$ the unit circle $\{\lambda\in\C:|\lambda|=1\}$,
\item $\D$ the unit disk $\{\lambda\in\C:|\lambda|<1\}$,
\item $\D_{\rho}$ the disk $\{\lambda\in\C:|\lambda|<\rho\}$,
\item $\D_{\rho}^*$ the punctured disk $\D_{\rho}\setminus\{0\}$,
\item $\A_{\rho}$ the annulus $\{\lambda\in\C:\rho^{-1}<|\lambda|<\rho\}$.
\end{itemize}
A loop is a smooth map from the unit circle to a matrix Lie group.
\begin{itemize}
\item If $G$ is a matrix Lie group (respectively a Lie algebra), $\Lambda G$ denotes the group (respectively the algebra) of smooth maps $\Phi:\S^1\to G$.
\item $\Lambda_+ SL(2,\C)\subset \Lambda SL(2,\C)$ is the subgroup of smooth maps $B:\S^1\to SL(2,\C)$ which extend holomorphically to the unit disk $\D$.
\item $\Lambda_+^{\R}SL(2,\C)$ is the set of $B\in \Lambda_+ SL(2,\C)$ such that
$B(0)$ is upper triangular with real elements on the diagonal.
\end{itemize}
\begin{theorem}[Iwasawa decomposition]
\label{thm-Iwasawa}
The multiplication $\Lambda SU(2)\times \Lambda_+^{\R} SL(2,\C)\to\Lambda SL(2,\C)$
is a diffeomorphism. The unique splitting of an element $\Phi\in\Lambda SL(2,\C)$ as
$\Phi=FB$ with $F\in\Lambda SU(2)$ and $B\in\Lambda_+^{\R} SL(2,\C)$
is called Iwasawa decomposition. $F$ is called the unitary factor of $\Phi$ and
denoted $\Uni(\Phi)$.
$B$ is called the positive factor and denoted $\Pos(\Phi)$.
\end{theorem}
\subsubsection{The matrix model of $\R^3$}
In the DPW method, one identifies $\R^3$ 
with the Lie algebra $\su(2)$ by
$$x=(x_1,x_2,x_3)\in\R^3\longleftrightarrow X=\frac{-\ii}{2}\matrix{
-x_3&x_1+\ii x_2\\x_1-\ii x_2 &x_3}\in\mathfrak{su}(2)$$
The group $SU(2)$ acts as linear isometries on $\su(2)$ by
$H\cdot X=HXH^{-1}$.
\subsubsection{The recipe}
\label{section-background-DPW-recipe}
The input data for the DPW method is a quadruple $(\Sigma,\xi,z_0,\Phi_0)$ where:
\begin{itemize}
\item $\Sigma$ is a Riemann surface.
\item $\xi$ is a $\Lambda\sl(2,\C)$-valued holomorphic 1-form
on $\Sigma$ called the DPW potential. More precisely,
$$\xi(z,\lambda)=\matrix{\alpha(z,\lambda)&\lambda^{-1}\beta(z,\lambda)\\
\gamma(z,\lambda)&-\alpha(z,\lambda)}$$
where $\alpha$, $\beta$, $\gamma$ are holomorphic 1-forms on $\Sigma$ with
respect to the $z$ variable and depend holomorphically on $\lambda$
in the disk $\D_{\rho}$ for some $\rho>1$.
\item $z_0\in\Sigma$ is a base point.
\item $\phi_0\in\Lambda SL(2,\C)$ is an initial condition.
\end{itemize}
Given this data, the DPW method is the following procedure.
Let $\wtSigma$ be the universal cover of $\Sigma$ and $\wtz_0\in\wtSigma$
be an arbitrary element in the fiber of $z_0$.
\begin{myenumerate}
\item Solve the Cauchy problem on $\wtSigma$:
\begin{equation}
\label{eq-cauchy}\left\{\begin{array}{l}
d\Phi(z,\lambda)=\Phi(z,\lambda)\xi(z,\lambda)\\
\Phi(\wtz_0,\lambda)=\phi_0(\lambda)\end{array}\right.\end{equation}
to obtain a solution
$\Phi:\wtSigma\to\Lambda SL(2,\C)$.
(The lift of $\xi$ to $\wtSigma$ is still denoted $\xi$.)
\item Compute, for $z\in\wtSigma$, the Iwasawa decomposition $(F(z,\cdot),B(z,\cdot))$ of
 $\Phi(z,\cdot)$.
\item Define $f:\wtSigma\to\su(2)\sim\R^3$ by the Sym-Bobenko formula:
\begin{equation}
\label{eq-sym}
f(z)=\ii\frac{\partial F}{\partial\lambda}(z,1)F(z,1)^{-1}=:\Sym(F(z,\cdot)).
\end{equation}
Then $f$ is a CMC-1 (banched) conformal immersion.
Its Gauss map is given by
\begin{equation}
\label{eq-normal}
N(z)=\frac{-\ii}{2}F(z,1)\matrix{1 &0 \\ 0 & -1}F(z,1)^{-1}=:\Nor(F(z,\cdot))
\end{equation}
In term of local conformal coordinates $z=x+\ii y$, the derivatives of $f$ are given by:
\begin{equation}
\label{eq-dfx}
\frac{\partial f}{\partial x}(z)=\frac{-\ii}{2}\;\mu(z) F(z,1)\matrix{0 &1\\1&0}F(z,1)^{-1}
\end{equation}
\begin{equation}
\label{eq-dfy}
\frac{\partial f}{\partial y}(z)=\frac{-\ii}{2}\;\mu(z) F(z,1)\matrix{0 &-\ii\\\ii&0}F(z,1)^{-1}
\end{equation}
where $\mu$ is the conformal factor of the immersion.
\medskip

We will use the following notation: if $u$ is a smooth function of $\lambda$ on the unit circle,
we denote $u^0$ the coefficient of $\lambda^0$ in its Fourier decomposition and
$u^{(-1)}$ the coefficient of $\lambda^{-1}$.
The elements of $\xi$ are denoted $\xi_{ij}$ for $1\leq i,j\leq 2$.
If the potential $\xi_t$ depends on some parameter $t$, its elements are
denoted $\xi_{t;ij}$, and the same conventions apply to all matrices.
With these notations, the conformal metric induced by the immersion is
given by
\begin{equation}
\label{eq-metric}
ds=\mu |dz|=2(B_{11}^0)^2|\xi_{12}^{(-1)}|.
\end{equation} 
\end{myenumerate}
\begin{remark}
What the DPW method actually does is construct the moving frame
$(\mu^{-1}\frac{\partial f}{\partial x},\mu^{-1}\frac{\partial f}{\partial y},N)$, which is
encoded in the unitary matrix $F(z,1)$.
The Sym-Bobenko formula is a magical trick to recover the immersion $f$ directly.
\end{remark}
\subsubsection{The Monodromy Problem}
Assume that $\Sigma$ is not simply connected so its universal cover $\wtSigma$
is not trivial.
Let $\mbox{Deck}(\wtSigma/\Sigma)$ be the group of fiber-preserving diffeomorphisms of $\wtSigma$.
For $\sigma\in\mbox{Deck}(\wtSigma/\Sigma)$, let
$$\boM(\Phi,\sigma)(\lambda)=\Phi(\sigma(z),\lambda)\Phi(z,\lambda)^{-1}$$ 
be the monodromy of $\Phi$ with respect to $\sigma$,
which is independent of $z\in\wtSigma$.
The standard condition which ensures that the immersion $f$ descends to
a well defined immersion on $\Sigma$
is the following system of equations, called the Monodromy Problem:
\begin{equation}
\label{pb-monodromy}
\forall \sigma\in\mbox{Deck}(\wtSigma/\Sigma)\quad
\left\{\begin{array}{lc}
\boM(\Phi,\sigma)\in\Lambda SU(2)\qquad &(i)\\
\boM(\Phi,\sigma)(1)=\pm I_2\qquad&(ii)\\
\frac{\partial}{\partial\lambda}\boM(\Phi,\sigma)(1)=0\qquad &(iii)
\end{array}\right.\end{equation}
\subsubsection{Gauging}
\begin{definition}
A gauge on $\Sigma$ is a map $G:\Sigma\to\Lambda_+ SL(2,\C)$ such that $G(z,\lambda)$ depends holomorphically on $z\in\Sigma$ and $\lambda\in\D_{\rho}$,
and $G(z,0)$ is upper triangular (with no restriction on its diagonal elements).
\end{definition}
Let $\Phi$ be a solution of $d\Phi=\Phi\xi$ and $G$ be a gauge.
Let $\whPhi=\Phi\times G$. It is standard that $\Phi$ and $\whPhi$
define the same immersion $f$.
The gauged potential is
$$\whxi:=\whPhi^{-1}d\whPhi=G^{-1}\xi G+G^{-1}d G$$
and is denoted $\xi\cdot G$, the dot denoting the action of the gauge
group on the potential.
Gauging does not change the Monodromy of $\Phi$ (provided $G$ is well defined in $\Sigma$): 
$$\boM(\Phi\times G,\sigma)=\Phi(\sigma(z))G(\sigma(z))G(z)^{-1}\Phi(z)^{-1}=\boM(\Phi,\sigma).$$
We will need the following standard result, which is an easy computation:
\begin{proposition}
\label{prop-gauging}
Let $G(z,\lambda)$ be a gauge and $\whxi=\xi\cdot G$.
Then (using the notations introduced at the end of Section \ref{section-background-DPW-recipe})
$$\whxi_{12}^{\,(-1)}=(G_{22}^0)^2\xi_{12}^{(-1)}$$
$$\whxi_{21}^{\,0}=(G_{11}^0)^2\xi_{21}^0.$$
\end{proposition}
\subsubsection{The Regularity Problem}
\begin{definition}
\label{def-regular}
We say that $\xi$ is regular at $z\in\Sigma$ if $\xi_{12}^{(-1)}(z)\neq 0$.
\end{definition}
In view of Equation \eqref{eq-metric}, this ensures that the immersion $f$ is unbranched at $z$.
In many cases, $\Sigma$ is a compact Riemann surface $\SSigma$
minus a finite number of points, and the potential $\xi$ extends meromorphically
to $\SSigma$.
\begin{definition}
\label{def-removable}
Let $p\in\SSigma$ be a pole of $\xi$. We say that $p$ is a removable singularity
if there exists a local gauge $G$ in a neighborhood of $p$ such that $\xi\cdot G$
extends holomorphically at $p$.
\end{definition}
This ensures that the immersion $f$ extends analytically at $p$. 
\medskip

In many cases, the meromorphic 1-form $\xi_{12}^{(-1)}$ has zeros in $\SSigma$.
(This is always the case if the genus of $\SSigma$ is greater than one). If $p$ is a zero of of $\xi_{12}^{(-1)}$ and we want $f$ to be unbranched at $p$, then $\xi$ must have a pole at $p$, $p$ must be a removable singularity and
$\xi\cdot G$ must be regular at $p$.
 
 \subsubsection{Dressing and isometries}
Let $H\in\Lambda SL(2,\C)$.
Let $\Phi(z,\lambda)$ be a solution of $d\Phi=\Phi\xi$.
Then $\wtPhi=H\Phi$ is a solution
of $d\wtPhi=\wtPhi\xi$. This is called ``dressing'' and amounts to change the
initial value $\phi_0$. In general, the effect of dressing on the immersion $f$
is not explicit. However if $H\in\Lambda SU(2)$ then by uniqueness in
the Iwasawa decomposition we have
$$\wtF(z,\lambda)=H(\lambda) F(z,\lambda)$$
 and by the Sym-Bobenko formula,
\begin{equation}
\label{eq-dressing}
\wtf(z)=H(1) f(z) H(1)^{-1}+\ii\frac{\partial H}{\partial\lambda}(1) H(1)^{-1}.
\end{equation}
Hence $\wtf=H\cdot f$ where
the action of $\Lambda SU(2)$ as rigid motions of $\su(2)\sim\R^3$ is given by
\begin{equation}
\label{eq-actionSU2}
H\cdot X=H(1)XH(1)^{-1}+\ii\frac{\partial H}{\partial\lambda}(1) H(1)^{-1}.
\end{equation}
\subsubsection{Duality}
\label{section-background-DPW-duality}
Let $\xi$ be a DPW potential, $\Phi$ a solution of Problem \eqref{eq-cauchy}
and
$$H(\lambda)=\matrix{0&\frac{\ii}{\sqrt{\lambda}}\\\ii\sqrt{\lambda}&0}.$$
I define the dual potential $\wcxi$ and its dual solution $\wcPhi$ (the terminology is not standard) by
$$\wcxi=H^{-1}\xi H\quad\mbox{ and }\quad
\wcPhi=H^{-1}\Phi H$$
which are both independent of the choice of the square root in $H$, so are well defined.
Explicitely,
$$H^{-1}\matrix{\alpha&\lambda^{-1}\beta\\\gamma &-\alpha}H=
\matrix{-\alpha&\lambda^{-1}\gamma\\\beta&\alpha}$$
so we see that duality essentially exchanges $\beta$ and $\gamma$.
We will take advantage of this to prove Corollary \ref{cor-regularity-infty} in 
Appendix \ref{appendix-regularity}.
Duality changes the immersion in the following explicit way:
it is easy to see that if $(F,B)$ is the Iwasawa decomposition of $\Phi$, then
the Iwasawa decomposition of $\wcPhi$ is
$$\wcF:=\Uni(\wcPhi)=H^{-1}FH\quad\mbox{ and }\quad
\wcB:=\Pos(\wcPhi)=H^{-1}B H.$$
If $f=\Sym(\Uni(f))$, the Sym-Bobenko formula gives after an easy computation:
$$\wcf(z):=\Sym(\Uni(\wcf))=H(1)^{-1}\left[\frac{-\ii}{2}\matrix{-1&0\\0&1}+f(z)+N(z)\right]H(1).$$
In other words, in euclidean coordinates,
$$\wcf(z)=\sigma(f(z)+N(z)+(0,0,1))\quad\mbox{ where }\quad
\sigma(x_1,x_2,x_3)=(x_1,-x_2,-x_3).$$
So up to a rigid motion, the dual (branched) immersion $\wcf$ is the parallel surface at distance one
to $f$.
Formula \eqref{eq-normal} gives the Gauss map of $\wcf$:
$$\wcN(z)=-\sigma(N(z)).$$
\begin{remark}
In this paper, we use the notations $\wtxi$, $\whxi$ and $\wcxi$ to denote various
transformations undergone by the potential $\xi$, including gauging, rescaling, pullback and dressing (so $\wtxi$ does not necessarily mean dressing and $\wcxi$ does not necessarily mean dual). These notations will always be consistently applied to all quantities derived from the potential: $\wtPhi$, $\wtF$, $\wtB$, $\wtf$ are derived from $\wtxi$.
\end{remark}
\subsubsection{The standard sphere}
\label{section-background-DPW-sphere}
We denote $\xiS$ the standard DPW potential for the sphere in $\C$:
$$\xiS(z,\lambda)=\matrix{0&\lambda^{-1}\\0&0}dz.$$
The following computation shows that $\infty$ is a regular removable singularity:
$$\xiS\cdot\matrix{z&0\\-\lambda&z^{-1}}=\matrix{0&\lambda^{-1}\\0&0}\frac{dz}{z^2}.$$
Let $\PhiS$ the solution of $d\PhiS=\PhiS\xiS$ with initial condition
$\PhiS(0,\lambda)=I_2$ in $\C$:
$$\PhiS(z,\lambda)=\matrix{1&\lambda^{-1}z\\0&1}.$$
The Iwasawa decomposition of $\PhiS(z)$ is explicitely given by
$$\FS(z,\lambda)=\Uni(\PhiS(z,\lambda))=
\frac{1}{\sqrt{1+|z|^2}}\matrix{1 & \lambda^{-1}z\\-\lambda\overline{z} & 1}$$
$$\BS(z,\lambda)=\Pos(\PhiS(z,\lambda))=
\frac{1}{\sqrt{1+|z|^2}}\matrix{1 & 0 \\ \lambda\overline{z} & 1+|z|^2}.$$

The Sym-Bobenko formula gives
$$\fS(z)=\frac{1}{1+|z|^2}\left(2\Re(z),2\Im(z),-2|z|^2\right)=(0,0,-1)+\pi_S^{-1}(z)$$
$$\NS(z)=\frac{-1}{1+|z|^2}\left(2\Re(z),2\Im(z),1-|z|^2\right)=-\pi_S^{-1}(z)$$
where $\pi_S:\S^2\to\C\cup\{\infty\}$ is the stereographic projection from the south pole.
\subsubsection{The infinitesimal catenoid}
\label{section-background-DPW-catenoid}
We denote $\xiC$ the dual potential to $\xiS$:
$$\xiC(z,\lambda)=\matrix{0&0\\1&0}dz$$
The dual solution is
$$\PhiC(z,\lambda)=\matrix{1&0\\z&1}$$
Its Iwasawa decomposition is the standard QR decomposition (from which we could derive the Iwasawa decomposition of $\PhiS$ by duality):
$$\FC(z,\lambda)=\Uni(\PhiC(z,\lambda))=\frac{1}{\sqrt{1+|z|^2}}\matrix{1 & -\overline{z}\\z &1}$$
$$B^C(z,\lambda)=\Pos(\PhiC(z,\lambda))=\frac{1}{\sqrt{1+|z|^2}}\matrix{1+|z|^2 & \overline{z}\\ 0 &1}$$
Since $\FC(z,\lambda)$ does not depend on $\lambda$, the Sym-Bobenko formula
gives $\fC\equiv 0$: the immersion degenerates into a point, as expected by duality since
the parallel surface at distance one to a unit sphere degenerates into its center.
The normal is still well defined and is given by:
$$\NC(z)=\frac{1}{1+|z|^2}\left(2\Re(z),-2\Im(z),|z|^2-1\right).$$
We will use this potential to model catenoidal necks.
(The catenoid is of course not a CMC-1 surface: it is a minimal surface.)
\subsection{Principal solution}
\label{section-background-principal}
We will formulate the Monodromy Problem using the notion of principal solution
(see Chapter 3.4 in \cite{teschl}).
Let $\Sigma$ be a Riemann surface, $\xi$ a matrix valued holomorphic 1-form on $\Sigma$
and $\gamma:[0,1]\to\Sigma$ a path (not necessarily closed).
The principal solution of $\xi$ with respect to $\gamma$, denoted $\boP(\xi,\gamma)$
is the value at $\gamma(1)$ of the analytical continuation along $\gamma$ of the solution $\Phi$ of
$d\Phi=\Phi\xi$ with initial condition $\Phi(\gamma(0))=I_2$.
More precisely, let $\wtSigma$ be the universal cover of $\Sigma$ and
$\wtgamma:[0,1]\to\wtSigma$ be an arbitrary lift of $\gamma$.
Let $\Phi$ be the solution on $\wtSigma$ of
$d\Phi=\Phi\xi$ with initial condition $\Phi(\wtgamma(0))=I_2$.
Then $\Phi(\wtgamma(1))$ does not depend on the choice of the lift $\wtgamma$
and is denoted $\boP(\xi,\gamma)$.
Equivalently, we can define $\boP(\xi,\gamma)$ as follows: let $Y$ be the solution on $[0,1]$ of the Cauchy Problem
$$\left\{\begin{array}{l}
Y'(s)=Y(s)\xi(\gamma(s))\gamma'(s)\\
Y(0)=I_2\end{array}\right.$$
Then $\boP(\xi,\gamma)=Y(1)$.
(The relation between the two definitions is $Y(s)=\Phi(\wtgamma(s))$.)
\medskip

The principal solution has the following properties, which follow easily from its definition:
\begin{myenumerate}
\item $\boP(\xi,\gamma)$ only depends on the homotopy class of
$\gamma$.
\item The principal solution is a morphism for the product of paths:
If $\gamma_1$ and $\gamma_2$ are two paths such that $\gamma_1(1)=\gamma_2(0)$
then
\begin{equation}
\label{eq-principal-morphism}
\boP(\xi,\gamma_1\gamma_2)=\boP(\xi,\gamma_1)\boP(\xi,\gamma_2).
\end{equation}
\item If $\psi:\Sigma_1\to\Sigma_2$ is a differentiable map, $\xi$ is a matrix-valued 1-form on $\Sigma_2$ and $\gamma$ is a path in $\Sigma_1$, then
\begin{equation}
\label{eq-principal-pullback}
\boP(\psi^*\xi,\gamma)=\boP(\xi,\psi\circ\gamma).
\end{equation}
\end{myenumerate}
We have the following formula to compute the derivative of the principal solution:
\begin{proposition}
\label{prop-principal-derivee}
Let $\xi_t$ be a family of holomorphic matrix-valued 1-forms on a Riemann surface
$\Sigma$, depending $C^1$ on some parameter $t$.
Let $\gamma$ be a path in $\Sigma$ and $\wtgamma$ an arbitrary lift of $\gamma$ to
the universal cover $\wtSigma$. Let $\Phi_t$ be the solution of
$d\Phi_t=\Phi_t\xi_t$ in $\wtSigma$ with initial condition $\Phi_t(\wtgamma(0))=I_2$.
Then
$$\frac{d}{dt}\boP(\xi_t,\gamma)=\int_{\wtgamma}\Phi_t\frac{\partial\xi_t}{\partial t}\Phi_t^{-1}\times\boP(\xi_t,\gamma).$$
\end{proposition}
This is proven in Appendix A of \cite{nnoids}.
(In \cite{nnoids}, this result is formulated in terms of monodromy so we assume that $\gamma$
is a closed curve, but this is not used in the proof.)
\begin{remark}
\label{remark-principal}
If $\xi$ is holomorphic in a simply connected domain $\Omega$ and $z_1$ and $z_2$ are two points in $\Omega$, then $\boP(\xi,\gamma)$ does not depend on the choice of the path $\gamma$ from $z_1$ to $z_2$ so we will sometimes denote it $\boP(\xi,z_1,z_2)$.
\end{remark}

Returning to the DPW method, we now formulate the Monodromy Problem in terms of the principal solution.
The group $\mbox{Deck}(\wtSigma/\Sigma)$ is isomorphic to the fundamental group
$\pi_1(\Sigma,z_0)$ (Theorem 5.6 in \cite{forster}):
for $\sigma\in\mbox{Deck}(\wtSigma/\Sigma)$, let $\wtgamma$ be an arbitrary path
in $\wtSigma$ from $\wtz_0$ to $\sigma(\wtz_0)$ and $\gamma\in\pi_1(\Sigma,z_0)$
be the projection of $\wtgamma$. Then $\gamma$ is the image of $\sigma$.
(This isomorphism is not canonical: it depends on
the choice of $\wtz_0$.)
The monodromy of $\Phi$ with respect to $\sigma$ and the principal solution of $\xi$ with respect to $\gamma$ are related by
\begin{equation}
\label{eq-boM-boP}
\boM(\Phi,\sigma)(\lambda)=\Phi(\wtgamma(1),\lambda)\Phi(\wtgamma(0),\lambda)^{-1}=\Phi(\wtz_0,\lambda)\boP(\xi,\gamma)(\lambda)\Phi(\wtz_0,\lambda)^{-1}.
\end{equation}
In particular, if $\Phi(\wtz_0,\lambda)=I_2$ (which will be the case in this paper), the Monodromy Problem \eqref{pb-monodromy} is
equivalent to the following problem (which we still call the Monodromy Problem):
\begin{equation}
\label{pb-principal}
\forall \gamma\in\pi_1(\Sigma,z_0)\quad
\left\{\begin{array}{lc}
\boP(\xi,\gamma)\in\Lambda SU(2)\qquad &(i)\\
\boP(\xi,\gamma)(1)=\pm I_2\qquad&(ii)\\
\frac{\partial}{\partial\lambda}\boP(\xi,\gamma)(1)=0\qquad &(iii)
\end{array}\right.\end{equation}
We conclude this section with a standard result which we will use to study the restriction
of the immersion $f$ to a subdomain $\Omega$ of $\Sigma$:
\begin{proposition}
\label{prop-restriction}
Let $p:\wtSigma\to\Sigma$ be the universal covering map of $\Sigma$.
Let $\Omega$ be a connected domain in $\Sigma$ and $z_0\in\Omega$. Let
$\wtz_0\in p^{-1}(z_0)$ and $\wtOmega$ be the
component of $p^{-1}(\Omega)$ containing $\wtz_0$.
Assume that the inclusion $i:\Omega\to\Sigma$ induces an injective
morphism $i_*:\pi_1(\Omega,z_0)\to\pi_1(\Sigma,z_0)$.
Then $\wtOmega$ is simply connected, so the restriction
$p:\wtOmega\to\Omega$ is a universal covering map of $\Omega$.
\end{proposition}
\subsection{Functional spaces}
\label{section-background-functional}
In this section, we introduce functional spaces for functions of the variable $\lambda\in\S^1$.
\subsubsection{The Banach algebra $\boW$}
\label{section-background-functional-boW}
We decompose a smooth function $f:\S^1\to\C$ in Fourier
series
$$f(\lambda)=\sum_{i\in\Z}c_i\lambda^i$$
Fix some $\rho>1$ and define
$$||f||=\sum_{i\in\Z}|c_i|\rho^{|i|}$$
Let $\boW$ be the space of functions $f$ with finite norm.
This is a Banach algebra (classically called the Wiener algebra when $\rho=1$).
Functions in $\boW$ extend holomorphically to the annulus $\A_{\rho}$ and
satisfy $|f(\lambda)|\leq ||f||$ for all $\lambda\in\A_{\rho}$.
\medskip

We define $\boW^{\pos}$, $\boW^{\pospos}$, $\boW^{\neg}$ and $\boW^{\negneg}$ as the subspaces of functions $f$ such that $c_i=0$
for $i<0$, $i\leq 0$, $i>0$ and $i\geq 0$, respectively.
Functions in $\boW^{\pos}$ extend holomorphically to the disk $\D_{\rho}$ and
satisfy $|f(\lambda)|\leq ||f||$ for all $\lambda\in\D_{\rho}$.
We also write $\boW^0$ for the subspace of constant functions, so we have
a direct sum $\boW=\boW^{\negneg}\oplus\boW^0\oplus\boW^{\pospos}$. A function $f$
will be decomposed as $f=f^{\negneg}+f^0+f^{\pospos}$, where $f^0=c_0$.
We also denote $f^{(-1)}=c_{-1}$ the coefficient of $\lambda^{-1}$ in the Fourier series of $f$.
\medskip

We define the star operator by
$$f^*(\lambda)=\overline{f\left(\frac{1}{\overline{\lambda}}\right)}=\sum_{i\in\Z}\overline{c_{-i}}\lambda^i$$
The involution $f\mapsto f^*$ exchanges $\boW^{\pos}$ and $\boW^{\neg}$.
We have $\lambda^*=\lambda^{-1}$ and $c^*=\overline{c}$ if $c$ is a constant.
\subsubsection{Linear operators}
The value $\lambda=1$ plays a special role in the DPW method because the Sym-Bobenko
formula is evaluated at $\lambda=1$. We shall need the following result in the case
$\mu=1$:
\begin{proposition}
\label{prop-boL}
For $\mu\in\C$, define
$$\boL_{\mu}(f)(\lambda)=\frac{f(\lambda)-f(\mu)}{\lambda-\mu}.$$
\begin{myenumerate}
\item If $\mu\in\D_{\rho}$, then $\boL_{\mu}:\boW^\pos\to\boW^\pos$ is a bounded operator with norm at most
$\frac{1}{\rho-|\mu|}$.
Consequently, any $f\in\boW^{\pos}$ can be decomposed in a unique way as
$$f(\lambda)=f(\mu)+(\lambda-\mu)\wtf(\lambda)\quad\mbox{ with } \wtf=\boL(f)\in\boW^\pos.$$
\item If $\mu\in\A_{\rho}$, then $\boL_{\mu}:\boW\to\boW$ is a bounded operator with norm at most $\max\{\frac{1}{\rho-|\mu|},\frac{1}{|\mu|-\rho^{-1}}\}$.
\end{myenumerate}
\end{proposition}
Proof: 
\begin{myenumerate}
\item Let $f\in\boW^\pos$.
Writing $f(\lambda)=\sum_{i=0}^\infty c_i \lambda^i$, we have:
$$\boL_{\mu}(f)(\lambda)=\sum_{i=0}^\infty c_i\frac{\lambda^i-\mu^i}{\lambda-\mu}=\sum_{i=0}^\infty
c_i\sum_{j=0}^{i-1}\lambda^j\mu^{i-1-j}$$
$$||\boL_{\mu}(f)||\leq \sum_{i=0}^\infty |c_i|\sum_{j=0}^{i-1}\rho^j|\mu|^{i-1-j}=\sum_{i=0}^{\infty}|c_i|\frac{\rho^i-|\mu|^i}{\rho-|\mu|}\leq \sum_{i=0}^{\infty}\frac{|c_i|\rho^i}{\rho-|\mu|}=\frac{||f||}{\rho-|\mu|}.$$
Hence $\boL_{\mu}(f)\in\boW^{\pos}$ and $|||\boL_{\mu}|||\leq \frac{1}{\rho-|\mu|}$.
\item Let $f\in\boW^\neg$. Then $f^*\in\boW^\pos$ and
$$\boL_{\mu}(f)^*(\lambda)
=\frac{f^*(\lambda)-\overline{f(\mu)}}{\lambda^{-1}-\overline{\mu}}
=-\frac{\lambda}{\overline{\mu}}\left(\frac{f^*(\lambda)-f^*(\frac{1}{\overline{\mu}})}{\lambda-\frac{1}{\overline{\mu}}}\right)
=-\frac{\lambda}{\overline{\mu}}\boL_{1/\overline{\mu}}(f^*)(\lambda).$$
Hence $\boL_{\mu}(f)\in\boW$ and since $*$ is an isometry of $\boW$ and using Point 1:
$$||\boL_{\mu}(f)||=||\boL_{\mu}(f)^*||
\leq\frac{\rho}{|\mu|}||\boL_{1/\overline{\mu}}(f^*)||
\leq \frac{\rho}{|\mu|} \frac{||f^*||}{(\rho-|\mu|^{-1})}=\frac{||f||}{|\mu|-\rho^{-1}}.$$
So the restriction of $\boL_{\mu}$ to $\boW^{\neg}$ has norm at most
$\frac{1}{|\mu|-\rho^{-1}}$.
\cqfd
\end{myenumerate}
\medskip

To prove that our operators are isomorphism, we will use the following elementary results:
\begin{definition}
\label{def-matrixtype}
Let $E$ be a Banach space. We say an operator $\boL=(\boL_1,\cdots,\boL_n):E^n\to E^n$ is
 of matrix type if there exists a matrix $A\in\boM_n(\C)$ such that
$$\forall (x_1,\cdots,x_n)\in E^n,\quad
\boL_i(x_1,\cdots,x_n)=\sum_{j=1}^n a_{ij} x_j$$
\end{definition}
Clearly, an operator of matrix type is invertible if its matrix $A$ is invertible.
\begin{proposition}
\label{prop-fredholm}
Let $E$ be a banach space and $F$, $G$ be finite dimensional normed vector space
of the same dimension.
Let $\boL=(\boL_1,\boL_2):E\times F\to E\times G$ be a bounded linear operator. Assume that
the restriction of $\boL_1$ to $E\times\{0\}$ is an automorphism of $E$ and $\boL$ is injective.
Then $\boL$ is an isomorphism.
\end{proposition}
Proof: Proposition \ref{prop-fredholm} can be proved by elementary means. Here is a short proof using the theory of Fredholm operators. Let $\boL':E\times F\to E\times G$ be the operator defined by
$\boL'(x,y)=(\boL_1(x,0),0)$.
Then $\mbox{Ker}(\boL')=\{0\}\times F$ and $\mbox{Im}(\boL')=E\times \{0\}$ so
$\boL'$ is a Fredholm operator of index $0$.
Now $\boL-\boL'$ is a finite rank operator (hence compact), so $\boL$ is a Fredholm operator of index $0$.
Hence $\boL$ injective implies that $\boL$ is an isomorphism.\cqfd
\begin{proposition}
\label{prop-perturbation-isomorphism}
Let $E,F$ be Banach spaces and $(\boL_n)_{n\in\N}$ be a sequence of bounded linear
operators from $E$ to $F$ converging to $\boL_{\infty}$. If $\boL_{\infty}:E\to F$ is
an isomorphism, then $\boL_n$ is an isomorphism for $n$ large enough.
\end{proposition}
Proof: for $n$ large enough, one has
$$|||\boL_{\infty}^{-1}\boL_n-id_{E}|||\leq |||\boL_{\infty}^{-1}|||\cdot|||\boL_n-\boL_{\infty}|||<1.$$
By the contraction mapping principle, $\boL_{\infty}^{-1}\boL_n$ is an automorphism of $E$ and Proposition \ref{prop-perturbation-isomorphism} follows.
\cqfd
\subsubsection{Holomorphic maps in Banach spaces}
\label{section-background-functional-holomorphy}
There is a theory of holomorphic maps between complex Banach spaces which
retains many features of the theory of functions of several complex variables.
A good reference is \cite{chae}.
Let $E$ and $F$ be Banach spaces. A map $f:\Omega\subset E\to F$ is analytic in
$\Omega$ if it admis a convergent ``power'' series expansion about any
point $x_0\in\Omega$ of the form
$$f(x)=\sum_{n=0}^{\infty} A_n (x-x_0)^n$$
where $A_n\in L(E^n,F)$ is a bounded, symmetric $n$-linear operator and $y^n$ denotes
the $n$-uple $(y,\cdots,y)$. We have the following fundamental results:
\begin{itemize}
\item If $f$ is analytic in $\Omega$ then $f$ is $C^{\infty}$ in $\Omega$ in the
Frechet sense (Theorem 11.12 in \cite{chae}).
\item If $E$ and $F$ are complex Banach spaces, then
$f$ is analytic in $\Omega$ if and only if $f$ is Frechet-differentiable in $\Omega$
(Graves-Taylor-Hille-Zorn, Theorem 14.6 in \cite{chae}). Such a map is called
holomorphic.
\item Hartog theorem on separate holomorphy remains true in the case of complex Banach spaces:
a map of a finite number of variables which is holomorphic with respect to each variable (the others being fixed) is holomorphic (Theorem 14.27 in \cite{chae}).
\end{itemize}
We will not need the Graves-Taylor-Hille-Zorn Theorem because the maps that we consider in this paper admit a power series expansion, where ``power'' refers to the product of the Banach algebra $\boW$.
For ${\bf a}=(a_1,\cdots,a_n)\in\C^n$ and ${\bf r}=(r_1,\cdots,r_n)\in(0,\infty)^n$,
we denote
$D({\bf a},{\bf r})$ the polydisk $\prod_{i=1}^n D(a_i,r_i)$ in $\C^n$.
\begin{proposition}[substitution]
\label{prop-substitution}
Let $R>\rho$ and $f:\A_R\times D({\bf a},{\bf r})\to\C$ be a holomorphic function of 
$(n+1)$ variables $(\lambda,z_1,\cdots,z_n)$.
Let
$$B({\bf a},{\bf r})=\{(u_1,\cdots,u_n)\in \boW^n:\forall i\in[1,n], ||u_i-a_i||<r_i\}$$
where we identify $a_i$ with a constant function in $\boW$.
Define for $(u_1,\cdots,u_n)\in B({\bf a},{\bf r})$:
$$F(u_1,\cdots,u_n)(\lambda)=f(\lambda,u_1(\lambda),\cdots,u_n(\lambda)).$$
Then
$F:B({\bf a},{\bf r})\subset\boW^n\to\boW$ is holomorphic.
\end{proposition}
This is proven in Appendix B of \cite{nnoids} by expanding $f$ in power series.
Thanks to Proposition \ref{prop-substitution}, to prove that a map is holomorphic,
we can first consider the case where the variables are complex numbers, and then
substitute functions of $\lambda$.
\section{Strategy}
\label{section-strategy}
\subsection{The case  $N=1$}
\label{section-strategy-nnoids}
In the case of one sphere, Theorem \ref{thm-main} is proved in \cite{nnoids},
using a DPW potential of the following form on the Riemann sphere:
$$\xi=\matrix{0 & \lambda^{-1}dz\\ t(\lambda-1)^2\omega(z,\lambda) & 0}\quad
\mbox{ with }\quad
\omega(z,\lambda)=\sum_{i=1}^n\left(\frac{a_i(\lambda)}{(z-p_i(\lambda))^2}+\frac{b_i(\lambda)}{z-p_i(\lambda)}\right)dz.$$
The parameters $a_i$, $b_i$ and $p_i$ are in the space $\boW^\pos$.
This potential is inspired from the potential used for 3-noids in \cite{schmitt-kilian-kobayashi-rossman}.
When $t=0$, we have $\xi=\xiS$, so this potential is a perturbation of the standard
spherical potential.
It is proven in \cite{nnoids} that the Monodromy Problem can be solved by an Implicit 
Function argument at $t=0$. Moreover, the potential can be locally gauged to a
potential with a simple pole at $p_i$ with a standard Delaunay residue, which implies
that the immersion has Delaunay ends by \cite{kilian-rossman-schmitt}.
We retain from this example how to grow Delaunay ends on the sphere by
perturbing the standard spherical potential.
\subsection{The construction in a nutshell}
\label{section-strategy-nutshell}
We want to construct a one-parameter family of compact Riemann surfaces $\SSigma_t$
and potentials $\xi_t$ for $t\in(0,\epsilon)$ by opening nodes on a Riemann surface
with nodes $\SSigma_0$. We want the principal solution of $\xi_t$ on paths which go
through a neck to extend continuously at $t=0$, so we request the regular
meromorphic potential $\xi_0$ to be holomorphic at the nodes. (If $\xi_0$ has a simple
pole at a node, the principal solution will diverge when reaching the node.)
We can ignore the rays of $\Gamma$ when defining $(\Sigma_0,\xi_0)$:
we have seen in Section \ref{section-strategy-nnoids} how to grow Delaunay ends by putting double poles
with principal parts of order $t$ in the potential.
\medskip

We want the immersion $f_t$ and the Gauss map $N_t$ to converge to well-defined
maps $f_0$ and $N_0$ on $\SSigma_0$: in particular they should be continuous
at the nodes.
$\SSigma_0$ should have one Riemann sphere for each $i\in[1,N]$, called a spherical part, on which $f_0$
will parametrize the sphere $\S^2(\v_i)$.
(The notation $\S^2(p)$ denotes the unit-sphere with center $p$.)
It should also have one Riemann sphere for each $(i,j)\in I$, called a catenoidal part, on which $f_0$
will be constant and equal to the tangency point between
$\S^2(\v_i)$ and $\S^2(\v_j)$.
On this Riemann sphere, $N_0$ should be the limit of
the gauss map $N_t$ of the catenoidal necks (a diffeomorphism to the sphere).
Observe that the normals of $\S^2(\v_i)$ and $\S^2(\v_j)$ at their tangency points are opposite, so the catenoidal part is required if we want $N_0$ to be continuous on
$\SSigma_0$.
\medskip

In the rest of this section, we explain how to construct $(\SSigma_0,\xi_0)$ so that
the Monodromy Problem for $\xi_0$ is solved.
Then in Section \ref{section-setup}, we open nodes to define $(\SSigma_t,\xi_t)$, throwing in a lot
of parameters which will be used to solve the Regularity and Monodromy Problems
for $t\neq 0$ in Sections \ref{section-regularity} to \ref{section-balancing}.
Finally, in Section \ref{section-geometry}, we prove that the resulting immersion $f_t$ has all the desired geometric properties.
\subsection{Notations}
\label{section-strategy-notations}
Without loss of generality, we assume (by a rotation) that all edges and rays of $\Gamma$ are non-vertical and (by a translation) that $\v_1=(0,0,-1)$.
\begin{itemize}
\item We define $I^*=\{(j,i):(i,j)\in I\}$. Two vertices $\v_i$ and $\v_j$ are adjacent if and only if $(i,j)\in I\cup I^*$.
\item We denote $\pi_S:\S^2\to\C\cup\{\infty\}$ the stereographic projection from the south pole.
\item For $(i,j)\in I\cup I^*$, we define $\pi_{ij}=\pi_S(\u_{ij})\in\C\setminus\{0\}$.
Since $\u_{ji}=-\u_{ij}$, we have $\pi_{ji}=\frac{-1}{\overline{\pi}_{ij}}$.
\item For $k\in[1,n]$, we define $\pi_k=\pi_S(\u_k)\in\C\setminus\{0\}$.
\end{itemize}

\subsection{The Riemann surface with nodes}
\label{section-strategy-nodes}
We define a compact Riemann surface with nodes $\SSigma_0$ as follows:
\begin{itemize}
\item
For each $i\in[1,N]$, consider a copy of the Rieman sphere $\CC=\C\cup\{\infty\}$,
denoted $\CC_i$ and called a spherical part.
\item 
For each $(i,j)\in I$, consider a copy of the Riemann sphere, denoted
$\CC_{ij}$ and called a catenoidal part.
\item
For each $(i,j)\in I$, identify a point $p_{ij}\in\C_i$ with a point $q_{ij}\in\C_{ij}$
and a point $p_{ji}\in\C_j$ with a point $q_{ji}\in\C_{ij}$ to create two nodes.
We will see later how to choose these points.
\end{itemize}
The points at infinity in $\CC_i$ and $\CC_{ij}$ are denoted respectively
$\infty_i$ and $\infty_{ij}$. The point $0$ in $\C_i$ and $\C_{ij}$ are denoted respectively $0_i$ and $0_{ij}$.
\subsection{The potential $\xi_0$: first guess}
We want to define a regular meromorphic potential $\xi_0$ on $\SSigma_0$ such
that the data $(\SSigma_0,\xi_0,0_1,I_2)$ yields by the DPW method a map $f_0$ and a Gauss map $N_0$, both continuous on $\SSigma_0$ and satisfying the following Ansatz:
\begin{equation}
\label{eq-ansatz}
\mbox{ For $i\in[1,N]$, $f_0$ is a translate of the standard spherical immersion $\fS$ in 
$\CC_i$.}
\end{equation}

The basic idea is to define $\xi_0$ on $\SSigma_0$ by
$$\xi_0=\left\{\begin{array}{l}
\xiS\;\mbox{on $\CC_i$ for $i\in[1,N]$}\\
\xiC\;\mbox{on $\CC_{ij}$ for $(i,j)\in I$.}\end{array}\right.$$
For $(i,j)\in I$, let $\Gamma_{ij}$ be a path from
$0_i$ to $0_j$ in in $\SSigma_0$, defined as the product of a path from $0_i$ to $p_{ij}$ in $\C_i$,
a path from $q_{ij}$ to $q_{ji}$ in $\C_{ij}$ and a path from $p_{ji}$ to $0_j$ in
$\C_j$. Also define $\Gamma_{ji}=\Gamma_{ij}^{-1}$.
Since $\xi_0$ is holomorphic in $\C_i$, $\C_{ij}$ and $\C_j$,
it is natural to define its principal solution with respect to $\Gamma_{ij}$ by
$$\boP(\xi_0,\Gamma_{ij})=\boP(\xi_0,0_i,p_{ij})\boP(\xi_0,q_{ij},q_{ji})\boP(\xi_0,q_{ji},0_j).$$
(In other words, if we consider the analytic continuation of the solution of $d\Phi_0=\Phi_0\xi_0$ along a path crossing a node, we simply require that $\Phi_0$
has the same value at the two points that are identified to create the node.
Theorem \ref{thm-principal} in Appendix \ref{appendix-principal} gives theoretical support for this definition.)
We would like to have the following Ansatz:
\begin{equation}
\label{eq-strategy-SU2}
\forall (i,j)\in I,\quad\boP(\xi_0,\Gamma_{ij})\in\Lambda SU(2).
\end{equation}
Since the fundamental group $\pi_1(\SSigma_0,0_1)$ is isomorphic to the
fundamental group of the graph $\Gamma$, any $\gamma\in\pi_1(\SSigma_0,0_1)$
can be written as a product
\begin{equation}
\label{eq-strategy-gamma}
\gamma=\prod_{j=1}^k \Gamma_{i_j i_{j+1}}
\end{equation}
with $i_1=i_{k+1}=1$. 
By Equation \eqref{eq-principal-morphism},
\begin{equation}
\label{eq-strategy-boPgamma}
\boP(\xi_0,\gamma)=\prod_{j=1}^k\boP(\xi_0,\Gamma_{i_j i_{j+1}}).
\end{equation}
So Ansatz \eqref{eq-strategy-SU2} implies that $\boP(\xi_0,\gamma)\in\Lambda SU(2)$
as required for the Monodromy Problem \eqref{pb-principal}.
\medskip

Unfortunately, \eqref{eq-strategy-SU2} does not hold for this choice of $\xi_0$: a computation gives:
\begin{eqnarray*}
\boP(\xi_0,\Gamma_{ij})&=&\PhiS(p_{ij},\lambda)\PhiC(q_{ij},\lambda)^{-1}\PhiC(q_{ji},\lambda)\PhiS(p_{ji},\lambda)^{-1}\\
&=&\matrix{1+\lambda^{-1}p_{ij}(q_{ji}-q_{ij})&\lambda^{-1}(p_{ij}-p_{ji})
+\lambda^{-2}p_{ij}p_{ji}(q_{ij}-q_{ji})\\q_{ji}-q_{ij}&1+\lambda^{-1}p_{ji}(q_{ij}-q_{ji})}
.\end{eqnarray*}
Whatever the choice of the points $p_{ij}$, $q_{ij}$, $p_{ji}$ and $q_{ji}$, this
matrix is not in $\Lambda SU(2)$ because there are only non-positive powers of $\lambda$.
So Ansatz \eqref{eq-strategy-SU2} does not hold. And worse, if $\gamma$ is a non trivial loop,
$\boP(\xi_0,\gamma)\not\in\Lambda SU(2)$ for the same reason
(unless really miraculous cancelations happen in the product \eqref{eq-strategy-boPgamma}, which does not seem to be the case). So the Monodromy Problem cannot be solved for this choice of
$\xi_0$ (unless $\Gamma$ is a tree so there is no Monodromy Problem for $\xi_0$).
\subsection{The potential $\xi_0$: second guess}
\label{section-strategy-potential}
The next idea (which works) is to gauge the catenoidal potential $\xiC$ in $\CC_{ij}$
by a well-chosen gauge $G_{ij}$:
we define $\xi_0$ by
\begin{equation}
\label{eq-xi0}
\xi_0=\left\{\begin{array}{l}
\xiS\quad\mbox{ in $\CC_i$ for $i\in[1,N]$}\\
\xiC\cdot G_{ij}\quad\mbox{in $\CC_{ij}$ for $(i,j)\in I$.}\end{array}\right.
\end{equation}
We now have
\begin{equation}
\label{eq-strategy-boPGammaij}
\boP(\xi_0,\Gamma_{ij})(\lambda)=\PhiS(p_{ij},\lambda)\left[\PhiC(q_{ij},\lambda)G_{ij}(q_{ij},\lambda)\right]^{-1}
\left[\PhiC(q_{ji},\lambda)G_{ij}(q_{ji},\lambda)\right]\PhiS(p_{ji},\lambda)^{-1}.
\end{equation}
In the next section, we choose the gauge $G_{ij}$ satisfying the following Ansatz:
\begin{equation}
\label{eq-strategy-doubleSU2}
\forall (i,j)\in I\cup I^*,\quad
\PhiS(p_{ij},\cdot)\left[\PhiC(q_{ij},\cdot)G_{ij}(q_{ij},\cdot)\right]^{-1}\in\Lambda SU(2)\end{equation}
By Equation \eqref{eq-strategy-boPGammaij}, Ansatz \eqref{eq-strategy-doubleSU2} implies Ansatz \eqref{eq-strategy-SU2}.
\subsection{Choosing the gauge $G_{ij}$}
\label{section-strategy-gauge}
Since we want $f_0$ to be a translate of $\fS$ on each $\CC_i$,
it is necessary to take
$$p_{ij}=\pi_{ij}\quad\mbox{ for $(i,j)\in I\cup I^*$}.$$
We propose to take
\begin{equation}
\label{eq-qij}
q_{ij}=\frac{1}{\overline{\pi}_{ij}}\quad\mbox{ for $(i,j)\in I\cup I^*$}.
\end{equation}
(Other choices are possible. This choice yields extra nice properties for the potential $\xi_0$.)
With these values, Ansatz \eqref{eq-strategy-doubleSU2} is equivalent to
$$\forall q\in\{q_{ij},q_{ji}\},\quad
\BS(\smallfrac{1}{\overline{q}},\cdot)G_{ij}(q,\cdot)^{-1}\BC(q,\cdot)^{-1}\in\Lambda SU(2).$$
Since this product is also in $\Lambda_+SL(2,\C)$, it must have
the form $\minimatrix{e^{\ii\theta_q}&0\\0&e^{-\ii\theta_q}}$ for some constant complex unitary number
$e^{\ii\theta_q}$.
So Ansatz \eqref{eq-strategy-doubleSU2} is equivalent to
\begin{eqnarray*}
\forall q\in\{q_{ij},q_{ji}\},\quad
G_{ij}(q,\lambda)&=&\BC(q,\lambda)^{-1}\matrix{e^{-\ii\theta_q}&0\\0&e^{\ii\theta_q}}\BS(\smallfrac{1}{\overline{q}},\lambda)\\
&=&\frac{1}{(1+|q|^2)|q|}\matrix{e^{-\ii\theta_q}|q|^2-\lambda\overline{q}^2e^{\ii\theta_q}&
-e^{\ii\theta_q}\overline{q}(1+|q|^2)\\ \lambda e^{\ii\theta_q}\overline{q}(1+|q|^2)& e^{\ii\theta_q}(1+|q|^2)^2}
\end{eqnarray*}
The simplest choice is to take $e^{\ii\theta_q}=\frac{q}{|q|}$, which gives the equation
\begin{equation}
\label{eq-Gijq}
\forall q\in\{q_{ij},q_{ji}\},\quad
G_{ij}(q,\lambda)=\matrix{\frac{(1-\lambda)\overline{q}}{1+|q|^2}&-1\\\lambda &\frac{1+|q|^2}{\overline{q}}}.\end{equation}
Let
\begin{equation}
\label{eq-muij}
\mu_{ij}=\frac{q_{ij}+q_{ji}}{2}.
\end{equation}
The following gauge does the job:
\begin{equation}
\label{eq-Gij}
G_{ij}(z,\lambda)=\matrix{\frac{(1-\lambda)}{2(z-\mu_{ij})}&-1\\\lambda&2(z-\mu_{ij})}.
\end{equation}
The gauged potential is simple enough: a computation gives
\begin{equation}
\label{eq-xiC.Gij0}
\xiC\cdot G_{ij}(z,\lambda)=
\matrix{\frac{\lambda-1}{2(z-\mu_{ij})}&1\\\frac{1-\lambda^2}{4(z-\mu_{ij})^2}&\frac{1-\lambda}{2(z-\mu_{ij})}}dz.
\end{equation}
In particular, at $\lambda=1$, it simplifies to $\minimatrix{0&1\\0&0}dz$, which will help
in solving Items (ii) and (iii) of the Monodromy Problem \eqref{pb-principal}.
This is a nice property that we get thanks to our choice of $q_{ij}$.
\subsection{Checking it works}
\label{section-strategy-immersion}
Assume that the graph $\Gamma$ has length-2 edges.
We now check that the ``immersion'' $f_0$ obtained with the potential $\xi_0$
and the initial data $(z_0,\phi_0)=(0_1,I_2)$ satisfies Ansatz \eqref{eq-ansatz}.
(Ansatz \eqref{eq-strategy-doubleSU2}
is satisfied by construction of $G_{ij}$.) The computations below will be used in
Section \ref{section-monodromy-edges} when solving the Monodromy Problem.
Let
$$U(q,\lambda)=\PhiS(\smallfrac{1}{\overline{q}},\lambda)[\PhiC(q,\lambda)G_{ij}(q,\lambda)]^{-1}.$$
Using Equation \eqref{eq-Gijq}, we obtain for $q\in\{q_{ij},q_{ji}\}$:
\begin{eqnarray*}
U(q,\lambda)&=&\matrix{1&\frac{1}{\lambda\overline{q}}\\0&1}
\matrix{\frac{1+|q|^2}{\overline{q}}&1\\-\lambda&\frac{(1-\lambda)\overline{q}}{1+|q|^2}}
\matrix{1&0\\-q&1}\\
&=&\frac{1}{1+|q|^2}\matrix{(1-\lambda^{-1})q&\lambda^{-1}+|q|^2\\-\lambda-|q|^2&
(1-\lambda)\overline{q}}\in\Lambda SU(2).
\end{eqnarray*}
By substitution of $q_{ij}=\frac{1}{\overline{\pi_{ij}}}$ and $q_{ji}=-\pi_{ij}$, we obtain
for $(i,j)\in I$:
\begin{equation}
\label{eq-Uqij}
U(q_{ij},\lambda)=\frac{1}{1+|\pi_{ij}|^2}\matrix{
(1-\lambda^{-1})\pi_{ij}&1+\lambda^{-1}|\pi_{ij}|^2\\
-1-\lambda|\pi_{ij}|^2&(1-\lambda)\overline{\pi}_{ij}}
\end{equation}
$$U(q_{ji},\lambda)=\frac{1}{1+|\pi_{ij}|^2}\matrix{
(\lambda^{-1}-1)\pi_{ij}&\lambda^{-1}+|\pi_{ij}|^2\\
-\lambda-|\pi_{ij}|^2&(\lambda-1)\overline{\pi}_{ij}}.$$
By Equation \eqref{eq-strategy-boPGammaij}, we obtain:
\begin{equation}
\label{eq-strategy-boPGammaijfinal}
\boP(\xi_0,\Gamma_{ij})(\lambda)=U(q_{ij},\lambda)U(q_{ji},\lambda)^{-1}
=\frac{1}{1+|\pi_{ij}|^2}\matrix{
\lambda^{-1}|\pi_{ij}|^2+\lambda&
\pi_{ij}(\lambda^{-2}-1)\\
\overline{\pi}_{ij}(1-\lambda^2)&\lambda|\pi_{ij}|^2+\lambda^{-1}}.
\end{equation}
This implies
\begin{equation}
\label{eq-strategy-boPGammaij1}
\boP(\xi_0,\Gamma_{ij})(1)=I_2
\end{equation}
\begin{equation}
\label{eq-strategy-boPGammaij2}
\ii\frac{\partial}{\partial\lambda}\boP(\xi_0,\Gamma_{ij})(1)=
\frac{\ii}{1+|\pi_{ij}|^2}\matrix{1-|\pi_{ij}|^2&-2\pi_{ij}\\
-2\overline{\pi}_{ij}&|\pi_{ij}|^2-1}=2\pi_S^{-1}(\pi_{ij})=2\u_{ij}=\v_j-\v_i.
\end{equation}
(In the last equality, we have used the fact that the edges have length $2$.)
If we decompose an arbitrary element $\gamma$ of $\pi_1(\SSigma_0,0_1)$
as in Equation \eqref{eq-strategy-gamma}, we obtain from Equations 
\eqref{eq-strategy-boPgamma}, \eqref{eq-strategy-boPGammaij1} and \eqref{eq-strategy-boPGammaij2}:
$$\boP(\xi_0,\gamma)(1)=I_2$$
$$\ii\frac{\partial}{\partial\lambda}\boP(\xi_0,\gamma)(1)=
\ii\sum_{j=1}^k\frac{\partial}{\partial\lambda}\boP(\xi_0,\Gamma_{i_ji_{j+1}})(1)
=\sum_{j=1}^k (\v_{i_{j+1}}-\v_{i_j})=0.$$
Hence the Monodromy Problem \eqref{pb-principal} is solved, so the DPW method
produces a well defined map $f_0$ on $\SSigma_0$. (We do not call it an immersion
because it is constant in the catenoidal parts.)
By our choice of the initial condition, we have $f_0=\fS$ in $\CC_1$, and by
Equation \eqref{eq-strategy-boPGammaij2}:
$$f_0=\fS+\v_i-\v_1\quad\mbox{ in $\CC_i$}$$
so Ansatz \eqref{eq-ansatz} is satisfied.
\section{Setup}
\label{section-setup}
In this section, we define a family of compact Riemann surfaces $\SSigma_{t,\x}$ and meromorphic DPW potentials
$\xi_{t,\x}$, depending on the complex parameter $t$ in a neighborhood of $0$ and many other parameters that we will introduce.
The vector of these parameters (put in some arbitrary order) is denoted $\x$.
The parameters involved in the definition of the Riemann surface are complex numbers,
while the parameters involved in the definition of the DPW potential are functions of $\lambda$
in the space $\boW^{\pos}$.
The parameter vector $\x$ is in a neighborhood of
a central value denoted $\x_0$. We will solve equations using the Implicit Function Theorem at the point $(t,\x)=(0,\x_0)$.
When solving the Monodromy Problem, we will restrict the parameter $t$ to positive real values.
\subsection{Opening nodes}
\label{section-setup-openingnodes}
For $(i,j)\in I\cup I^*$, we introduce a complex parameter $p_{ij}$ in a neighborhood of $\pi_{ij}$ and a non-zero complex parameter
$r_{ij}$. The point $q_{ij}$ is fixed and given by Equation \eqref{eq-qij}.
We define a compact Riemann surface with nodes as explained in Section \ref{section-strategy-nodes}, except that now the points $p_{ij}$ are parameters so we denote it $\SSigma_{0,\x}$ instead of
$\SSigma_0$.

\medskip
We fix a positive $\varepsilon<1$ small enough such that for all $i\in[1,N]$,
the disks $D(\pi_{ij},2\varepsilon)$ for $j\in E_i$ and the disks $D(\pi_k,2\varepsilon)$
for $k\in R_i$  are disjoint and do not contain $0$.
We assume that $|p_{ij}-\pi_{ij}|<\varepsilon$ for all $(i,j)\in I\cup I^*$,
so for all $i\in[1,N]$, the disks $D(p_{ij},\varepsilon)$ for $j\in E_i$ are disjoint.
The disks $D(q_{ij},\varepsilon)$ and $D(q_{ji},\varepsilon)$ are disjoint because
$q_{ji}=-\frac{1}{\overline{q_{ij}}}$ and $\varepsilon<1$.
\medskip

To open nodes, we introduce, for each $(i,j)\in I$, the following standard local complex coordinates
\begin{itemize}
\item $v_{ij}=z-p_{ij}$ in the disk $D(p_{ij},\varepsilon)$ in $\CC_i$,
\item $w_{ij}=z-q_{ij}$ in the disk $D(q_{ij},\varepsilon)$ in $\CC_{ij}$,
\item $v_{ji}=z-p_{ji}$ in the disk $D(p_{ji},\varepsilon)$ in $\CC_j$,
\item  $w_{ji}=z-q_{ji}$ in the disk $D(q_{ji},\varepsilon)$ in $\CC_{ij}$.
\end{itemize}
For $(i,j)\in I\cup I^*$, we take
$t_{ij}=t\,r_{ij}.$
As explained in Section \ref{section-background-openingnodes}, we remove the
disks $|v_{ij}|\leq\frac{|t_{ij}|}{\varepsilon}$ and $|w_{ij}|\leq \frac{|t_{ij}|}{\varepsilon}$.
We Identify
each point $z$ in the annulus $\frac{|t_{ij}|}{\varepsilon}<|v_{ij}|<\varepsilon$
with the point $z'$ in the annulus $\frac{|t_{ij}|}{\varepsilon}<|w_{ij}|<\varepsilon$
such that
$$v_{ij}(z)w_{ij}(z')=t_{ij}=t\, r_{ij}.$$
This defines for $t\neq 0$ a genuine compact Riemann surface which we denote $\SSigma_{t,\x}$.
\begin{remark}
All nodes open at the same time and the parameter $r_{ij}$ controls the speed at which
the node $p_{ij}\sim q_{ij}$ opens. The parameter $r_{ij}$ is related to the weight
$\tau_{ij}$.
\end{remark}
 \subsection{Definition of the DPW potential}
 \label{section-setup-potential}
We first perturb the gauge $G_{ij}$ introduced in Section \ref{section-strategy-gauge}.
For $(i,j)\in I$, we introduce three parameters
$g_{ij},h_{ij},m_{ij}$ in $\boW^\pos$ with respective central values
$1$, $1$ and $\mu_{ij}$.
We define a gauge $G_{ij,\x}$ in $\CC_{ij}$ by
 $$G_{ij,\x}(z,\lambda)=\matrix{\frac{(1-\lambda)}{2g_{ij}(\lambda)(z-m_{ij}(\lambda))}&-h_{ij}(\lambda)\\
\frac{\lambda}{h_{ij}(\lambda)}&2g_{ij}(\lambda)(z-m_{ij}(\lambda))}.$$
At the central value, we have $G_{ij,\x_0}=G_{ij}$.
We define a meromorphic regular potential $\zeta_{0,\x}$ on $\SSigma_{0,\x}$
as follows:
$$\zeta_{0,\x}=\left\{\begin{array}{ll}
\xiS&\mbox{ in $\CC_i$ for $i\in[1,N]$}\\
\xiC\cdot  G_{ij,\x}&\mbox{ in $\CC_{ij}$ for $(i,j)\in I$}.
\end{array}\right.$$
Explicitly, a computation gives: 
\begin{equation}
\label{eq-xiC.Gij}
\xiC\cdot G_{ij,\x}(z,\lambda)=\matrix{\frac{(\lambda-1)(2g_{ij}-h_{ij})}{2g_{ij}(z-m_{ij})}
&h_{ij}(2g_{ij}-h_{ij})\\\frac{(\lambda-1)(\lambda h_{ij}-h_{ij}-2\lambda g_{ij})}{4 h_{ij}g_{ij}^2(z-m_{ij})^2}&\frac{(1-\lambda)(2g_{ij}-h_{ij})}{2g_{ij}(z-m_{ij})}}dz.
\end{equation}
At the central value, we have $\zeta_{0,\x_0}=\xi_0$.
For $t$ in a neighborhood of $0$, the meromorphic regular potential $\zeta_{0,\x}$ defines a meromorphic potential $\zeta_{t,\x}$ on $\SSigma_{t,\x}$.
(We apply Theorem \ref{thm-openingnodes} to each element of $\zeta_{0,\x}$.)
We are not done yet: we still have to introduce poles for the Delaunay ends and we can prescribe the periods of our potential around the nodes.
We introduce the following parameters:
\begin{itemize}
\item $a_{ij},b_{ij},c_{ij}$ for $(i,j)\in I\cup I^*$,
\item $a_{mij},b_{mij},c_{mij}$ for $(i,j)\in I$,
\item $a_k,b_k,p_k$ for $k\in [1,n]$.
\end{itemize}
All these parameters are functions of $\lambda$ in the space $\boW^\pos$.
We give their central value in Section \ref{section-setup-parameters}.
We define the meromorphic regular potentials $\chi_{0,\x}$ and $\Theta_{0,\x}$ on
$\SSigma_{0,\x}$ by:
\begin{itemize}
\item For $i\in[1,N]$ and $z\in\C_i$:
\begin{equation}
\label{eq-chi0Ci}
\chi_{0,\x}(z,\lambda)=
(\lambda-1)\sum_{k\in R_i}\matrix{0&0\\\frac{a_k(\lambda)}{(z-p_k)^2}+\frac{b_k(\lambda)}{(z-p_k)}&0}dz
+\sum_{j\in E_i}\matrix{a_{ij}(\lambda)&\lambda^{-1}b_{ij}(\lambda)\\c_{ij}(\lambda)&-a_{ij}(\lambda)}\frac{dz}{z-p_{ij}}
\end{equation}
$$\Theta_{0,\x}(z,\lambda)=0.$$
\item For $(i,j)\in I$ and $z\in\C_{ij}$:
\begin{eqnarray}
\label{eq-chi0Cij}
\chi_{0,\x}(z,\lambda)&=&
-\matrix{a_{ij}(\lambda)&\lambda^{-1}b_{ij}(\lambda)\\c_{ij}(\lambda)&-a_{ij}(\lambda)}\frac{dz}{z-q_{ij}}
-\matrix{a_{ji}(\lambda)&\lambda^{-1}b_{ji}(\lambda)\\c_{ji}(\lambda)&-a_{ji}(\lambda)}\frac{dz}{z-q_{ji}}\\
\nonumber
&&+\matrix{0&0\\\frac{c_{ij}(\lambda)+c_{ji}(\lambda)}{z-m_{ij}(\lambda)}&0}dz
\end{eqnarray}
\begin{equation}
\label{eq-Theta0}
\Theta_{0,\x}=G_{ij,\x}(z,\lambda)^{-1}
\matrix{\frac{a_{mij}(\lambda)}{z-m_{ij}}&\frac{\lambda^{-1}b_{mij}(\lambda)}{(z-m_{ij})^2}+\frac{\lambda^{-1}c_{mij}(\lambda)}{z-m_{ij}}\\0&\frac{-a_{mij}(\lambda)}{z-m_{ij}}}G_{ij,\x}(z,\lambda)dz.
\end{equation}
\end{itemize}
For $t$ in a neighborhood of $0$, the meromorphic regular potentials $\chi_{0,\x}$ and $\Theta_{0,\x}$ define two meromorphic potentials $\chi_{t,\x}$ and $\Theta_{t,\x}$ on $\SSigma_{t,\x}$
(using Theorem \ref{thm-openingnodes} again).
We define the meromorphic potential $\xi_{t,\x}$ on $\SSigma_{t,\x}$ by
\begin{equation}
\label{eq-xit}
\xi_{t,\x}=\zeta_{t,\x}+ t(\lambda-1)(\chi_{t,\x}+\Theta_{t,\x}).
\end{equation}
\begin{remark}  The potential $\xi_{t,\x}$ looks awfully complicated, so
let me explain the purpose of each term in its definition.
Each parameter is introduced to solve a certain problem:
\begin{myenumerate}
\item The gauge $G_{ij}$ was chosen in Section \ref{section-strategy-gauge} so that
the principal solution of $\xi_0$ along $\Gamma_{ij}$ is in $\Lambda SU(2)$.
The parameters $g_{ij}$, $h_{ij}$ and
$m_{ij}$ involved in the definition of the gauge $G_{ij,\x}$ will be used in Section \ref{section-monodromy-edges}
to solve that problem when $t\neq 0$ using the Implicit Function Theorem.
\item The first term in \eqref{eq-chi0Ci} creates the desired Delaunay ends, as explained in Section \ref{section-strategy-nnoids}.
\item The second term in \eqref{eq-chi0Ci} prescribes the period of $\chi_{t,\x}$
around the nodes. The parameters
$a_{ij}$, $b_{ij}$ and $c_{ij}$ will be used in Section \ref{section-monodromy-nodes} to solve the Monodromy Problem with respect to cycles around the nodes.
\item The first term in \eqref{eq-chi0Cij} is forced by the fact that the residues of $\chi_{0,\x}$
at $p_{ij}$ and $q_{ij}$ must be opposite (see Definition \ref{def-regulardifferential}). Same for the second term.
\item The third term in \eqref{eq-chi0Cij} is there so that $\chi_{t,\x;21}$ is holomorphic at
$\infty_{ij}$, which will be useful when solving the Regularity Problem at $\infty_{ij}$.
\item The parameters $a_{mij}$, $b_{mij}$, $c_{mij}$ in the definition of $\Theta_{0,\x}$ will be used in Section \ref{section-regularity} to solve the Regularity Problem at $m_{ij}$, namely ensure that $\xi_{t,\x}\cdot G_{ij,\x}^{-1}$
is holomorphic at $m_{ij}$. There is no need to compute explicitely the matrix product in the
definition of $\Theta_{0,\x}$: the matrices $G_{ij,\x}^{-1}$ and $G_{ij,\x}$ will cancel when we gauge by
$G_{ij,\x}^{-1}$. The only thing that matters is that $\Theta_{0,\x}$ has poles only
at $m_{ij}$ and $\infty_{ij}$, which is clear from the definition of $G_{ij,\x}$.
\end{myenumerate}
\end{remark}
We collect some immediate properties of the potential $\xi_{t,\x}$ in
the following
\begin{proposition}
\label{prop-properties}
\begin{myenumerate}
\item At the central value, $\xi_{0,\x_0}=\xi_0$, where $\xi_0$ is given by
Equation \eqref{eq-xi0}.
\item If $t\neq 0$, $\xi_{t,\x}$ has poles at the following points:
\begin{itemize}
\item $\infty_i$ in $\CC_i$ for $i\in [1,N]$,
\item $m_{ij}$ and $\infty_{ij}$ in $\CC_{ij}$ for $(i,j)\in I$,
\item $p_k$ in $\C_i$ for $i\in[1,N]$ and $k\in R_i$.
\end{itemize}
\item At $\lambda=1$, $\xi_{t,\x}$ has the following form:
\begin{equation}
\label{eq-xit1}
\xi_{t,\x}(z,1)=\matrix{0&\beta_{t,\x}(z)\\0&0}\end{equation}
where $\beta_{t,\x}$ is a meromorphic 1-form on $\SSigma_{t,\x}$ with no periods
around the nodes.
\end{myenumerate}
\end{proposition}
Proof:
the only point which does not follow directly from the definitions is Point 3.
At $\lambda=1$, we have $\xi_{t,\x}(z,1)=\zeta_{t,\x}(z,1)$.
By definition of $\zeta_{0,\x}$ and Equation \eqref{eq-xiC.Gij}, we have
$$\zeta_{0,\x}(z,1)=\matrix{0&\beta_{0,\x}(z)\\0&0}\quad\mbox{ with }\quad
\beta_{0,\x}=\left\{\begin{array}{cc}
dz&\mbox{ in $\CC_i$}\\
h_{ij}(1)(2g_{ij}(1)-h_{ij}(1))dz&\mbox{ in $\CC_{ij}$}\end{array}\right.$$
Observe that $\beta_{0,\x}$ is holomorphic at the nodes.
By Theorem \ref{thm-openingnodes}, the regular differential $\beta_{0,\x}$
defines a unique meromorphic 1-form $\beta_{t,\x}$ on $\SSigma_{t,\x}$
with no periods around the nodes.\cqfd
\subsection{Parameters of the construction}
\label{section-setup-parameters}
The following table gives all parameters of the construction, together with their index range,
functional space, central value, and the section in which each parameter is
used to solve an equation. The parameters appear in the order in which they are used.
\medskip

\begin{center}
\begin{tabular}{|c|c|c|c|c|}\hline
{\scriptsize Parameter} &{\scriptsize Index} & {\scriptsize Space} & {\scriptsize Central value} & {\scriptsize Section}\\\hline
$a_{mij}$& $I$ &$\boW^\pos$ & $-4\lambda\tau_{ij}$ & \ref{section-regularity}\\\hline
$b_{mij}$& $I$ &$\boW^\pos$ & $-2(\lambda-1)\lambda\tau_{ij}$ & \ref{section-regularity}\\\hline
$c_{mij}$& $I$ & $\boW^\pos$ & $0$ & \ref{section-regularity}\\\hline
$a_{ij}$& $I\cup I^*$&$\boW^\pos$ & $2\tau_{ij}$ & \ref{section-monodromy-nodes}\\\hline
$b_{ij}$& $I\cup I^*$&$\boW^\pos$ & $2\rho_{ij}^{-1}\tau_{ij}$ & \ref{section-monodromy-nodes}\\\hline
$c_{ij}$& $I\cup I^*$&$\boW^\pos$ & $-2(\lambda-1)\rho_{ij}\tau_{ij}$ & \ref{section-monodromy-nodes}\\\hline
$r_{ij}$&$I\cup I^*$& $\C$ & $-\rho_{ij}^{-2}\tau_{ij}$ & \ref{section-monodromy-nodes}\\\hline
\end{tabular}
\begin{tabular}{|c|c|c|c|c|}\hline
{\scriptsize Parameter}&{\scriptsize Index} & {\scriptsize Space} & {\scriptsize Central value} & {\scriptsize Section}\\\hline
$a_k$ &$[1,n]$&$\boW^\pos$ & $\tau_k$&\ref{section-monodromy-ends}\\\hline
$b_k$ &$[1,n]$&$\boW^\pos$ & $-2\rho_k\tau_k$&\ref{section-monodromy-ends}\\\hline
$p_k$ &$[1,n]$&$\boW^\pos$ & $\pi_k$ &\ref{section-monodromy-ends}\\\hline
$g_{ij}$ &$I$&$\boW^\pos$&$1$&\ref{section-monodromy-edges}\\\hline
$h_{ij}$ &$I$&$\boW^\pos$&$1$&\ref{section-monodromy-edges}\\\hline
$m_{ij}$ &$I$&$\boW^\pos$&$\mu_{ij}$&\ref{section-monodromy-edges}\\\hline
$p_{ij}$&$I\cup I^*$&$\C$&$\pi_{ij}$&\ref{section-monodromy-edges}\\\hline
\end{tabular}\end{center}
\medskip

In the fourth column, the greek letters refer to constants depending only
on the weighted graph $\Gamma$:
\begin{itemize}
\item For $(i,j)\in I$, $\tau_{ij}$ is the weight of the edge $e_{ij}$ and $\tau_{ji}=\tau_{ij}$.
\item For $(i,j)\in I\cup I^*$, $\pi_{ij}=\pi_S(\u_{ij})$. We have $\pi_{ji}=\frac{-1}{\overline{\pi_{ij}}}$.
\item For $(i,j)\in I\cup I^*$, $\rho_{ij}=\displaystyle\frac{\overline{\pi_{ij}}}{1+|\pi_{ij}|^2}$.
We have $\rho_{ji}=-\rho_{ij}$.
\item For $(i,j)\in I$, $\mu_{ij}=\mu_{ji}=\displaystyle\frac{q_{ij}+q_{ji}}{2}=\frac{1-|\pi_{ij}|^2}{2\overline{\pi_{ij}}}$.
\item For $k\in[1,n]$: $\tau_k$ is the weight of the ray $r_k$, $\pi_k=\pi_S(\u_k)$ and $\rho_k=\displaystyle\frac{\overline{\pi_k}}{1+|\pi_k|^2}$.
\end{itemize}
Until Section \ref{section-balancing}, we do not assume that the graph $\Gamma$
is balanced nor has length-2 edges.
\section{The Regularity Problem at $m_{ij}$}
\label{section-regularity}
We want the following poles of $\xi_{t,\x}$ to be removable singularities
(see Definition \ref{def-removable}):
\begin{itemize}
\item the points $\infty_{i}$ for $i\in[1,N]$,
\item the points $\infty_{ij}$ for $(i,j)\in I$,
\item the points $m_{ij}$ for $(i,j)\in I$.
\end{itemize}
We call this the Regularity Problem.
In this section, we solve the Regularity Problem at $m_{ij}$.
The Regularity Problems at $\infty_{ij}$ and $\infty_i$ are solved in Sections
\ref{section-regularity2} and \ref{section-balancing}.
We define in $\C_{ij}$:
\begin{equation}
\label{eq-whxiij}
\whxi_{ij,t,\x}=\xi_{t,\x}\cdot G_{ij,\x}^{-1}.
\end{equation}
Our goal is to adjust the parameters $a_{mij}$, $b_{mij}$ and $c_{mij}$ so that $\whxi_{ij,t,\x}$
extends holomorphically at $m_{ij}$.
We will see in Section \ref{section-geometry-immersion} that $\whxi_{ij,t,\x}$ is regular at $m_{ij}$.
\subsection{Order of $\whxi_{ij,t,\x}$ at $m_{ij}$ and $\infty_{ij}$}
\label{section-regularity-order}
The following terminology will be convenient.
Let $f$ be a meromorphic function or 1-form on a Riemann surface $\Sigma$.
Let $p\in\Sigma$ and $z$ be a local coordinate in a neighborhood of $p$
such that $z(p)=0$. The order of $f$ at $p$, denoted $\Ord_p(f)$, is the largest
$a\in\Z\cup\{\infty\}$ such that $z^{-a}f$ is holomorphic at $p$.
(So $a>0$ means that $f$ has a zero of multiplicity $a$ at $p$,
$a<0$ means that $f$ has a pole of multiplicity $-a$ at $p$
and $a=\infty$ means $f\equiv 0$.)
If $F$ is a $2\times 2$ matrix of meromorphic functions or 1-forms, we define
$$\Ord_p(F)=(\Ord_p(F_{ij}))_{1\leq i,j\leq 2}\in\boM_2(\Z\cup\{\infty\}).$$
It is straightforward to check that
$$\Ord_p(F+G)\geq\min(\Ord_p(F),\Ord_p(G))$$
$$\Ord_p(F\times G)\geq \Ord_p(F)\star\Ord_p(G)$$
where $\star$ is the ``tropical'' matrix product obtained by replacing
$(+,\times)$ by $(\min,+)$ in the usual matrix product,
$\min(A,B)=(\min(A_{ij},B_{ij}))_{1\leq i,j\leq 2}$ and
$A\geq B$ means that $A_{ij}\geq B_{ij}$ for all $ i,j\in[1,2]$.
\begin{proposition}
\label{prop-order}
$\whxi_{ij,t,\x}$ has order at least $\minimatrix{-1&-2\\0&-1}$ at $m_{ij}$ and
at least $\minimatrix{-1&-1\\-2&-1}$ at $\infty_{ij}$.
\end{proposition}
(In other words, $\whxi_{ij,t,\x}$ has poles of multiplicity at most
$\minimatrix{1&2\\0&1}$ at $m_{ij}$ and at most $\minimatrix{1&1\\2&1}$ at $\infty_{ij}$).
\medskip

Proof: we compute the order at $m_{ij}$ of each term in the definition of $\xi_{t,\x}$.
Recall that by definition (see Theorem \ref{thm-openingnodes}), $\zeta_{t,\x}-\zeta_{0,\x}$ is
holomorphic at $m_{ij}$, so by Equation \eqref{eq-xiC.Gij}:
$$\Ord_{m_{ij}}(\zeta_{t,\x})\geq\min\left[\Ord_{m_{ij}}(\zeta_{0,\x}),\minimatrix{0&0\\0&0}\right]
\geq\minimatrix{-1&0\\-2&-1}.$$
In the same way, by Equation \eqref{eq-chi0Cij}:
$$\Ord_{m_{ij}}(\chi_{t,\x})\geq\minimatrix{0&0\\-1&0}$$
By Definition of $G_{ij,\x}$,
$$\Ord_{m_{ij}}(G_{ij,\x})\geq \minimatrix{-1&0\\0&1}.$$
By Equation \eqref{eq-Theta0}:
$$\Ord_{m_{ij}}(\Theta_{0,\x})\geq \minimatrix{1&0\\0&-1}\star\minimatrix{-1&-2\\\infty&-1}\star
\minimatrix{-1&0\\0&1}=\minimatrix{-1&0\\-2&-1}.$$
Hence
$$\Ord_{m_{ij}}(\xi_{t,\x})\geq \minimatrix{-1&0\\-2&-1}.$$
$$\Ord_{m_{ij}}(G_{ij,\x}\xi_{t,\x}G_{ij,\x}^{-1})\geq\minimatrix{-1&0\\0&1}\star
\minimatrix{-1&0\\-2&-1}\star\minimatrix{1&0\\0&-1}=\minimatrix{-1&-2\\0&-1}.$$
$$\Ord_{m_{ij}}(G_{ij,\x}dG_{ij,\x}^{-1})\geq
\minimatrix{-1&0\\0&1}\star\minimatrix{0&\infty\\\infty&-2}
=\minimatrix{-1&-2\\0&-1}.$$
So
$$\Ord_{m_{ij}}(\whxi_{ij,t,\x})\geq \minimatrix{-1&-2\\0&-1}.$$
In the same way:
$$\Ord_{\infty_{ij}}(\zeta_{t,\x})\geq\minimatrix{-1&-2\\0&-1}.$$
Thanks to the third term in Equation \eqref{eq-chi0Cij}:
$$\Ord_{\infty_{ij}}(\chi_{t,\x})\geq\minimatrix{-1&-1\\0&-1}$$
$$\Ord_{\infty_{ij}}(G_{ij,\x})\geq\minimatrix{1&0\\0&-1}.$$
The term $\Theta_{t,\x}$ is more delicate: to get the right order, it is necessary to write
$$\Theta_{t,\x}=\Theta_{0,\x}+\Xi_{t,\x}$$
where $\Xi_{t,\x}$ is holomorphic in a neighborhood of $\infty_{ij}$.
Then by definition of $\Theta_{0,\x}$:
\begin{equation}
\label{whThetaij0}
G_{ij,\x}\Theta_{0,\x}G_{ij,\x}^{-1}=\matrix{\frac{a_{mij}(\lambda)}{z-m_{ij}}&\frac{\lambda^{-1}b_{mij}(\lambda)}{(z-m_{ij})^2}+\frac{\lambda^{-1}c_{mij}(\lambda)}{z-m_{ij}}\\0&\frac{-a_{mij}(\lambda)}{z-m_{ij}}}dz\end{equation}
so
$$\Ord_{\infty_{ij}}(G_{ij,\x}\Theta_{0,\x}G_{ij,\x}^{-1})\geq\minimatrix{-1&-1\\\infty&-1}.$$
Regarding the other terms,
$$\Ord_{\infty_{ij}}\left[G_{ij,\x}\left(\zeta_{t,\x}+t(\lambda-1)(\chi_{t,\x}+\Xi_{t,\x})\right)G_{ij,\x}^{-1}\right]\geq \minimatrix{1&0\\0&-1}\star\minimatrix{-1&-2\\0&-1}\star\minimatrix{-1&0\\0&1}=\minimatrix{-1&0\\-2&-1}.$$
$$\Ord_{\infty_{ij}}(G_{ij,\x}dG_{ij,\x}^{-1})\geq\minimatrix{1&0\\0&-1}\star\minimatrix{-2&\infty\\\infty&0}=\minimatrix{-1&0\\-2&-1}.$$
So
$$\Ord_{\infty_{ij}}(\whxi_{ij,t,\x})\geq \minimatrix{-1&-1\\-2&-1}.$$
\cqfd
\subsection{Solution of the Regularity Problem at $m_{ij}$}
\label{section-regularity-solution}
\begin{proposition}
\label{prop-regularity}
For $t$ in a neighborhood of $0$, there exists unique values of the parameters
$$(a_{mij},b_{mij},c_{mij})_{(i,j)\in I}\in(\boW^{\pos})^{3I}$$
depending holomorphically on $t$ and the remaining parameters, such that for all $(i,j)\in I$,
$\whxi_{ij,t,\x}$ extends holomorphically at $m_{ij}$.
\end{proposition}
Proof: define for an arbitrary potential $\xi$:
$$\boR_{ij}(\xi)(\lambda)=\frac{1}{2\pi\ii}\int_{C(\mu_{ij},\varepsilon)}
\left[\xi_{11},\lambda(z-m_{ij}(\lambda))\xi_{12},\lambda\xi_{12}\right].$$
By Proposition \ref{prop-order}, $\whxi_{ij,t,\x}$ has a pole of multiplicity at most
$\minimatrix{1&2\\0&1}$ at $m_{ij}$ so we want to solve the equation
\begin{equation}
\label{eq-boRwhxi}
\boR_{ij}(\whxi_{ij,t,\x})=0.
\end{equation}
Note that $m_{ij}$ depends on $\lambda$, which is why we compute the residue
$\boR_{ij}(\xi)$ as an integral on the fixed circle $C(\mu_{ij},\varepsilon)$.
In this section, we restrict the variable $z$ to the fixed annulus
$\frac{\varepsilon}{2}<|z-\mu_{ij}|<2\varepsilon$ in $\C_{ij}$ where the potential
$\whxi_{ij,t,\x}$ is holomorphic. Let
$$\whzeta_{ij,t,\x}=\zeta_{t,\x}\cdot G_{ij,\x}^{-1}$$
$$\whchi_{ij,t,\x}=G_{ij,\x}\chi_{t,\x} G_{ij,\x}^{-1}$$
$$\whTheta_{ij,t,\x}=G_{ij,\x}\Theta_{t,\x} G_{ij,\x}^{-1}.$$
Then
$$\whxi_{ij,t,\x}=\whzeta_{ij,t,\x}+t(\lambda-1)(\whchi_{ij,t,\x}+\whTheta_{ij,t,\x}).$$
Define for $t\neq 0$:
$$\wczeta_{ij,t,\x}=\frac{1}{t}(\whzeta_{ij,t,\x}-\whzeta_{ij,0,\x})$$
Then $\wczeta_{ij,t,\x}$ extends holomorphically at $t=0$. Let
$$\wtzeta_{ij,t,\x}=\boL_1(\wczeta_{ij,t,\x})$$
where $\boL_1$ is the operator introduced in Proposition \ref{prop-boL} with $\mu=1$,
and we apply it to each element of the matrix $\wczeta_{ij,t,\x}$.
In other words,
$$\wtzeta_{ij,t,\x}(z,\lambda)=\frac{\wczeta_{ij,t,\x}(z,\lambda)-\wczeta_{ij,t,\x}(z,1)}{\lambda-1}.$$
Since $\whzeta_{ij,0,\x}=\xiC$ does not depend on $\lambda$:
$$\whzeta_{ij,t,\x}(z,\lambda)=\whzeta_{ij,t,\x}(z,1)+t(\lambda-1)\wtzeta_{ij,t,\x}(z,\lambda).$$
Hence
$$\whxi_{ij,t,\x}(z,\lambda)=\whzeta_{ij,t,\x}(z,1)+t(\lambda-1)\left[\wtzeta_{ij,t,\x}(z,\lambda)+\whchi_{ij,t,\x}(z,\lambda)+\whTheta_{ij,t,\x}(z,\lambda)\right].$$
Now $G_{ij,\x}(z,1)$ and $\xiC$ are holomorphic at $m_{ij}(1)$,
so $\zeta_{t,\x}(z,1)$ and $\whzeta_{ij,t,\x}(z,1)$ are holomorphic at $m_{ij}(1)$.
Hence
$$\boR_{ij}(\whzeta_{ij,t,\x}(\cdot,1))=0.$$
So to solve Equation \eqref{eq-boRwhxi}, it suffices to solve the following equation:
\begin{equation}
\label{eq-boRwhTheta}
\boR_{ij}(\whTheta_{ij,t,\x})=-\boR_{ij}(\wtzeta_{ij,t,\x}+\whchi_{ij,t,\x}).
\end{equation}
Each term in Equation \eqref{eq-boRwhTheta} is a smooth map from the space of parameters to $(\boW^\pos)^3$ (by composition of various bounded linear operators).
Moreover, the right member does not depend on the parameters $a_{mi'j'}$, $b_{mi'j'}$ and $c_{mi'j'}$ for $(i',j')\in I$, and
the left member depends linearly on $a_{mij}$, $b_{mij}$ and $c_{mij}$.
When $t=0$, we have by Equation \eqref{whThetaij0}:
$$\boR_{ij}(\whTheta_{ij,0,\x})=(a_{mij},b_{mij},c_{mij}).$$
By Proposition \ref{prop-perturbation-isomorphism}, for $t$ small enough, the linear operator
$$(a_{mij},b_{mij},c_{mij})_{(i,j)\in I}\mapsto (\boR(\whTheta_{ij,t,\x}))_{(i,j)\in I}$$
remains an automorphism of $(\boW^\pos)^{3I}$. This means that Equation \eqref{eq-boRwhTheta}
for $(i,j)\in I$ uniquely determines $(a_{mij},b_{mij},c_{mij})_{(i,j)\in I}
\in(\boW^\pos)^{3I}$.
\cqfd
\begin{remark} We see in this proof that although $\Theta_{0,\x}$ is not explicit,
the parameters $a_{mij}$, $b_{mij}$, $c_{mij}$ are determined without having to
invert a linear operator. In a previous version of this work, the term $\Theta_{0,\x}$
was defined explicitly, but then solving the Regularity Problem at $m_{ij}$ required
a quite tedious computation.
\end{remark}
\subsection{Computation of $b_{mij}^0$ and $c_{mij}^0$ at $t=0$}
Let $\x'$ be the collection of the remaining parameters, so $a_{mij}$, $b_{mij}$ and
$c_{mij}$ are now holomorphic functions of $(t,\x')$.
In principle, one can compute explicitely the right member of Equation
\eqref{eq-boRwhTheta} at $t=0$ and obtain the values of $a_{mij}$, $b_{mij}$
and $c_{mij}$ at $t=0$ in function of $\x'$. In particular, this is how the central value
of these parameters (as indicated in Section \ref{section-setup-parameters}) was
computed. We omit this computation because that result will not be used.
We shall need the following easier result:
\begin{proposition}
\label{prop-bmij0}
At $t=0$, we have for $\x'$ in a neighborhood of $\x'_0$:
$$b_{mij}^0(0,\x')=\frac{1}{4(g_{ij}^0)^2}\left[\frac{b_{ij}^0}{m_{ij}^0-q_{ij}}+\frac{b_{ji}^0}{m_{ij}^0-q_{ji}}-\frac{r_{ij}}{(m_{ij}^0-q_{ij})^2}-\frac{r_{ji}}{(m_{ij}^0-q_{ji})^2}\right]$$
$$c_{mij}^0(0,\x')=\frac{1}{4(g_{ij}^0)^2}\left[-\frac{b_{ij}^0}{(m_{ij}^0-q_{ij})^2}-\frac{b_{ji}^0}{(m_{ij}^0-q_{ji})^2}+\frac{2\,r_{ij}}{(m_{ij}^0-q_{ij})^3}+\frac{2\,r_{ji}}{(m_{ij}^0-q_{ji})^3}\right].$$
\end{proposition}
Proof:
Equation \eqref{eq-boRwhTheta} at $t=0$ and $\lambda=0$ gives:
\begin{equation}
\label{eq-bmij0}
b_{mij}^0(0,\x')=-\Res_{m_{ij}^0}\left[(z-m_{ij}^0)\left(
\wtzeta_{ij,0,\x';12}^{(-1)}
+\whchi_{ij,0,\x';12}^{(-1)}\right)\right]
\end{equation}
\begin{equation}
\label{eq-cmij0}
c_{mij}^0(0,\x')=-\Res_{m_{ij}^0}\left[
\wtzeta_{ij,0,\x';12}^{(-1)}
+\whchi_{ij,0,\x';12}^{(-1)}
\right].
\end{equation}
By Proposition \ref{prop-gauging} with $G=G_{ij,\x'}^{-1}$, using the notations
introduced at the end of Section \ref{section-background-DPW-recipe}:
$$\whchi_{ij,0,\x';12}^{(-1)}=
(G_{ij,\x';11}^0)^2\chi_{0,\x';12}^{(-1)}
=\frac{1}{(2g_{ij}^0(z-m_{ij}^0))^2}\left(-\frac{b_{ij}^0\,dz}{z-q_{ij}}-\frac{b_{ji}^0\,dz}{z-q_{ji}}\right)$$
We have
$$\wtzeta_{ij,0,x';12}^{(-1)}=-\wczeta_{ij,0,x';12}^{(-1)}=-\frac{\partial}{\partial t}\whzeta_{ij,t,x';12}^{(-1)}|_{t=0}.$$
By Proposition \ref{prop-gauging} again:
$$\whzeta_{ij,t,\x';12}^{(-1)}=\frac{1}{(2g_{ij}^0(z-m_{ij}^0))^2}\zeta_{t,\x';12}^{(-1)}.$$
By Theorem \ref{thm-openingnodes-derivee}:
$$\frac{\partial}{\partial t}\zeta_{t,\x';12}^{(-1)}|_{t=0}=-\frac{r_{ij}\,dz}{(z-q_{ij})^2}-\frac{r_{ji}\,dz}{(z-q_{ji})^2}.$$
Hence
$$\wtzeta_{ij,0,x';12}^{(-1)}=\frac{1}{(2g_{ij}^0(z-m_{ij}^0))^2}\left(\frac{r_{ij}\,dz}{(z-q_{ij})^2}+\frac{r_{ji}\,dz}{(z-q_{ji})^2}\right).$$
Proposition \ref{prop-bmij0} follows by computing the residues in Equations \eqref{eq-bmij0} and \eqref{eq-cmij0},
using the following elementary residue computation for $c_{mij}^0$:
$$\Res_m\left(\frac{1}{(z-m)^2(z-q)^k}\right)=\frac{-k}{(m-q)^{k+1}}.$$
\cqfd
\section{The Monodromy Problem}
\label{section-monodromy}
\subsection{Definition of domains and paths}
\label{section-monodromy-paths}
We define the domain $\Omega_{t,\x}$ as $\SSigma_{t,\x}$ from which we remove the following sets:
\begin{itemize}
\item $\{\infty_i\}$ for $i\in[1,N]$,
\item $\{\infty_{ij}\}$ and the closed disk $\overline{D}(\mu_{ij},\smallfrac{\varepsilon}{2})$ in $\CC_{ij}$ for $(i,j)\in I$,
\item The closed disk $\overline{D}(\pi_k,\smallfrac{\varepsilon}{2})$ in $\CC_i$ for
$i\in[1,N]$ and $k\in R_i$.
\end{itemize}
If $\x$ is close enough to $\x_0$, the potential $\xi_{t,\x}$ is holomorphic in $\Omega_{t,\x}$.
Also, $\Omega_{t,\x}$ does not depend on $\lambda$, as required for the DPW method.
(This is why we removed a disk centered at $\pi_k$ and not just $p_k$, which depends on $\lambda$.) We first construct an immersion on $\Omega_{t,\x}$. In Section
\ref{section-geometry-immersion}, we will extend it analytically to $\SSigma_{t,\x}$ minus $n$ points corresponding to the Delaunay ends.
\medskip

We define the following fixed domains (independent of $t$, $\x$ and $\lambda$):
$$\Omega_i=\{z\in\C_i:\forall j\in E_i,|z-\pi_{ij}|>\smallfrac{\varepsilon}{2}
\;\mbox{ and }\;
\forall k\in R_i,|z-\pi_k|>\smallfrac{\varepsilon}{2}
\}\quad\mbox{ for $i\in[1,N]$}$$
$$\Omega_{ij}=\{z\in\C_{ij}:
|z-q_{ij}|>\smallfrac{\varepsilon}{2},\;
|z-q_{ji}|>\smallfrac{\varepsilon}{2}\;\mbox{ and }\;
|z-\mu_{ij}|>\smallfrac{\varepsilon}{2}
\} \quad\mbox{ for $(i,j)\in I$}.$$
For $(t,\x)$ close enough to $(0,\x_0)$, the domains
$\Omega_i$ and $\Omega_{ij}$ are included in $\Omega_{t,\x}$.
We fix an arbitrary base point $O_{ij}$ in $\Omega_{ij}$.
\begin{remark}
We would like to take $O_{ij} =0_{ij}$ but if $|\pi_{ij}|=1$, then $\mu_{ij}=0$
so $0_{ij}\not\in \Omega_{ij}$.
Note that we could have assumed without loss of generality that all edges
are non-horizontal, in which case $|\pi_{ij}|\neq 1$ so we could take $O_{ij}=0_{ij}$
which makes some computations slightly simpler.
\end{remark}
\medskip

In the rest of this section, we assume that $t$ is a positive real number and
$(t,\x)$ is close enough to $(0,\x_0)$.
If $(i,j)\in I^*$, we denote $\C_{ij}=\C_{ji}$, $\Omega_{ij}=\Omega_{ji}$ and $O_{ij}=O_{ji}$.
We define the following paths in $\Omega_{t,\x}$ (see Figure \ref{fig-paths}):
\begin{figure}
\includegraphics[width=15cm]{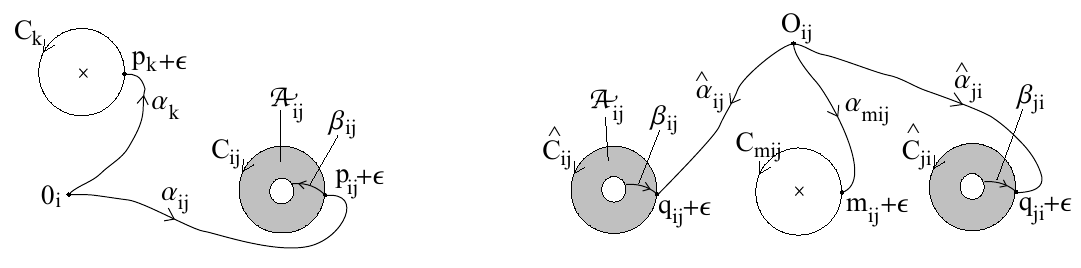}
\caption{Left: paths in $\C_i$. Right: paths in $\C_{ij}$. The shaded annulus $\boA_{ij}$ in
$\C_i$ is identified with the shaded annulus $\boA_{ij}$ in $\C_{ij}$ (via an inversion)
when opening nodes.}
\label{fig-paths}
\end{figure}
\begin{itemize}
\item For $(i,j)\in I\cup I^*$:
\begin{itemize}
\item $\alpha_{ij}$ is a path from $0_i$ to $p_{ij}+\varepsilon$ in
$\Omega_i$, depending continuously on $p_{ij}$. (For example, we can take a fixed path from $0_i$ to $\pi_{ij}$, followed by the segment from $\pi_{ij}$ to $p_{ij}$, which is included in $\Omega_i$ if $|p_{ij}-\pi_{ij}|<\frac{\varepsilon}{2}$.)
\item $C_{ij}$ is the circle in $\Omega_i$ parametrized by $s\mapsto p_{ij}+\varepsilon e^{2\pi\ii s}$.
\item $\whC_{ij}$ is the circle in $\Omega_{ij}$ parametrized by
$s\mapsto q_{ij}+\varepsilon e^{2\pi\ii s}$.
\item $\boA_{ij}$ is the closed annulus bounded by the circles $C_{ij}$ and $\whC_{ij}$
in $\SSigma_{t,\x}$.
\item $\gamma_{ij}=\alpha_{ij}C_{ij}\alpha_{ij}^{-1}\in\pi_1(\Omega_i,0_i)$.
\item $\beta_{ij}$ is a path from $p_{ij}+\varepsilon$ to
to $q_{ij}+\varepsilon$ inside the annulus $\boA_{ij}$, defined as follows using the coordinates $v_{ij}$ and
$w_{ij}$ introduced in Section \ref{section-setup-openingnodes}:
$$v_{ij}(\beta_{ij}(s))=\varepsilon^{1-2s}(t_{ij})^s,\quad s\in[0,1].$$
It goes from the point $v_{ij}=\varepsilon$ to the point $v_{ij}=\frac{t_{ij}}{\varepsilon}$,
which is equivalent to $w_{ij}=\varepsilon$.
Note that the definition of $\beta_{ij}$ depends on the choice of $\arg(t_{ij})$.
Let $r_{ij,0}\neq 0$ be the central value of the parameter $r_{ij}$.
We restrict the parameter $r_{ij}$ to the disk $D(r_{ij,0},|r_{ij,0}|)$,
so $r_{ij}$ does not vanish and we can choose a continuous determination of $\arg(r_{ij})$.
Then since we restricted $t$ to positive values, $\arg(t_{ij})=\arg(r_{ij})$ is well defined, so
$\beta_{ij}$ is well defined and depends continuously on all parameters.
\item $\whalpha_{ij}$ is a fixed path from $O_{ij}$ to $q_{ij}+\varepsilon$ in $\Omega_{ij}$.
\item $\delta_{ij}=\alpha_{ij}\beta_{ij}\whalpha_{ij}^{-1}$ is a path from $0_i$ to $O_{ij}$.
\end{itemize}
\item For $(i,j)\in I$:
\begin{itemize}
\item $\Gamma_{ij}=\delta_{ij}\delta_{ji}^{-1}$ is a path from $0_i$ to $0_j$.
\item $\alpha_{mij}$ is a path from $O_{ij}$ to $m_{ij}+\varepsilon$, depending continuously
on $m_{ij}$.
\item $C_{mij}$ is the circle in $\Omega_{ij}$ parametrized by $s\mapsto m_{ij}+\varepsilon e^{2\pi\ii s}$.
\item $\whgamma_{mij}=\alpha_{mij}C_{mij}\alpha_{mij}^{-1}\in\pi_1(\Omega_{ij},O_{ij})$.
\item $\gamma_{mij}=\delta_{ij}\whgamma_{mij}\delta_{ij}^{-1}\in\pi_1(\Omega_{t,\x},0_i)$.
\end{itemize}
\item For $i\in[1,N]$ and $k\in R_i$:
\begin{itemize}
\item $\alpha_k$ is a path from $0_i$ to $p_k+\varepsilon$ in $\Omega_i$, depending continuously on $p_k$.
\item $C_k$ is the circle in $\Omega_i$ parametrized by $s\mapsto p_k+\varepsilon e^{2\pi\ii s}$.
\item $\gamma_k=\alpha_k C_k\alpha_k^{-1}\in\pi_1(\Omega_i,0_i)$.
\end{itemize}
\end{itemize}
\begin{proposition}
\label{prop-paths}
Any element $\gamma$ in $\pi_1(\Omega_{t,\x},0_1)$ is homotopic to a product of paths or inverse of paths in the following list:
\begin{myenumerate}
\item $\gamma_{ij}$ for $(i,j)\in I\cup I^*$,
\item $\gamma_k$ for $k\in[1,n]$,
\item $\gamma_{mij}$ for $(i,j)\in I$,
\item $\Gamma_{ij}$ for $(i,j)\in I$.
\end{myenumerate}
\end{proposition}
Proof: we denote $\sim$ the homotopy between paths. Let $\gamma\in\pi_1(\Omega_{t,\x},0_1)$. Without loss of generality, we may assume that $\gamma$ is represented by a smooth regular curve which is transverse to the
circles $C_{ij}$ and $\whC_{ij}$ for $(i,j)\in I\cup I^*$.
Without loss of generality, we may also assume that $\gamma$ always intersects a circle
$C_{ij}$ at the point $p_{ij}+\varepsilon$ and a circle $\whC_{ij}$ at the point $q_{ij}+\varepsilon$.
Then we can write
$$\gamma=\prod_{k=1}^r c_k$$
where $c_1(0)=c_r(1)=0_1$ and all other end-points of the paths $c_k$ are either a point
$p_{ij}+\varepsilon$ or a point $q_{ij}+\varepsilon$ with $(i,j)\in I\cup I^*$.
Moreover, the path $c_k$ is included in a domain $\Omega_i$ or $\Omega_{ij}$ if $k$ is odd and an annulus $\boA_{ij}$ if $k$ is even.
\begin{claim}
\label{claim-paths}
$\gamma$ is homotopic to a finite product of paths in the following list:
\begin{itemize}
\item the elements of $\pi_1(\Omega_i,0_i)$ for $i\in[1,N]$,
\item the elements of $\pi_1(\Omega_{ij},O_{ij})$ for $(i,j)\in I$,
\item the paths $\delta_{ij}$ or their inverses for $(i,j)\in I\cup I^*$.
\end{itemize}
\end{claim}
Proof: we define the paths $a_0,\cdots,a_r$ as follows.
For odd $k\in[1,r]$:
\begin{itemize}
\item If $c_k\subset\Omega_i$:
\begin{itemize}
\item If $k=1$ then $c_k(0)=0_1$ so we set $a_0\equiv 0_1$ (a constant path).
\item If $1<k\leq r$, there exists $j\in E_i$ such that $c_k(0)=p_{ij}+\varepsilon$. We set
$a_{k-1}=\alpha_{ij}$.
\item If $k=r$ then $c_k(1)=0_1$ so we set $a_k\equiv 0_1$.
\item If $1\leq k<r$, there exists $\ell\in E_i$ such that $c_k(1)=p_{i\ell}+\varepsilon$.
We set $a_k=\alpha_{i\ell}$.
\end{itemize}
We have $a_{k-1}c_ka_k^{-1}\in\pi_1(\Omega_i,0_i)$.
\item If $c_k\subset\Omega_{ij}$:
\begin{itemize}
\item If $c_k(0)=q_{ij}+\varepsilon$, we set $a_{k-1}=\whalpha_{ij}$.
\item If $c_k(0)=q_{ji}+\varepsilon$, we set $a_{k-1}=\whalpha_{ji}$.
\item If $c_k(1)=q_{ij}+\varepsilon$, we set $a_k=\whalpha_{ij}$.
\item If $c_k(1)=q_{ji}+\varepsilon$, we set $a_k=\whalpha_{ji}$.
\end{itemize}
We have $a_{k-1}c_ka_k^{-1}\in\pi_1(\Omega_{ij},O_{ij})$.
\end{itemize}
We have
$$\gamma\sim \prod_{k=1}^r (a_{k-1}c_ka_k^{-1})$$
For odd $k\in[1,r]$, the path $a_{k-1}c_k a_k^{-1}$ is in the list of Claim \ref{claim-paths}
by construction. So consider some even $k\in[1,r]$.
There exists $(i,j)\in I\cup I^*$ such that $c_k\subset\boA_{ij}$.
\begin{itemize}
\item If $c_k(0)=c_k(1)=p_{ij}+\varepsilon$, then since the fundamental group
$\pi_1(\boA_{ij},p_{ij}+\varepsilon)$ is generated by $C_{ij}$, there exists $\ell\in\Z$ such that
 $c_k\sim C_{ij}^{\ell}$. Then
 $$a_{k-1}c_ka_k^{-1}\sim \alpha_{ij} C_{ij}^{\ell}\alpha_{ij}^{-1}\sim\gamma_{ij}^{\ell}
 \in\pi_1(\Omega_i,0_i).$$
 \item If $c_k(0)=c_k(1)=q_{ij}+\varepsilon$ then since
 $\pi_1(\boA_{ij},q_{ij}+\varepsilon)$ is generated by $\whC_{ij}$, there exists 
 $\ell\in\Z$ such that
 $c_k\sim \whC_{ij}^{\ell}$. Then
 $$a_{k-1}c_ka_k^{-1}\sim \whalpha_{ij} \whC_{ij}^{\ell}\whalpha_{ij}^{-1}\sim\whgamma_{ij}^{\ell}\in\pi_1(\Omega_{ij},O_{ij}).$$
 \item If $c_k(0)=p_{ij}+\varepsilon$ and $c_k(1)=q_{ij}+\varepsilon$ then
 $c_k\beta_{ij}^{-1}\in\pi_1(\boA_{ij},p_{ij}+\varepsilon)$ so there exists $\ell\in\Z$ such
 that $c_k\beta_{ij}^{-1}\sim C_{ij}^{\ell}$. Then
 $$a_{k-1}c_ka_k^{-1}\sim \alpha_{ij}C_{ij}^{\ell}\beta_{ij}\whalpha_{ij}^{-1}
 \sim (\alpha_{ij}C_{ij}^{\ell}\alpha_{ij}^{-1})(\alpha_{ij}\beta_{ij}\whalpha_{ij}^{-1})
 \sim\gamma_{ij}^{\ell}\delta_{ij}.$$
  \item If $c_k(0)=q_{ij}+\varepsilon$ and $c_k(1)=p_{ij}+\varepsilon$ then
  $\beta_{ij}c_k\in\pi_1(\boA_{ij},p_{ij}+\varepsilon)$ so there exists $\ell\in\Z$ such that
  $\beta_{ij}c_k\sim C_{ij}^{\ell}$. Then
  $$a_{k-1}c_ka_k^{-1}\sim\whalpha_{ij}\beta_{ij}^{-1}C_{ij}^{\ell}\alpha_{ij}^{-1}
  \sim(\whalpha_{ij}\beta_{ij}^{-1}\alpha_{ij}^{-1})(\alpha_{ij}C_{ij}^{\ell}\alpha_{ij}^{-1})
  \sim \delta_{ij}^{-1}\gamma_{ij}^{\ell}.$$
\end{itemize}
\cqfd

We return to the proof of Proposition \ref{prop-paths}. Recall that the fundamental group of
the $n$-punctured plane is the free group with $n$ generators.
\begin{itemize}
\item The fundamental group $\pi_1(\Omega_i,0_i)$ is generated by the curves
$\gamma_{ij}$ for $j\in E_i$ and $\gamma_k$ for $k\in R_i$.
\item The fundamental group $\pi_1(\Omega_{ij},O_{ij})$ is generated by
the curves $\whgamma_{ij}$, $\whgamma_{ji}$ and $\whgamma_{mij}$.
Conjugating by $\delta_{ij}$, we see that any element of $\pi_1(\Omega_{ij},O_{ij})$
can be written as a product of $\gamma_{ij}$, $\gamma_{ji}$, $\gamma_{mij}$,
$\delta_{ij}$ and their inverses.
\end{itemize}
Consequently, $\gamma$ is homotopic to a product of paths or inverse of paths in the following list:
\begin{itemize}
\item $\gamma_{ij}$ for $(i,j)\in I\cup I^*$,
\item $\gamma_{mij}$ for $(i,j)\in I$,
\item $\gamma_k$ for $k\in[1,n]$,
\item $\delta_{ij}$ for $(i,j)\in I\cup I^*$.
\end{itemize}
Each path in the first three points is a closed curve based at some point $0_i$ for $i\in[1,N]$, so a path
$\delta_{ij}$ in the decomposition of $\gamma$ must be followed by either the path
$\delta_{ij}^{-1}$, in which case they cancel, or the path $\delta_{ji}^{-1}$, which gives the path $\Gamma_{ij}$ if $(i,j)\in I$ and the path $\Gamma_{ji}^{-1}$ if $(i,j)\in I^*$.\cqfd
\subsection{Formulation of the Monodromy Problem}
\label{section-monodromy-formulation}
Note that the paths in the list of Proposition \ref{prop-paths} are not in $\pi_1(\Omega_{t,\x},0_1)$
(in fact $\Gamma_{ij}$ is not even a closed path)
but they are much more convenient to work with than a set of generators
of $\pi_1(\Omega_{t,\x},0_1)$.
This is why we formulate the Monodromy Problem using the principal solution
(see Section \ref{section-background-principal}).
\begin{proposition}
\label{prop-formulation}
Let $\wtOmega_{t,\x}$ be the universal cover of $\Omega_{t,\x}$. Let $\Phi_{t,\x}$ be the
solution of $d\Phi_{t,\x}=\Phi_{t,\x}\xi_{t,\x}$ in $\wtOmega_{t,\x}$ with initial condition
$\Phi_{t,\x}(0_1)=I_2$. Assume that the Regularity Problem for $\xi_{t,\x}$ at $m_{ij}$ is solved for $(i,j)\in I$ and that:
\begin{equation}
\label{pb-monodromy-nodes}
\mbox{For $(i,j)\in I\cup I^*$,}\qquad
\left\{
\begin{array}{ll}
\boP(\xi_{t,\x},\gamma_{ij})\in\Lambda SU(2)&\quad (i)\\
\boP(\xi_{t,\x},\gamma_{ij})(1)=I_2&\quad (ii)\\
\frac{\partial}{\partial\lambda}
\boP(\xi_{t,\x},\gamma_{ij})(1)=0&\quad (iii)
\end{array}\right.\end{equation}
\begin{equation}
\label{pb-monodromy-ends}
\mbox{For $k\in[1,n]$,}\qquad
\left\{\begin{array}{ll}
\boP(\xi_{t,\x},\gamma_k)\in\Lambda SU(2)&\quad (i)\\
\boP(\xi_{t,\x},\gamma_k)(1)=I_2&\quad (ii)\\
\frac{\partial}{\partial\lambda}\boP(\xi_{t,\x},\gamma_k)(1)=0 &\quad (iii)
\end{array}\right.
\end{equation}
\begin{equation}
\label{pb-monodromy-edges}
\mbox{For $(i,j)\in I$,}\qquad
\left\{
\begin{array}{ll}
\boP(\xi_{t,\x},\Gamma_{ij})\in\Lambda SU(2)&\quad (i)\\
\boP(\xi_{t,\x},\Gamma_{ij})(1)=I_2&\quad (ii)\\
\frac{\partial}{\partial\lambda}\boP(\xi_{t,\x},\Gamma_{ij})(1)=-\ii(V_j-V_i)&\quad (iii)
\end{array}\right.\end{equation}
where $V_i$ represents the position of the vertex $\v_i$ in the $\su(2)$ model
of $\R^3$.
Then $\Phi_{t,\x}$ solves the Monodromy Problem \eqref{pb-monodromy}.
\end{proposition}
We call \eqref{pb-monodromy-nodes} the Monodromy Problem around the nodes,
\eqref{pb-monodromy-ends} the Monodromy Problem at the ends
and \eqref{pb-monodromy-edges} the Monodromy Problem along the edges
(even if $\boP(\xi_{t,\x},\Gamma_{ij})$ is definitely not a monodromy since
$\Gamma_{ij}$ is not a closed curve.)
\medskip

Proof: the Monodromy Problem \eqref{pb-monodromy} for $\Phi_{t,\x}$ is equivalent to Problem
\eqref{pb-principal} for the principal solution of $\xi_{t,\x}$.
Since the Regularity Problem at $m_{ij}$ is solved, we have for all $(i,j)\in I$:
$$\boP(\xi_{t,\x},\gamma_{mij})=I_2.$$
Let $\gamma$ be an arbitrary element of $\pi_1(\Omega_{t,\x},0_1)$. 
By Proposition \ref{prop-paths}, we may write
$$\gamma=\prod_{k=1}^r c_k$$ where for all $k$, $c_k$ or $c_k^{-1}$ is in the list of
Proposition \ref{prop-paths}.
By Equation \eqref{eq-principal-morphism},
$$\boP(\xi_{t,\x},\gamma)=\prod_{k=1}^r \boP(\xi_{t,\x},c_k).$$
Hence $\boP(\xi_{t,\x},\gamma)$ satisfies Equations (i) and (ii) of Problem \eqref{pb-principal}.
Regarding Equation (iii), let us denote, for $k\in[1,r]$, $0_{i_k}=c_k(0)$ and $0_{i_{k+1}}=c_k(1)$.
Using $\boP(\xi_{t,\x},c_k)(1)=I_2$ and Claim \ref{claim-formulation} below:
$$\frac{\partial}{\partial\lambda}\boP(\xi_{t,\x},\gamma)(1)
=\sum_{k=1}^r \frac{\partial}{\partial\lambda}\boP(\xi_{t,\x},c_k)(1)
=-\ii\sum_{k=1}^r(V_{i_{k+1}}-V_{i_k})=0.$$
So Problem \eqref{pb-principal} is solved.
\cqfd
\begin{claim}
\label{claim-formulation}
For $k\in[1,r]$, we have:
\begin{equation}
\label{eq-claim-formulation}
\frac{\partial}{\partial\lambda}\boP(\xi_{t,\x},c_k)(1)=-\ii(V_{i_{k+1}}-V_{i_k}).
\end{equation}
\end{claim}
Proof: if $c_k$ is of type 1,2 or 3 in the list of Proposition \ref{prop-paths},
then $c_k$ is a closed curve so
$i_k=i_{k+1}$ and Equation \eqref{eq-claim-formulation} is true. If $c_k$ is of type 4, then 
Equation \eqref{eq-claim-formulation} is true by Equation (iii) of Problem \eqref{pb-monodromy-edges}.
If $c_k^{-1}$ is in the list of Proposition \ref{prop-paths}, then we write
$$\frac{\partial}{\partial\lambda}\boP(\xi_{t,\x},c_k)(1)
=-\frac{\partial}{\partial\lambda}\boP(\xi_{t,\x},c_k)^{-1}(1)
=-\frac{\partial}{\partial\lambda}\boP(\xi_{t,\x},c_k^{-1})(1)
=\ii(V_{i_k}-V_{i_{k+1}}).$$\cqfd
\begin{remark}
Problems \eqref{pb-monodromy-nodes}, \eqref{pb-monodromy-ends} and
\eqref{pb-monodromy-edges} are stronger than required for the Monodromy
Problem.
The advantage of this stronger formulation is that it only involves
``short'' curves, and the three Problems can be solved essentially independently from
each other using the Implicit Function Theorem.
Also, this stronger formulation yields $\Phi_{t,\x}(\wt0_i,\cdot)\in\Lambda SU(2)$ which
is a strong asset to study the resulting immersion.
\end{remark}
\subsection{The Monodromy Problem around the nodes}
\label{section-monodromy-nodes}
In this Section, we prove:
\begin{proposition}
\label{prop-monodromy-nodes}
Assume that the parameters $a_{mij},b_{mij},c_{mij}$ for $(i,j)\in I$
are given by Proposition \ref{prop-regularity}.
For $t$ small enough, there exists unique values of the parameters
$$\left(a_{ij},b_{ij},c_{ij},\Im(\rho_{ij}^2 r_{ij})\right)_{(i,j)\in I\cup I^*}\in\left((\boW^\pos)^3\times\R\right)^{I\cup I^*}$$
depending smoothly on $t$ and the remaining
parameters, such that the Monodromy Problem \eqref{pb-monodromy-nodes} with respect to $\gamma_{ij}$ is solved for all $(i,j)\in I\cup I^*$.
\end{proposition}
\subsubsection{Preliminaries}
In this section, $t$ is a complex parameter,
and until the end of Section \ref{section-monodromy-nodes-implicit},
the parameter vector $\x$ is free.
Fix a couple $(i,j)\in I\cup I^*$ and let
$$M_{ij}(t,\x)=\boP(\xi_{t,\x},\gamma_{ij}).$$
\begin{proposition}
\label{prop-monodromy-nodes-preliminaries}
\begin{myenumerate}
\item $(t,\x)\mapsto M_{ij}(t,\x)$ is a holomorphic map from a neighborhood of $(0,\x_0)$ to $SL(2,\boW)$.
\item For all $\x$ in a neighborhood of $\x_0$,
$M_{ij}(0,\x)=I_2$.
\item For all $(t,\x)$ in a neighborhood of $(0,\x_0)$, $M_{ij}(t,\x)(1)=I_2$.
\end{myenumerate}
\end{proposition}
Proof:\begin{myenumerate}
\item Using Proposition \ref{prop-substitution},
$(t,\x,z)\mapsto \xi_{t,\x}(z,\cdot)$ is holomorphic for $(t,\x)$ in a neighborhood of $(0,\x_0)$ and $z\in\Omega_i$, with values in $\sl(2,\boW)$.
Point 1 follows from standard Ordinary Differential Equation (ODE) theory.
\item
If $t=0$ then by definition, $\xi_{0,\x}=\xiS$ in $\C_i$
so $\boP(\xi_{0,\x},\gamma_{ij})=I_2$.
\item By Point 3 of Proposition \ref{prop-properties},
 we have if $\lambda=1$:
 $$\boP(\xi_{t,\x},\gamma_{ij})(1)=\matrix{1&\int_{\gamma_{ij}}\beta_{t,\x}\\0&1}=I_2.$$\cqfd
\end{myenumerate}
\medskip

Recall that $\exp$ is a diffeomorphism from a neighborhood of $0$ in the Lie algebra
$\sl(2,\C)$ (respectively $\su(2)$) to a neighborhood of $I_2$ in $SL(2,\C)$
(respectively $SU(2)$). The inverse diffeomorphism is denoted $\log$.
Define for $t\neq 0$:
$$\wcM_{ij}(t,\x)=\frac{1}{t}\log M_{ij}(t,\x).$$
By Point 2 of Proposition \ref{prop-monodromy-nodes-preliminaries},
$\wcM_{ij}(t,\x)$ extends holomorphically at $t=0$ and takes values in
$\sl(2,\boW)$.
Let
$$\wtM_{ij}(t,\x)=\boL_1(\wcM_{ij}(t,\x))$$
$$\whM_{ij}(t,\x)=\boL_1(\lambda\wtM_{ij}(t,\x))$$
where $\boL_1$ is the operator defined in Proposition \ref{prop-boL} with $\mu=1$.
(We apply $\boL_1$ to each element of the matrices
and $\lambda$ stands for the function $\lambda\mapsto\lambda$).
Explicitely, by Point 3 of Proposition \ref{prop-monodromy-nodes-preliminaries},
\begin{equation}
\label{eq-wtMij}
\wtM_{ij}(t,\x)(\lambda)=\frac{\wcM_{ij}(t,\x)}{\lambda-1}=\frac{1}{t(\lambda-1)}\log M_{ij}(t,\x)(\lambda).
\end{equation}
\begin{equation}
\label{eq-whMij}
\whM_{ij}(t,\x)(\lambda)=\frac{1}{\lambda-1}\left(\lambda\wtM_{ij}(t,\x)(\lambda)-\wtM_{ij}(t,\x)(1)\right).
\end{equation}
Since $\boL_1$ is a bounded linear operator, $\wtM_{ij}$ and $\whM_{ij}$ are holomorphic in a neighborhood
of $(0,\x_0)$ with values in $\sl(2,\boW)$.
\begin{proposition}
\item Problem \eqref{pb-monodromy-nodes} is equivalent for real $t\neq 0$ to the following Rescaled Monodromy Problem:
\begin{equation}
\label{pb-monodromy-rescaled}
\left\{\begin{array}{l}
\whM_{ij}(t,\x)\in \Lambda \su(2)\\
\wtM_{ij}(t,\x)(1)=0
\end{array}\right.\end{equation}
\end{proposition}
Proof: by Point 3 of Proposition \ref{prop-monodromy-nodes-preliminaries}, Item (ii) of Problem \eqref{pb-monodromy-nodes} is
automatically satisfied.
If $t\neq 0$, we have by Equation \eqref{eq-wtMij}:
$$\wtM_{ij}(t,\x)(1)=\frac{1}{t}\frac{\partial}{\partial\lambda} M_{ij}(t,\x)(1).$$
So Item (iii) of Problem \eqref{pb-monodromy-nodes} is equivalent to $\wtM_{ij}(t,\x)(1)=0$.
Assuming that this is true, we have by Equations \eqref{eq-wtMij} and \eqref{eq-whMij}:
$$\log M_{ij}(t,\x)(\lambda)=\frac{t(\lambda-1)^2}{\lambda}\whM_{ij}(t,\x)(\lambda).$$
Since $(\lambda-1)^2\lambda^{-1}$ is real on the unit circle,
Item (i) of Problem \eqref{pb-monodromy-nodes} is equivalent for real $t\neq 0$ to $\whM_{ij}(t,\x)\in \Lambda\su(2)$.
\cqfd
\medskip

We define the following maps:
$$\boE_{ij,1}(t,\x)=\whM_{ij;11}(t,\x)+\whM_{ij;11}(t,\x)^*\in\boW$$
$$\boE_{ij,2}(t,\x)=\lambda\left(\whM_{ij;12}(t,\x)+\whM_{ij;21}(t,\x)^*\right)\in\boW$$
$$\boE_{ij,3}(t,\x)=\left(\wtM_{ij;11}(t,\x)(1),\wtM_{ij;12}(t,\x)(1),\wtM_{ij;21}(t,\x)(1)\right)\in\C^3$$
By definition, we have $\boE_{ij,1}(t,\x)=\boE_{ij,1}(t,\x)^*$. So Problem \eqref{pb-monodromy-rescaled} is equivalent to the following problem, using the notations introduced
in Section \ref{section-background-functional-boW}:
\begin{equation}
\label{pb-monodromy-rescaled2}
\left\{\begin{array}{l}
\boE_{ij,1}(t,\x)^+=0\\
\Re(\boE_{ij,1}(t,\x)^0)=0\\
\boE_{ij,2}(t,\x)=0\\
\boE_{ij,3}(t,\x)=0\end{array}\right.\end{equation}
\subsubsection{Computation of $\wtM_{ij}$ and $\whM_{ij}$ at $t=0$.}
Define
$$A_{ij}(\x)=\boL_1(\xi_{0,\x;11}(q_{ij},\cdot))$$
$$B_{ij}(\x)(\lambda)=\lambda \xi_{0,\x;12}(q_{ij},\lambda)$$
$$C_{ij}(\x)=\boL_1(\xi_{0,\x;21}(q_{ij},\cdot)).$$
Clearly, $A_{ij}$, $B_{ij}$ and $C_{ij}$ are holomorphic functions of $\x$
with value in $\boW^{\pos}$.
Using Point 3 of Proposition \ref{prop-properties}, we have
\begin{equation}
\label{eq-xi0qij}
\xi_{0,\x}(q_{ij},\lambda)=\matrix{(\lambda-1)A_{ij}(\x)(\lambda)&\lambda^{-1}B_{ij}(\x)(\lambda)\\(\lambda-1)C_{ij}(\x)(\lambda)&(1-\lambda)A_{ij}(\x)(\lambda)}.
\end{equation}
\begin{proposition}
\label{prop-wtMij0}
At $t=0$, $\wtM_{ij}$ is explicitly given by
\begin{eqnarray*}
\wtM_{ij}(0,\x)(\lambda)&=&\frac{2\pi\ii}{\lambda}\matrix{\lambda a_{ij}+c_{ij}p_{ij}&
-2a_{ij}p_{ij}+b_{ij}-\lambda^{-1}c_{ij}p_{ij}^2\\
\lambda c_{ij}&-\lambda a_{ij}-c_{ij}p_{ij}}\\
&&+\frac{2\pi\ii\,r_{ij}}{\lambda}\matrix{-C_{ij}(\x)&
2A_{ij}(\x)+2\lambda^{-1}p_{ij}C_{ij}(\x)\\
0&C_{ij}(\x)}.\end{eqnarray*}
\end{proposition}
Proof: by Equation \eqref{eq-wtMij}, we have:
$$\wtM_{ij}(0,\x)=\frac{1}{(\lambda-1)}\frac{\partial M_{ij}}{\partial t}(0,\x).$$
Since $\xi_{0,\x}=\xiS$ in $\C_i$, we have by Proposition \ref{prop-principal-derivee}:
$$\frac{\partial M_{ij}}{\partial t}(0,\x)=\int_{\gamma_{ij}}
\PhiS\frac{\partial\xi_{t,\x}}{\partial t}|_{t=0}(\PhiS)^{-1}
=2\pi\ii\Res_{p_{ij}}\left[\PhiS\frac{\partial\xi_{t,\x}}{\partial t}|_{t=0}(\PhiS)^{-1}\right].$$
By Theorem \ref{thm-openingnodes-derivee}, we have in $\C_i$:
$$\frac{\partial\zeta_{t,\x}(z,\lambda)}{\partial t}|_{t=0}
=-\xi_{0,\x}(q_{ij},\lambda)\frac{r_{ij}\,dz}{(z-p_{ij})^2}.$$
Hence by Equations \eqref{eq-xit} and \eqref{eq-xi0qij}:
$$\frac{\partial\xi_{t,\x}(z,\lambda)}{\partial t}|_{t=0}=
-\matrix{(\lambda-1)A_{ij}&\lambda^{-1}B_{ij}\\(\lambda-1)C_{ij}&(1-\lambda)A_{ij}}
\frac{r_{ij}\,dz}{(z-p_{ij})^2}
+(\lambda-1)\matrix{a_{ij}&\lambda^{-1}b_{ij}\\ c_{ij}&-a_{ij}}\frac{dz}{z-p_{ij}}
+O\left((z-p_{ij})^0\right).$$
Proposition \ref{prop-wtMij0} follows from the following elementary
residue computations:
$$\Res_{p}\left[\PhiS(z)\matrix{a&\lambda^{-1}b\\c&-a}
\PhiS(z)^{-1}\frac{1}{(z-p)}\right]=\lambda^{-1}\matrix{
\lambda a+cp&-2ap+b-\lambda^{-1}cp^2\\
\lambda c&-\lambda a-cp}$$
$$\Res_{p}\left[\PhiS(z)\matrix{a&\lambda^{-1}b\\c&-a}
\PhiS(z)^{-1}\frac{1}{(z-p)^2}\right]=\lambda^{-1}
\matrix{c&-2a-2\lambda^{-1}cp\\0&-c}.$$
\cqfd
\medskip

Using Proposition \ref{prop-boL}, we decompose
an arbitrary parameter $x\in\boW^{\pos}$ as
\begin{equation}
\label{eq-x1xt}
x(\lambda)=x(1)+(\lambda-1)\widetilde{x}(\lambda)
\quad\mbox{ with }\quad \widetilde{x}_{ij}\in\boW^{\pos}.
\end{equation}
Recall that the parameters $p_{ij}$
and $r_{ij}$ are in $\boW^0$. Using Proposition \ref{prop-wtMij0} and the definition of $\whM_{ij}$, we obtain:
\begin{proposition}
\label{prop-whMij0}
At $t=0$, $\whM_{ij}$ is explicitly given by
\begin{eqnarray*}
\lefteqn{\whM_{ij}(0,\x)(\lambda)=2\pi\ii\matrix{a_{ij}(1)&\lambda^{-1}c_{ij}(1)p_{ij}^2
-2\lambda^{-1}r_{ij}p_{ij}C_{ij}(\x)(1)\\ c_{ij}(1)&-a_{ij}(1)}}\\
&&+
2\pi\ii\matrix{\lambda \wta_{ij}+\wtc_{ij}p_{ij}-r_{ij}\wtC_{ij}(\x)&
-2\wta_{ij}p_{ij}+\wtb_{ij}-\lambda^{-1}\wtc_{ij}p_{ij}^2+2r_{ij}\wtA_{ij}(\x)+2\lambda^{-1}r_{ij}p_{ij}\wtC_{ij}(\x)\\
\lambda \wtc_{ij}&-\lambda \wta_{ij}-\wtc_{ij}p_{ij}+r_{ij}\wtC_{ij}(\x)}.
\end{eqnarray*}
\end{proposition}
\subsubsection{Solving the Rescaled Monodromy Problem at $t=0$}
In this section, we solve Problem \eqref{pb-monodromy-rescaled} at $t=0$.
Observe that $A_{ij}(\x)$ and $C_{ij}(\x)$
only depend on the parameters $q_{ij}$, $q_{ji}$, $g_{ij}$, $h_{ij}$ and $m_{ij}$.
\begin{proposition}
\label{prop-monodromy-nodes-t0}
At $t=0$, Problem \eqref{pb-monodromy-rescaled} is equivalent to
\begin{equation}
\label{eq-solution-monodromy-nodes-t0}
\left\{\begin{array}{l}
a_{ij}=r_{ij}\left(C_{ij}(\x)(1)+(\lambda-1)\lambda^{-1}\wtC_{ij}(\x)^+\right)\\
b_{ij}=-2 r_{ij} A_{ij}(\x)\\
c_{ij}=-2(\lambda-1)\frac{\overline{p_{ij}}}{(1+|p_{ij}|^2)}r_{ij}C_{ij}(\x)^0\\
\Im\left(r_{ij}C_{ij}(\x)^0\right)=0\\
\end{array}\right.\end{equation}
\end{proposition}
Proof. assume that $\x$ is a solution of Problem \eqref{pb-monodromy-rescaled}
at $t=0$.
By Proposition \ref{prop-wtMij0}, $\boE_{ij,3}(0,\x)=0$ is equivalent to:
\begin{equation}
\label{eq-nodes-8}
\left\{\begin{array}{l}
c_{ij}(1)=0\\
a_{ij}(1)=r_{ij}C_{ij}(\x)(1)\\
b_{ij}(1)=-2 r_{ij} A_{ij}(\x)(1).\end{array}\right.
\end{equation}
By Proposition \ref{prop-whMij0},
\begin{equation}
\label{eq-nodes-1}
\boE_{ij,1}(0,\x)^+=2\pi\ii\left(\lambda\wta_{ij}+\wtc_{ij}^{\,+}p_{ij}-r_{ij}\wtC_{ij}(\x)^+\right)
\end{equation}
\begin{equation}
\label{eq-nodes-4}
\boE_{ij,1}(0,\x)^0=-4\pi\Im\left(\wtc_{ij}^{\,0}p_{ij}-r_{ij}\wtC_{ij}(\x)^0+a_{ij}(1)\right)
\end{equation}
$$\boE_{ij,2}(0,\x)=2\pi\ii\left(-2\lambda\wta_{ij}p_{ij}+\lambda\wtb_{ij}-\wtc_{ij}p_{ij}^2+2\lambda r_{ij}\wtA_{ij}(\x)+2 r_{ij}p_{ij}\wtC_{ij}(\x)-\wtc_{ij}^{\,*}\right).$$
By projection on $\boW^+$, $\boW^-$ and $\boW^0$:
\begin{equation}
\label{eq-nodes-2}
\boE_{ij,2}(0,\x)^+=2\pi\ii\left(-2\lambda\wta_{ij}p_{ij}+\lambda\wtb_{ij}-\wtc_{ij}^{\,+}p_{ij}^2+2\lambda r_{ij}\wtA_{ij}(\x)+2 r_{ij}p_{ij}\wtC_{ij}(\x)^+\right)
\end{equation}
\begin{equation}
\label{eq-nodes-3}
\boE_{ij,2}(0,\x)^-=-2\pi\ii (\wtc_{ij}^{\,+})^*
\end{equation}
\begin{equation}
\label{eq-nodes-5}
\boE_{ij,2}(0,\x)^0=2\pi\ii\left(-\wtc_{ij}^{\,0}p_{ij}^2+2r_{ij}p_{ij}\wtC_{ij}(\x)^0-2r_{ij}p_{ij}C_{ij}(\x)(1)-\overline{\wtc_{ij}^{\,0}}\right).
\end{equation}
Equations \eqref{eq-nodes-3}, \eqref{eq-nodes-1} and \eqref{eq-nodes-2} give:
$$\left\{\begin{array}{l}
\wtc_{ij}^{\,+}=0\\
\wta_{ij}=\lambda^{-1}r_{ij}\wtC_{ij}(\x)^+\\
\wtb_{ij}=-2 r_{ij}\wtA_{ij}(\x).\end{array}\right.$$
Observe that 
$$C_{ij}(\x)^0=C_{ij}(\x)(0)=C_{ij}(\x)(1)-\wtC_{ij}(\x)^0.$$
Equation \eqref{eq-nodes-4} gives
\begin{equation}
\label{eq-nodes-6}
\Im\left(p_{ij}\wtc_{ij}^{\,0}\right)=
-\Im\left(r_{ij}C_{ij}(\x)^0\right).\end{equation}
Equation \eqref{eq-nodes-5} multiplied by $\overline{p_{ij}}$ gives
\begin{equation}
\label{eq-nodes-7}
-|p_{ij}|^2 p_{ij}\wtc_{ij}^{\,0}-2r_{ij}|p_{ij}|^2 C_{ij}(\x)^0-\overline{p_{ij}\wtc_{ij}^{\,0}}=0.
\end{equation}
Taking the imaginary part of Equation \eqref{eq-nodes-7}
and using Equation \eqref{eq-nodes-6}, we obtain
$$\left(-|p_{ij}|^2+2|p_{ij}|^2+1\right)\Im\left(p_{ij}\wtc_{ij}^{\,0}\right)=0.$$
Hence $p_{ij}\wtc_{ij}^{\,0}\in\R$.
Equation \eqref{eq-nodes-7} then gives
$$\wtc_{ij}^{\,0}=\frac{-2 r_{ij}\overline{p_{ij}}C_{ij}(\x)^0}{1+|p_{ij}|^2}.$$
Collecting all results, we obtain \eqref{eq-solution-monodromy-nodes-t0}.
Conversely, assume that the parameters are given by
\eqref{eq-solution-monodromy-nodes-t0}.
Then \eqref{eq-nodes-8} is satisfied, and using Proposition \ref{prop-whMij0}, a computation gives
$$\whM_{ij}(0,\x)(\lambda)=2\pi\ii \frac{r_{ij}C_{ij}(\x)^0}{(1+|p_{ij}|^2)}\matrix{
1-|p_{ij}|^2&-2\lambda^{-1}p_{ij}\\-2\lambda\overline{p_{ij}}&
|p_{ij}|^2-1}\in\Lambda su(2).$$
\cqfd
\medskip

Using Proposition \ref{prop-monodromy-nodes-t0}, we can compute the central value of the parameters
$a_{ij}$, $b_{ij}$ and $c_{ij}$:
\begin{proposition}
\label{prop-aijbijcij}
Assume that $t=0$ and the parameters $p_{ij}$, $q_{ij}$, $r_{ij}$,
$g_{ij}$, $h_{ij}$ and $m_{ij}$ have their central value, as indicated in
Section \ref{section-setup-parameters}.
For $(i,j)\in I\cup I^*$, the Monodromy Problem with respect to $\gamma_{ij}$ is equivalent to
$$\left\{\begin{array}{l}
a_{ij}=2\tau_{ij}\\
b_{ij}=2\rho_{ij}^{-1}\tau_{ij}\\
c_{ij}=-2(\lambda-1)\rho_{ij}\tau_{ij}\end{array}\right.
$$
\end{proposition}
Proof: recalling the notations introduced at the end of Section \ref{section-setup-parameters}, we have:
\begin{equation}
\label{eq-2qij-2muij}
2(q_{ij}-\mu_{ij})=q_{ij}-q_{ji}=\frac{1}{\overline{\pi_{ij}}}+\pi_{ij}=\frac{1}{\rho_{ij}}.
\end{equation}
Equation \eqref{eq-xiC.Gij0} gives at the central value:
$$\xi_{0,\x_0}(q_{ij},\lambda)=\matrix{(\lambda-1)\rho_{ij}&1\\(1-\lambda^2)\rho_{ij}^2&(1-\lambda)\rho_{ij}}.$$
Hence
$$A_{ij}(\x_0)=\rho_{ij}\qquad\mbox{ and }\qquad C_{ij}(\x_0)=-(\lambda+1)\rho_{ij}^2$$
which gives
$$C_{ij}(\x_0)^0=-\rho_{ij}^2,\qquad
C_{ij}(\x_0)(1)=-2\rho_{ij}^2\qquad\mbox{ and }\qquad \wtC_{ij}(\x_0)=-\rho_{ij}^2.$$
Then
$$\Im(r_{ij}C_{ij}(\x_0)^0)=\Im(\tau_{ij})=0.$$
By substitution in Equation \eqref{eq-solution-monodromy-nodes-t0}, we obtain the values of $a_{ij}$, $b_{ij}$
and $c_{ij}$ as indicated in Proposition \ref{prop-aijbijcij}.
\cqfd
\subsubsection{Solving the rescaled Monodromy Problem for $t\neq 0$}
\label{section-monodromy-nodes-implicit}
Define
$$\boE_{ij}(t,\x)=\left(\boE_{ij,1}(t,\x)^+,\boE_{ij,2}(t,\x)^+,(\boE_{ij,2}(t,\x)^-)^*,\boE_{ij,3}(t,\x),\boE_{ij,2}(t,\x)^0,\Re(\boE_{ij,1}(t,\x)^0)\right)\in (\boW^+)^3\times\C^4\times \R.$$
Problem \eqref{pb-monodromy-rescaled2} is equivalent to $\boE_{ij}(t,\x)=0$.
\begin{proposition}
\label{prop-monodromy-nodes-differential}
For $(i,j)\in I\cup\I$,
the partial differential of $\boE_{ij}$ at $(0,\x_0)$ with respect to the variable
$$y_{ij}=\left(\lambda\wta_{ij},\lambda\wtb_{ij},\wtc_{ij}^{\,+},a_{ij}(1),b_{ij}(1),c_{ij}(1),\wtc_{ij}^{\,0},
\Im(\rho_{ij}^2 r_{ij})\right)$$
is an automorphism of $(\boW^+)^3\times\C^4\times \R$.
\end{proposition}
Proof: by Equations \eqref{eq-nodes-1}, \eqref{eq-nodes-2} and \eqref{eq-nodes-3},
the partial differential of $(\boE_{ij,1}^+,\boE_{ij,2}^+,(\boE_{ij,2}^-)^*)$ with respect
to $(\lambda\wta_{ij},\lambda\wtb_{ij},\wtc_{ij}^{\,+})$ is a matrix-type operator
(see Definition \ref{def-matrixtype}) from $(\boW^+)^3$ to itself with matrix
$$\left(\begin{array}{ccc}1&0&\pi_{ij}\\-2\pi_{ij}&1&-\pi_{ij}^2\\0&0&1\end{array}\right)
\in GL(3,\C)$$
so is an automorphism of $(\boW^+)^3$.
By Proposition \ref{prop-fredholm}, it suffices to prove that the partial differential of $\boE_{ij}$
with respect to $y_{ij}$ is injective.
Since at $t=0$, the map $y_{ij}\mapsto\boE_{ij}$ is affine, this is equivalent to
proving that $y_{ij}\mapsto\boE_{ij}$ is injective at $t=0$. This follows from Proposition \ref{prop-monodromy-nodes-t0} since we have found a unique solution $y_{ij}$.
\cqfd
\medskip

We now prove Proposition \ref{prop-monodromy-nodes}.
We decompose $\x=(\x',\x'',\x''')$ where
$$\x'=(a_{mij},b_{mij},c_{mij})_{(i,j)\in I}$$
$$\x''=(y_{ij})_{(i,j)\in I\cup \I}$$
and $\x'''$ denotes the remaining parameters.
Proposition \ref{prop-regularity} determines $\x'$ as a smooth function of $t$, $\x''$ and $\x'''$
so we write $\x'=\x'(t,\x'',\x''')$.
Define
$$\boE(t,\x)=(\boE_{ij}(t,\x))_{(i,j)\in I\cup \I}$$
$$\boF(t,\x'',\x''')=\boE(t,\x'(t,\x'',\x'''),\x'',\x''')$$
Since $\boE(0,\x)$ does not depend on $\x'$, we have $d_{\x'}\boE(0,\x_0)=0$ so by the chain rule,
$$d_{\x''}\boF(0,\x''_0,\x'''_0)=d_{\x''}\boE(0,\x_0).$$
By Proposition \ref{prop-monodromy-nodes-differential}, $d_{\x''}\boE(0,\x_0)$ is an automorphism (it has block diagonal form).
By the Implicit Function Theorem, for $(t,\x''')$ in a neighborhood of $(0,\x'''_0)$, there
exists a unique $\x''$, depending smoothly on $(t,\x''')$ such that
$\boF(t,\x''(t,\x'''),\x''')=0$.
\cqfd
\subsection{The Monodromy Problem at the ends}
\label{section-monodromy-ends}
\begin{proposition}
\label{prop-monodromy-ends}
Assume that the parameters
are as in Proposition \ref{prop-monodromy-nodes}.
For $t$ small enough, there exists unique values of the parameters
$$(a_k,b_k,p_k)_{1\leq k\leq n}\in(\boW^\pos)^{3n}$$
 depending smoothly on $t$ and the
remaining parameters, such that the Monodromy Problem \eqref{pb-monodromy-ends}
with respect to $\gamma_k$ is solved for $k\in[1,n]$ and the
following normalisations hold:
$$\forall k\in[1,n],\quad\Re( a_k^0)=\tau_k
\quad\mbox{ and }\quad p_k^0=\pi_k$$
Moreover, at $t=0$, we have
$$a_k=\tau_k,\qquad
b_k=\frac{-2\tau_k\overline{\pi_k}}{1+|\pi_k|^2}
\quad \mbox{ and }\quad
p_k=\pi_k.$$
\end{proposition}
This is proven in Proposition 3 of \cite{nnoids} using the Implicit Function Theorem,
in a way similar to Section \ref{section-monodromy-nodes}.
Of course, the potential is different in that paper, but the only properties of
the potential that are used are the following, which are satisfied by our potential $\xi_{t,\x}$:
\begin{equation}
\label{eq-properties-nnoids}
\begin{array}{ll}
\bullet&\mbox{ $\xi_{0,\x}=\xiS$ in the Riemann sphere containing $p_k$,}\\
\bullet &\mbox{for all $t$, $\xi_{t,\x}$ has a pole at $p_k$ with principal part}\\
&\displaystyle\xi_{t,\x}=t(\lambda-1)^2\matrix{0&0\\1&0}\left(\frac{a_k}{(z-p_k)^2}+\frac{b_k}{z-p_k}\right)dz
+O((z-p_k)^0).
\end{array}\end{equation}
\cqfd
\subsection{The Monodromy Problem along the edges.}
\label{section-monodromy-edges}
\begin{definition}
\label{def-tlogt}
Let $f(t)$ be a function of the real variable $t\geq 0$. We say that $f$ is a smooth function of $t$ and $t\log t$
if there exists a smooth function of two variables $g(t,s)$ defined in a neighborhood of
$(0,0)$ in $\R^2$ such that $f(t)=g(t,t\log t)$ for
$t>0$ and $f(0)=g(0,0)$.
\end{definition}
\begin{remark}
The function $t\log t$ extends continuously at $0$ but the extension is not differentiable at $0$ and is only of
H\"older class $C^{0,\alpha}$ for all $\alpha\in (0,1)$.
Therefore, a smooth function of $t$ and $t\log t$ is only of class $C^{0,\alpha}$.
\end{remark}
\begin{proposition}
\label{prop-monodromy-edges}
Assume that the parameters are as in Proposition \ref{prop-monodromy-ends}.
For $t>0$ small enough, there exists unique values of the parameters
$$(g_{ij},h_{ij},m_{ij},p_{ji})_{(i,j)\in I}\in((\boW^\pos)^3\times\C)^I$$
depending smoothly on $t$, $t\log t$
and the remaining parameters, such that the Monodromy Problem
\eqref{pb-monodromy-edges} with respect to $\Gamma_{ij}$ is solved for all $(i,j)\in I$, up to
one real equation (Equation (v) of Problem \eqref{pb-monodromy-edges2}) which we will solve in Section \ref{section-balancing} using the non-degeneracy hypothesis.
\end{proposition}
\subsubsection{Preliminaries}
Until the end of Section \ref{section-monodromy-edges-implicit},
the parameter vector $\x$ is free.
Define for $(i,j)\in I$ and $t>0$:
$$P_{ij}(t,\x)=\boP(\xi_{t,\x},\Gamma_{ij}).$$
\begin{proposition}
\label{prop-Pij}
For $(i,j)\in I$ and $\x$ in a neighborhood of $\x_0$,
$P_{ij}(t,\x)$ extends at $t=0$ as a smooth function of $t$, $t\log t$ and $\x$
with value in $SL(2,\boW)$.
Moreover, at $t=0$, we have:
$$P_{ij}(0,\x)(\lambda)=\PhiS(p_{ij},\lambda)\left[\PhiC(q_{ij},\lambda)G_{ij,\x}(q_{ij},\lambda)\right]^{-1}\left[\PhiC(q_{ji},\lambda)G_{ij,\x}(q_{ji},\lambda)\right]\PhiS(p_{ji},\lambda)^{-1}.$$
\end{proposition}
\begin{remark}
\label{remark-Pij0}
Observe that $P_{ij}(0,\x_0)$ is equal to $\boP(\xi_0,\Gamma_{ij})$
as given by Equation \eqref{eq-strategy-boPGammaij} in Section
\ref{section-strategy-potential}. This means that our definition of the principal solution
of $\xi_0$ on the noded Riemann surface $\SSigma_0$ was the right one.
\end{remark}
Proof: 
we consider the principal solution of $\xi_{t,\x}$ on each path in the definition of $\Gamma_{ij}$
in \ref{section-monodromy-paths}.
Let $(i,j)\in I\cup I^*$.
For ease of notation, we omit the $\lambda$ variable.
\begin{myenumerate}
\item By standard ODE theory, $\boP(\xi_{t,\x},\alpha_{ij})$ is a smooth function of $(t,\x)$ in a neighborhood
of $(0,\x_0)$. Since $\xi_{0,\x}=\xiS$ in $\Omega_i$, we have
(using the notation explained in Remark \ref{remark-principal}):
\begin{equation}
\label{eq-edges-principal-1}
\boP(\xi_{0,\x},\alpha_{ij})=\boP(\xiS,0_i,p_{ij}+\varepsilon)=\PhiS(p_{ij}+\varepsilon).
\end{equation}
\item By standard ODE theory, $\boP(\xi_{t,\x},\whalpha_{ij}^{-1})$ is a smooth function of $(t,\x)$ in a neighborhood
of $(0,\x_0)$. Since $\xi_{0,\x}=\xiC\cdot G_{ij,\x}$ in $\Omega_{ij}$, we have
\begin{equation}
\label{eq-edges-principal-2}
\boP(\xi_{0,\x},\whalpha_{ij}^{-1})=\boP(\xiC\cdot G_{ij},q_{ij}+\varepsilon,O_{ij})=
[\PhiC(q_{ij}+\varepsilon)
G_{ij,\x}(q_{ij}+\varepsilon)]^{-1}
\PhiC(O_{ij})G_{ij,\x}(O_{ij}).
\end{equation}

\item To evaluate $\boP(\xi_{t,\x},\beta_{ij})$, we use Theorem \ref{thm-principal} in
appendix \ref{appendix-principal}. We temporarily see $t$ as a complex number and fix
the value of $\x$.
The path $\beta_{ij}$ and the principal solution $\boP(\xi_{t,\x},\beta_{ij})$
depend on the choice of the argument of $t$.
For $t\neq 0$, let
$$F(t,\x)=\exp\left(-\frac{\log t_{ij}}{2\pi\ii}\log\boP(\xi_{t,\x},C_{ij})\right)\boP(\xi_{t,\x},\beta_{ij}).$$
By Theorem \ref{thm-principal}, $F(t,\x)$ is a well defined holomorphic function of $t$ and extends
holomorphically at $t=0$ with
\begin{equation}
\label{eq-edges-principal-F0}
F(0)=\boP(\xiS,p_{ij}+\varepsilon,p_{ij})\boP(\xiC\cdot G_{ij,\x},q_{ij},q_{ij}+\varepsilon).
\end{equation}
(We apply Theorem \ref{thm-principal} with $z=v_{ij}$, $t=t_{ij}$, $\xi_t=\xi_{t,\x}$, $\xi_0=\xiS$
and $\whxi_0=\xiC\cdot G_{ij,\x}$.)
For fixed $t$, $F(t,\x)$ depends holomorphically on $\x$ by Proposition \ref{prop-substitution}, so is a holomorphic function of $(t,\x)$ in a neighborhood of $(0,\x_0)$ by Hartog Theorem on separate holomorphy
(see Section \ref{section-background-functional-holomorphy}).
 \medskip
 
Since the path $C_{ij}$ lies in the fixed domain $\Omega_i$ where $\xi_{t,\x}$ is
holomorphic, and $\xi_{0,\x}=\xiS$,
the map $(t,\x)\mapsto\log\boP(\xi_{t,\x},C_{ij})$ is holomorphic in a neighborhood of $(0,\x_0)$ and vanishes at $t=0$. Hence we
can write
$$\log\boP(\xi_{t,\x},C_{ij})=t f(t,\x)$$
where $f$ is holomorphic in a neighborhood of $(0,\x_0)$.
\medskip

Next we restrict $t$ to positive values and recalling that $t_{ij}=t\,r_{ij}$, we define
$$g(t,s,\x)=\exp\left(\frac{(s +t\log r_{ij})}{2\pi\ii}f(t,\x)\right)F(t,\x).$$
Then
$$\boP(\xi_{t,\x},\beta_{ij})=g(t,t\log t,\x)$$
so $\boP(\xi_{t,\x},\beta_{ij})$ is a smooth function of $t$, $t\log t$ and $\x$.
By Equation \eqref{eq-edges-principal-F0}, we have:
\begin{equation}
\label{eq-edges-principal-3}
g(0,0,\x)=F(0,\x)=
\PhiS(p_{ij}+\epsilon)^{-1}\PhiS(p_{ij})
[\PhiC(q_{ij})G_{ij,\x}(q_{ij})]^{-1}\PhiC(q_{ij}+\varepsilon)G_{ij,\x}(q_{ij}+\varepsilon).
\end{equation}
\end{myenumerate}
By Equation \eqref{eq-principal-morphism}, we have
$$\boP(\xi_{t,\x},\delta_{ij})=\boP(\xi_{t,\x},\alpha_{ij})\boP(\xi_{t,\x},\beta_{ij})\boP(\xi_{t,\x},\whalpha_{ij}^{-1}).$$
Hence $\boP(\xi_{t,\x},\delta_{ij})$ is a smooth function of $t$, $t\log t$ and $\x$ and
its value at $t=0$ is obtained by multiplying Equations \eqref{eq-edges-principal-1}, \eqref{eq-edges-principal-3} and
\eqref{eq-edges-principal-2} in this order, which gives for $(i,j)\in I\cup I^*$:
\begin{equation}
\label{eq-edges-principal-4}
\boP(\xi_{0,\x},\delta_{ij})=\PhiS(p_{ij})
[\PhiC(q_{ij})G_{ij,\x}(q_{ij})]^{-1}\PhiC(O_{ij})G_{ij,\x}(O_{ij}).
\end{equation}
Let $(i,j)\in I$. By Equation \eqref{eq-principal-morphism}, we have
$$\boP(\xi_{t,\x},\Gamma_{ij})=\boP(\xi_{t,\x},\delta_{ij})\boP(\xi_{t,\x},\delta_{ji})^{-1}.$$
Proposition \ref{prop-Pij} follows from Equation \eqref{eq-edges-principal-4},
remembering that we have defined $O_{ji}=O_{ij}$.
\cqfd
\medskip
\subsubsection{Computation of $P_{ij}$ at the central value}
By Remark \ref{remark-Pij0}, $P_{ij}(0,\x_0)$ is equal to $\boP(\xi_0,\Gamma_{ij})$
which we have already computed in Section \ref{section-strategy-immersion}: it
is given by Equation \eqref{eq-strategy-boPGammaijfinal}.
To simplify the result, we introduce the matrix
$$H_{ij}(\lambda)=\FS(\pi_{ij},\lambda)\in\Lambda SU(2).$$
We compute
\begin{eqnarray}
\nonumber
H_{ij}(\lambda)\matrix{\lambda&0\\0&\lambda^{-1}}H_{ij}(\lambda)^{-1}
&=&\frac{1}{1+|\pi_{ij}|^2}
\matrix{1&-\lambda^{-1}\pi_{ij}\\\lambda\overline{\pi_{ij}}&1}
\matrix{\lambda&0\\0&\lambda^{-1}}
\matrix{1&\lambda^{-1}\pi_{ij}\\-\lambda\overline{\pi_{ij}}&1}\\
\nonumber
&=&\frac{1}{1+|\pi_{ij}|^2}
\matrix{\lambda^{-1}|\pi_{ij}|^2+\lambda&\pi_{ij}(\lambda^{-2}-1)\\
\overline{\pi}_{ij}(1-\lambda^2)&\lambda|\pi_{ij}|^2+\lambda^{-1}}\\
&=&P_{ij}(0,\x_0)
\label{eq-Ptij0}
\end{eqnarray}
\begin{remark}
A computation gives
$$\frac{-\ii}{2}H_{ij}(1)\matrix{-1&0\\0&1}H_{ij}(1)^{-1}
=\frac{-\ii}{2}\frac{1}{(1+|\pi_{ij}|^2)}\matrix{|\pi_{ij}|^2-1&2\pi_{ij}\\2\overline{\pi_{ij}}&1-|\pi_{ij}|^2}=\pi_S^{-1}(\pi_{ij})=\u_{ij}.$$
So $H_{ij}$ acts by conjugation on $\su(2)$ as a rigid motion whose linear part is a rotation which maps
the vertical vector $(0,0,1)$ to $\u_{ij}$. So essentially, what we are doing here is rotating the graph $\Gamma$ so that the vector $\u_{ij}$ becomes vertical.
\end{remark}
In view of Equation \eqref{eq-Ptij0}, we define for $t$ small enough:
$$\wtP_{ij}(t,\x)(\lambda)=\log\left[H_{ij}(\lambda)^{-1}P_{ij}(t,\x)(\lambda)
H_{ij}(\lambda)\matrix{\lambda^{-1}&0\\0&\lambda}\right]$$
so that $\wtP_{ij}(0,\x_0)=\log I_2=0$.
\subsubsection{Differential of $\wtP_{ij}$ at the central value}
\label{section-monodromy-edges-dwtPij}
\begin{proposition}
\label{prop-dwtPij}
The partial differential of $\wtP_{ij}$ with respect to $\x$ at $(0,\x_0)$ is given by the following formula:
\begin{eqnarray*}
\lefteqn{d_{\x}\wtP_{ij}(0,\x_0)=
\rho_{ij}\matrix{-1&\frac{\lambda^{-1}}{\overline{\pi_{ij}}}\\-\lambda\overline{\pi_{ij}}&1}dp_{ij}
+\rho_{ij}\matrix{1&\frac{-\lambda}{\overline{\pi_{ij}}}\\\lambda^{-1}\overline{\pi_{ij}}&-1}dp_{ji}}\\
&&+\matrix{2&\frac{-(\lambda+1)}{\overline{\pi_{ij}}}\\(\lambda+1)\overline{\pi_{ij}}&-2}dg_{ij}
+(\lambda-1)\matrix{0&\frac{-1}{\overline{\pi_{ij}}}\\\overline{\pi_{ij}}&0}dh_{ij}
+2(1-\lambda)\rho_{ij}\matrix{0&\frac{1}{\overline{\pi_{ij}}}\\\overline{\pi_{ij}}&0}dm_{ij}
\end{eqnarray*}
\end{proposition}
Proof:
By Proposition \ref{prop-Pij}, we have (omiting $\lambda$)
$$\wtP_{ij}(0,\x)=\log\left[H_{ij}^{-1}\PhiS(p_{ij})G_{ij,\x}(q_{ij})^{-1}
\PhiC(q_{ij})^{-1}\PhiC(q_{ji})G_{ij,\x}(q_{ji})\PhiS(p_{ji})^{-1}H_{ij}
\matrix{\lambda^{-1}&0\\0&\lambda}\right].$$
The parameters which appear in this formula are $p_{ij}$, $p_{ji}$, $g_{ij}$,
$h_{ij}$ and $m_{ij}$. To compute the differential of $\wtP_{ij}$,
we may assume that these parameters are complex numbers,
and then use Proposition \ref{prop-substitution}.
To compute the partial derivatives of $\wtP_{ij}$ with respect to each parameter, we use the following two identities: let $A$, $B(y)$, $C$ be three $SL(2,\C)$ matrices such that
$B(y)$ is a holomorphic function of $y$ in a neighborhood of $y_0$.
Assume that $AB(y_0)C=I_2$ and define $M(y)=\log(AB(y)C)$ in a neighborhood
of $y_0$.
Then
\begin{equation}
\label{identity-left}
M'(y_0)=AB'(y_0)B(y_0)^{-1}A^{-1}
\end{equation}
\begin{equation}
\label{identity-right}
M'(y_0)=C^{-1}B(y_0)^{-1}B'(y_0)C.
\end{equation}
In the following computations, we write $p=\pi_{ij}$ for short.
Using Identity \eqref{identity-left} with
$$A=H_{ij}^{-1},\qquad
B(p_{ij})=\PhiS(p_{ij})$$ we obtain
\begin{eqnarray*}
\frac{\partial\wtP_{ij}}{\partial p_{ij}}(0,\x_0)
&=&\frac{1}{1+|p|^2}\matrix{1&-\lambda^{-1}p\\\lambda\overline{p}&1}
\matrix{0&\lambda^{-1}\\0&0}
\matrix{1&-\lambda^{-1}p\\0&1}
\matrix{1&\lambda^{-1}p\\-\lambda\overline{p}&1}\\
&=&\frac{1}{1+|p|^2}\matrix{-\overline{p}&\lambda^{-1}\\
-\lambda\overline{p}^2&\overline{p}}
\end{eqnarray*}
Using Identity \eqref{identity-right} with
$$B(p_{ji})=\PhiS(p_{ji})^{-1},\qquad
C=H_{ij}\matrix{\lambda^{-1}&0\\0&\lambda}=\frac{1}{\sqrt{1+|p|^2}}
\matrix{\lambda^{-1}&p\\-\overline{p}&\lambda}$$
we obtain
\begin{eqnarray*}
\frac{\partial\wtP_{ij}}{\partial p_{ji}}(0,\x_0)&=&
\frac{1}{1+|p|^2}\matrix{\lambda&-p\\\overline{p}&\lambda^{-1}}
\matrix{1&\frac{-\lambda^{-1}}{\overline{p}}\\0&1}\matrix{0&-\lambda^{-1}\\0&0}
\matrix{\lambda^{-1}&p\\-\overline{p}&\lambda}\\
&=&\frac{1}{1+|p|^2}\matrix{\overline{p}&-\lambda\\\lambda^{-1}\overline{p}^2&-\overline{p}}
\end{eqnarray*}
The parameter $g_{ij}$ appears in both $G_{ij,\x}(q_{ij})$ and
$G_{ij,\x}(q_{ji})$.
We compute, from the definition of $G_{ij,\x}$:
$$G_{ij,\x_0}(z)^{-1}\frac{\partial G_{ij,\x}(z)}{\partial g_{ij}}|_{\x=\x_0}=
-\frac{\partial G_{ij,\x}^{-1}(z)}{\partial g_{ij}}|_{\x=\x_0}G_{ij,\x_0}(z)
=\matrix{\lambda-1&2(z-\mu_{ij})\\\frac{\lambda-\lambda^2}{2(z-\mu_{ij})}&1-\lambda}.$$
Using Identity \eqref{identity-left} with
$$A=H_{ij}^{-1}\PhiS(\pi_{ij})=\BS(\pi_{ij}),\qquad
B(g_{ij})=G_{ij,\x}(q_{ij})^{-1}$$
 and Identity \eqref{identity-right} with
$$B(g_{ij})=G_{ij,\x}(q_{ji}),\qquad
C=\PhiS(\pi_{ji})^{-1}H_{ij}\matrix{\lambda^{-1}&0\\0&\lambda}
=\frac{1}{\sqrt{1+|p|^2}}\matrix{0&\frac{1+|p|^2}{\overline{p}}\\-\overline{p}&\lambda}$$
we obtain, using Equation \eqref{eq-2qij-2muij}:
\begin{eqnarray*}
\frac{\partial \wtP_{ij}}{\partial g_{ij}}(0,\x_0)&=&
\frac{1}{1+|p|^2}\matrix{1&0\\\lambda\overline{p}&1+|p|^2}
\matrix{1-\lambda&\frac{-1-|p|^2}{\overline{p}}\\\frac{(\lambda^2-\lambda)\overline{p}}{1+|p|^2}&\lambda-1}\matrix{1+|p|^2&0\\-\lambda\overline{p}&1}\\
&&+\frac{1}{1+|p|^2}\matrix{\lambda&\frac{-1-|p|^2}{\overline{p}}\\\overline{p}&0}
\matrix{\lambda-1&\frac{-1-|p|^2}{\overline{p}}\\\frac{(\lambda^2-\lambda)\overline{p}}{1+|p|^2}&1-\lambda}
\matrix{0&\frac{1+|p|^2}{\overline{p}}\\-\overline{p}&\lambda}\\
&=&\matrix{1&\frac{-1}{\overline{p}}\\\lambda\overline{p}&-1}
+\matrix{1&\frac{-\lambda}{\overline{p}}\\\overline{p}&-1}
\end{eqnarray*}
The computations of the partial derivatives with respect to $h_{ij}$ and $m_{ij}$ are
similar:
$$G_{ij,\x_0}(z)^{-1}\frac{\partial G_{ij,\x}(z)}{\partial h_{ij}}|_{\x=\x_0}=
\matrix{-\lambda&-2(z-\mu_{ij})\\\frac{\lambda^2-\lambda}{2(z-\mu_{ij})}&\lambda}$$
\begin{eqnarray*}
\frac{\partial \wtP_{ij}}{\partial h_{ij}}(0,\x_0)&=&
\frac{1}{1+|p|^2}
\matrix{1&0\\\lambda\overline{p}&1+|p|^2}
\matrix{\lambda&\frac{1+|p|^2}{\overline{p}}\\\frac{(\lambda-\lambda^2)\overline{p}}{1+|p|^2}&-\lambda}
\matrix{1+|p|^2&0\\-\lambda\overline{p}&1}\\
&&+\frac{1}{1+|p|^2}
\matrix{\lambda&\frac{-1-|p|^2}{\overline{p}}\\\overline{p}&0}
\matrix{-\lambda&\frac{1+|p|^2}{\overline{p}}\\\frac{(\lambda-\lambda^2)\overline{p}}{1+|p|^2}&\lambda}
\matrix{0&\frac{1+|p|^2}{\overline{p}}\\-\overline{p}&\lambda}\\
&=&\matrix{0&\frac{1}{\overline{p}}\\\lambda\overline{p}&0}
+\matrix{0&\frac{-\lambda}{\overline{p}}\\-\overline{p}&0}
\end{eqnarray*}
$$G_{ij,\x_0}(z)^{-1}\frac{\partial G_{ij,\x}(z)}{\partial m_{ij}}|_{\x=\x_0}=
\matrix{\frac{1-\lambda}{z-\mu_{ij}}&-2\\\frac{\lambda^2-\lambda}{2(z-\mu_{ij})^2}&
\frac{\lambda-1}{z-\mu_{ij}}}$$
\begin{eqnarray*}
\frac{\partial \wtP_{ij}}{\partial m_{ij}}(0,\x_0)&=&
\frac{1}{1+|p|^2}
\matrix{1&0\\\lambda\overline{p}&1+|p|^2}
\matrix{\frac{2(\lambda-1)\overline{p}}{1+|p|^2}& 2\\
\frac{2(\lambda-\lambda^2)\overline{p}^2}{(1+|p|^2)^2}&\frac{2(1-\lambda)\overline{p}}{1+|p|^2}}
\matrix{1+|p|^2&0\\-\lambda\overline{p}&1}\\
&&+\frac{1}{1+|p|^2}
\matrix{\lambda&\frac{-1-|p|^2}{\overline{p}}\\\overline{p}&0}
\matrix{\frac{2(\lambda-1)\overline{p}}{1+|p|^2}&-2\\\frac{2(\lambda^2-\lambda)\overline{p}^2}{(1+|p|^2)^2}&\frac{2(1-\lambda)\overline{p}}{1+|p|^2}}
\matrix{0&\frac{1+|p|^2}{\overline{p}}\\-\overline{p}&\lambda}\\
&=&\frac{2}{1+|p|^2}\matrix{-\overline{p}&1\\-\lambda\overline{p}^2&\overline{p}}
+\frac{2}{1+|p|^2}\matrix{\overline{p}&-\lambda\\\overline{p}^2&-\overline{p}}
\end{eqnarray*}
\cqfd

We define
$$\whP_{ij}(t,\x)=\boL_1(\wtP_{ij}(t,\x))$$
where $\boL_1$ is the operator introduced in Proposition \ref{prop-boL} with 
$\mu=1$. Explicitely,
$$\whP_{ij}(t,\x)(\lambda)=\frac{1}{\lambda-1}\left(\wtP_{ij}(t,\x)(\lambda)-\wtP_{ij}(t,\x)(1)\right).$$
Since $\boL_1$ is a bounded linear operator, $\whP_{ij}$ is a smooth map of
$t$, $t\log t$ and $\x$ with values in $\sl(2,\boW)$.
Using Notation \eqref{eq-x1xt} for the parameters $g_{ij}$ and
remembering that the parameters $p_{ij}$ and $p_{ji}$ are in $\C$, we obtain from Proposition \ref{prop-dwtPij}:
\begin{proposition}
\label{prop-dwhPij}
The differential of $\whP_{ij}$ is given by the following formula\begin{eqnarray*}
\lefteqn{d_{\x}\whP_{ij}(0,\x_0)=
\rho_{ij}\matrix{0&\frac{-\lambda^{-1}}{\overline{\pi_{ij}}}\\-\overline{\pi_{ij}}&0}dp_{ij}
+\rho_{ij}\matrix{0&\frac{-1}{\overline{\pi_{ij}}}\\-\lambda^{-1}\overline{\pi_{ij}}&0}dp_{ji}
+\matrix{0&\frac{-1}{\overline{\pi_{ij}}}\\\overline{\pi_{ij}}&0}dg_{ij}(1)}\\
&&
+\matrix{2&\frac{-(\lambda+1)}{\overline{\pi_{ij}}}\\(\lambda+1)\overline{\pi_{ij}}&-2}
d\wtg_{ij}
+\matrix{0&\frac{-1}{\overline{\pi_{ij}}}\\\overline{\pi_{ij}}&0}dh_{ij}
-2\rho_{ij}\matrix{0&\frac{1}{\overline{\pi_{ij}}}\\\overline{\pi_{ij}}&0}dm_{ij}
\end{eqnarray*}
\end{proposition}
\subsubsection{Reformulation of the Monodromy Problem}
We define for $(t,\x)$ in a neighborhood of $(0,\x_0)$:
$$\boE_{\Gamma ij,1}(t,\x)=
\lambda\whP_{ij,11}(t,\x)-\whP_{ij,11}^*(t,\x))$$
$$\boE_{\Gamma ij,2}(t,\x)=
\lambda\whP_{ij,12}(t,\x)-\whP_{ij,21}^*(t,\x).$$
\begin{proposition}
Problem \eqref{pb-monodromy-edges} is equivalent to
\begin{equation}
\label{pb-monodromy-edges2}
\left\{\begin{array}{ll}
\boE_{\Gamma ij,1}(t,\x)^+=0&\quad(i)\\
\boE_{\Gamma ij,2}(t,\x)=0&\quad(ii)\\
\wtP_{ij;11}(t,\x)(1)=0 & \quad(iii)\\
\whP_{ij;12}(t,\x)(1)=0&\quad (iv)\\
\Re\left[\whP_{ij;11}(t,\x)(1)\right]=\smallfrac{1}{2}(\ell_{ij}-2)&\quad(v)
\end{array}\right.
\end{equation}
\end{proposition}
Proof:
\begin{myenumerate}
\item We have by definition of $\wtP_{ij}$:
$$P_{ij}(t,\x)(1)=I_2\Leftrightarrow \wtP_{ij}(t,\x)(1)=0.$$
By Point 3 of Proposition \ref{prop-properties}:
$$P_{ij}(t,\x)(1)=\exp\left[\matrix{0&1\\0&0}\int_{\Gamma_{ij}}\beta_{t,\x}\right]$$
This gives:
$$\wtP_{ij}(t,\x)(1)=
H_{ij}(1)^{-1}\matrix{0&1\\0&0}H_{ij}(1)\int_{\Gamma_{ij}}\beta_{t,\x}
=\frac{1}{1+|\pi_{ij}|^2}\matrix{-\overline{\pi_{ij}}&1\\-\overline{\pi_{ij}}^2&\overline{\pi_{ij}}}\int_{\Gamma_{ij}}\beta_{t,\x}$$
Since $\pi_{ij}\neq 0$,
$$\wtP_{ij}(t,\x)(1)=0\Leftrightarrow\wtP_{ij;11}(t,\x)(1)=0.$$
So Item (ii) of Problem \eqref{pb-monodromy-edges} is equivalent to Item (iii) of Problem \eqref{pb-monodromy-edges2}.
\item Assume from now on that Item (ii) of Problem \eqref{pb-monodromy-edges}
is satisfied. Then:
$$\wtP_{ij}(t,\x)=(\lambda-1)\whP_{ij}(t,\x)$$
so for $k,\ell\in\{1,2\}$,
\begin{eqnarray*}
\wtP_{ij;k\ell}(t,\x)+\wtP_{ij;\ell k}(t,\x)^*&=&(\lambda-1)\whP_{ij;k\ell}(t,\x)+(\lambda^{-1}-1)\whP_{ij;\ell k}(t,\x)^*\\
&=&(1-\lambda^{-1})\left(\lambda\whP_{ij;k\ell}(t,\x)-\whP_{ij;\ell k}(t,\x)^*\right).\end{eqnarray*}
Hence
\begin{eqnarray*}
P_{ij}(t,\x)\in\Lambda\su(2)&\Leftrightarrow&
\left\{\begin{array}{l}
\wtP_{ij;11}(t,\x)+\wtP_{ij;11}(t,\x)^*=0\\
\wtP_{ij;12}(t,\x)+\wtP_{ij;21}(t,\x)^*=0\end{array}\right.\\
&\Leftrightarrow&
\left\{\begin{array}{l}
\lambda\whP_{ij;11}(t,\x)-\whP_{ij;11}(t,\x)^*=0\\
\lambda\whP_{ij;12}(t,\x)-\whP_{ij;21}(t,\x)^*=0\end{array}\right.\\
&\Leftrightarrow&
\left\{\begin{array}{l}
\boE_{\Gamma ij,1}(t,\x)=0\\
\boE_{\Gamma ij,2}(t,\x)=0.\end{array}\right.
\end{eqnarray*}
We have
$$\lambda\boE_{ij,1}(t,\x)^*=\lambda\left(\lambda^{-1}\whP_{ij;11}(t,\x)^*-\whP_{ij;11}(t,\x)\right)=-\boE_{ij,1}(t,\x).$$
Hence
$$\boE_{ij,1}(t,\x)^+=-\lambda\left(\boE_{ij,1}(t,\x)^*\right)^\pos=-\lambda\left(\boE_{ij,1}(t,\x)^\neg\right)^*.$$
So
$$\boE_{ij,1}(t,\x)=0\Leftrightarrow \boE_{ij,1}(t,\x)^+=0.$$
Hence Item (i) of Problem \eqref{pb-monodromy-edges} is equivalent to Items (i) and (ii) of Problem \eqref{pb-monodromy-edges2}.
\item Assume from now on that Items (i) and (ii) of Problem \eqref{pb-monodromy-edges}
are satisfied. Then
$$\whP_{ij}(t,\x)(1)=\frac{\partial}{\partial\lambda} \wtP_{ij}(t,\x)(1)=
H_{ij}(1)^{-1}\frac{\partial}{\partial\lambda} P_{ij}(t,\x)(1)H_{ij}(1)+\matrix{-1&0\\0&1}.$$
On the other hand, recalling that $\ell_{ij}$ is the length of the edge $e_{ij}$, we have
$$\v_j-\v_i=\ell_{ij}\u_{ij}=-\ell_{ij}\NS(\pi_{ij}).$$
By Equation \eqref{eq-normal},
$$V_j-V_i=\frac{\ii}{2}\ell_{ij}\FS(\pi_{ij},1)\matrix{1&0\\0&-1}\FS(\pi_{ij},1)^{-1}
=\frac{\ii}{2}\ell_{ij}H_{ij}(1)\matrix{1&0\\0&-1}H_{ij}(1)^{-1}.$$
So
\begin{eqnarray*}
\frac{\partial}{\partial\lambda}P_{ij}(t,\x)(1)=-\ii(V_j-V_i)
&\Leftrightarrow&
H_{ij}(1)^{-1}\frac{\partial}{\partial\lambda}P_{ij}(t,\x)(1)H_{ij}(1)=
\smallfrac{1}{2}\ell_{ij}\matrix{1&0\\0&-1}\\
&\Leftrightarrow&
\whP_{ij}(t,\x)(1)=\smallfrac{1}{2}(\ell_{ij}-1)\matrix{1&0\\0&-1}.
\end{eqnarray*}
Since $\whP_{ij}(t,\x)(1)\in\Lambda\su(2)$,
Item (iii) of Problem \eqref{pb-monodromy-edges} is equivalent to Items (iv) and (v) of Problem \eqref{pb-monodromy-edges2}.
\cqfd
\end{myenumerate}
\subsubsection{Solving the Monodromy Problem for $t\neq 0$}
\label{section-monodromy-edges-implicit}
At the central value, we have
$\wtP_{ij}(0,\x_0)=\whP_{ij}(0,\x_0)=0$, so Items (i) to (iv)
of Problem \eqref{pb-monodromy-edges2} are satisfied, and
Item (v) is equivalent to $\ell_{ij}=2$.
So we leave aside Item (v) which will be solved using the non-degeneracy
hypothesis in Section \ref{section-balancing}.
We define for $(t,\x)$ in a neighborhood of $(0,\x_0)$:
$$\boE_{\Gamma ij,3}(t,\x)=\wtP_{ij,11}(t,\x)(1)\in\C$$
$$\boE_{\Gamma ij,4}(t,\x)=\whP_{ij,12}(t,\x)(1)\in\C.$$
\begin{proposition}
\label{prop-edges-isomorphism}
Let $L$ be the partial differential of
$$\left(\lambda^{-1}\boE_{\Gamma ij,1}^+,\lambda^{-1}\boE_{\Gamma ij,2}^+,
(\boE_{\Gamma ij,2}^\neg)^*,\boE_{\Gamma ij,3},\boE_{\Gamma ij,4}\right)$$
with respect to
$$\left(\wtg_{ij},h_{ij},m_{ij},g_{ij}(1),p_{ji}\right)$$
at $(0,\x_0)$. Then $L$ is an automorphism of $(\boW^\pos)^3\times\C^2$.
\end{proposition}
Proof: using Proposition \ref{prop-dwtPij}, we obtain:
\begin{equation}
\label{eq-edges-1}
d_{\x}\boE_{\Gamma ij,3}(0,\x_0)=\rho_{ij}(dp_{ji}-dp_{ij})+2dg_{ij}(1).
\end{equation}
Using Proposition \ref{prop-dwhPij}, we obtain:
$$d_{\x}\boE_{\Gamma ij,1}(0,\x_0)=2\lambda d\wtg_{ij}-2d\wtg_{ij}^*$$
\begin{equation}
\label{eq-edges-2}
\lambda^{-1}d_{\x}\boE_{\Gamma ij,1}(0,\x_0)^+=2d\wtg_{ij}
\end{equation}
\begin{eqnarray*}
d_{\x}\boE_{\Gamma ij,2}(0,\x_0)&=&
\frac{-\lambda}{\overline{\pi_{ij}}}\left[\rho_{ij}(\lambda^{-1}dp_{ij}+dp_{ji}+2dm_{ij})+dg_{ij}(1)+(\lambda+1)d\wtg_{ij}+dh_{ij}\right]\\
&&+\pi_{ij}\left[\overline{\rho_{ij}}(d\overline{p_{ij}}+\lambda d\overline{p_{ji}}+2dm_{ij}^*)-d\overline{g_{ij}(1)}-(\lambda^{-1}+1)d\wtg_{ij}^*-dh_{ij}^*\right]
\end{eqnarray*}
\begin{equation}
\label{eq-edges-3}
\lambda^{-1}d_{\x}\boE_{\Gamma ij,2}(0,\x_0)^+=
\frac{-1}{\overline{\pi_{ij}}}\left[\rho_{ij}(dp_{ji}+2dm_{ij})+dg_{ij}(1)+(\lambda+1)d\wtg_{ij}+dh_{ij}\right]+\pi_{ij}\overline{\rho_{ij}}d\overline{p_{ji}}
\end{equation}
\begin{equation}
\label{eq-edges-4}
\left(d_{\x}\boE_{\Gamma ij,2}(0,\x_0)^\neg\right)^*=
\frac{-\overline{\rho_{ij}}}{\pi_{ij}}d\overline{p_{ij}}+\overline{\pi_{ij}}\left[\rho_{ij}(dp_{ij}+2dm_{ij})-dg_{ij}(1)
-(\lambda+1)d\wtg_{ij}-dh_{ij}\right]
\end{equation}
\begin{equation}
\label{eq-edges-5}
d_{\x}\boE_{\Gamma ij,4}(0,\x_0)=\frac{-1}{\overline{\pi_{ij}}}\left[
\rho_{ij}(dp_{ij}+dp_{ji}+2dm_{ij}(1))+dg_{ij}(1)+2d\wtg_{ij}(1)+dh_{ij}(1)\right].
\end{equation}
By Equations \eqref{eq-edges-1}, \eqref{eq-edges-3} and \eqref{eq-edges-4},
the partial differential of
$(\lambda^{-1}\boE_{\Gamma ij,1}^+,\lambda^{-1}\boE_{\Gamma ij,2}^+,
(\boE_{\Gamma ij,2}^\neg)^*)$
with respect to
$(\wtg_{ij},h_{ij},m_{ij})$ is a matrix-type operator (see Definition \ref{def-matrixtype}) with
matrix
$$\left(\begin{array}{ccc}2&0&0\\\frac{-\lambda-1}{\overline{\pi_{ij}}}&
\frac{-1}{\overline{\pi_{ij}}}&\frac{-2\rho_{ij}}{\overline{\pi_{ij}}}\\
-(\lambda+1)\overline{\pi_{ij}}&-\overline{\pi_{ij}}&2\rho_{ij}\overline{\pi_{ij}}
\end{array}\right).$$
This matrix has constant determinant $-8\rho_{ij}$ so is invertible in
$\boM_3(\boW^\pos)$.
By Proposition \ref{prop-fredholm}, it suffices to prove that $L$ is injective.
Let us solve formally the system $L=0$
to express the differential of all parameters in function of $dp_{ij}$.
Equation \eqref{eq-edges-2} gives
\begin{equation}
\label{eq-edges-6}
d\wtg_{ij}=0.
\end{equation}
Equation \eqref{eq-edges-3} evaluated at $\lambda=1$ substracted from Equation \eqref{eq-edges-5} gives
$$\frac{-\rho_{ij}}{\overline{\pi_{ij}}}dp_{ij}-\pi_{ij}\overline{\rho_{ij}}d\overline{p_{ji}}=0$$
from which we obtain
\begin{equation}
\label{eq-dpji}
dp_{ji}=\frac{-1}{\overline{\pi_{ij}}^2}d\overline{p_{ij}}.
\end{equation}
Equations \eqref{eq-edges-1}, \eqref{eq-edges-6} and \eqref{eq-dpji} give
\begin{equation}
\label{eq-dgij}
dg_{ij}=dg_{ij}(1)=\frac{\rho_{ij}}{2}\left(dp_{ij}+\frac{1}{\overline{\pi_{ij}}^2}d\overline{p_{ij}}\right).
\end{equation}
By substitution of Equations \eqref{eq-edges-6}, \eqref{eq-dpji} and \eqref{eq-dgij} in \eqref{eq-edges-3} and
\eqref{eq-edges-4}, we obtain after simplification the system
$$\left\{\begin{array}{ll}
dh_{ij}+2\rho_{ij}dm_{ij}=\frac{-3\rho_{ij}}{2}dp_{ij}+\frac{\rho_{ij}}{2\overline{\pi_{ij}}^2}d\overline{p_{ij}}\\
-dh_{ij}+2\rho_{ij}dm_{ij}=\frac{-\rho_{ij}}{2}dp_{ij}+\frac{3\rho_{ij}}{2\overline{\pi_{ij}}^2}d\overline{p_{ij}}
\end{array}\right.$$
whose solution is
\begin{equation}
\label{eq-dhij}
dh_{ij}=\frac{-\rho_{ij}}{2}\left(dp_{ij}+\frac{1}{\overline{\pi_{ij}}^2}d\overline{p_{ij}}\right).
\end{equation}
\begin{equation}
\label{eq-dmij}
dm_{ij}=\frac{1}{2}\left(-dp_{ij}+\frac{1}{\overline{\pi_{ij}}^2}d\overline{p_{ij}}\right)
\end{equation}

Setting $dp_{ij}=0$, we obtain that $L$ is injective, so is an automorphism by Proposition \ref{prop-fredholm}.
\cqfd

\medskip

We now prove Proposition \ref{prop-monodromy-edges}.
We decompose $\x=(\x',\x'',\x''')$ where $\x'$ is the vector of the parameters
which have already been determined in Propositions \ref{prop-regularity},
\ref{prop-monodromy-nodes} and \ref{prop-monodromy-ends},
$$\x''=\left(\wtg_{ij},h_{ij},m_{ij},g_{ij}(1),p_{ji}\right)_{(i,j)\in I}$$
and $\x'''$ denotes the remaining parameters.
Then $\x'$ is a smooth function of $t$, $\x''$ and $\x'''$
so we write $\x'=\x'(t,\x'',\x''')$.
Define
$$\boE(t,\x)=\left(\lambda^{-1}\boE_{\Gamma ij,1}^+,\lambda^{-1}\boE_{\Gamma ij,2}^+,
(\boE_{\Gamma ij,2}^\neg)^*,\boE_{\Gamma ij,3},\boE_{\Gamma ij,4}\right)_{(i,j)\in I}(t,\x)$$
$$\boF(t,\x'',\x''')=\boE(t,\x'(t,\x'',\x'''),\x'',\x''')$$
Since $\boE(0,\x)$ does not depend on $\x'$, we have $d_{\x'}\boE(0,\x_0)=0$ so by the chain rule,
$$d_{\x''}\boF(0,\x''_0,\x'''_0)=d_{\x''}\boE(0,\x_0).$$
By Proposition \ref{prop-edges-isomorphism}, $d_{\x''}\boE(0,\x_0)$ is an automorphism (it has block diagonal
form).
By Proposition \ref{prop-Pij}, $\boF(t,\x'',\x''')$ is a smooth function of $t$ and $t\log t$,
so we can write
$$\boF(t,\x'',\x''')=\boG(t,t\log t,\x'',\x''')$$
where $\boG$ is a smooth function of all its arguments.
By the Implicit Function Theorem, for $(t,s,\x''')$ in a neighborhood of $(0,0,\x'''_0)$, there
exists a unique $\y''$, depending smoothly on $(t,s,\x''')$ such that
$\boG(t,s,\y''(t,s,\x'''),\x''')=0$.
We define $\x''(t,\x''')=\y''(t,t\log t,\x''')$.
Then $\boF(t,\x''(t,\x'''),\x''')=0$, so Problem \eqref{pb-monodromy-edges} is solved.
Moreover, $\x''$ is a smooth function of $t$, $t\log t$ and $\x'''$.
\cqfd
\begin{remark}
The parameters $g_{ij}$, $h_{ij}$, $m_{ij}$ and $p_{ji}$ are now
function of $(t,\x''')$. The differential of these parameters with respect to
$\x'''$ at $(0,\x_0)$ are obtained by solving $L=0$ so are given by Equations \eqref{eq-dgij},
\eqref{eq-dhij}, \eqref{eq-dmij} and \eqref{eq-dpji}.
In particular, they only depend on $dp_{ij}$.
\end{remark}
\begin{remark} As we have seen, when solving equations depending smoothly
on $t$ and $t\log t$, we apply the Implicity Function Theorem with the variables
$t$ and $s$ and then we substitute $s=t\log t$, so the solution depends 
smoothly on $t$ and $t\log t$. This remark applies to all our subsequent uses
of the Implicit Function Theorem.
\end{remark}
\section{The Regularity Problem at $\infty_{ij}$}
\label{section-regularity2}
In this section, we consider again the potential $\whxi_{ij,t,\x}=\xi_{t,\x}\cdot G_{ij,\x}^{-1}$
introduced in Section \ref{section-regularity}.
We have $\whxi_{ij,0,\x}=\xiC$ in $\CC_{ij}$ and by Proposition \ref{prop-order}, $\whxi_{ij,t,\x}$ has a pole of multiplicity at
most $\minimatrix{1&1\\2&1}$ at $\infty_{ij}$.
In this section, we solve the following problem for $(i,j)\in I$:
\begin{equation}
\label{pb-regularity2}
\left\{\begin{array}{l}
\Res_{\infty_{ij}}\whxi_{ij,t,\x;12}^{(-1)}=0\\
\Re\left[\Res_{\infty_{ij}}(z\,\whxi_{ij,t,\x;12}^{(-1)})\right]=0
\end{array}\right.
\end{equation}
If the Monodromy Problem is solved and Problem \eqref{pb-regularity2}
is solved, Theorem \ref{thm-regularity-infty} in Appendix \ref{appendix-regularity}
tells us that $\infty_{ij}$ is a removable singularity.
Note that Problem \eqref{pb-regularity2} amounts to only three real equations:
this is what remains of the Regularity Problem at $\infty_{ij}$
when the Monodromy Problem is solved.
\medskip

Assume that the parameters are as in Proposition \ref{prop-monodromy-edges}.
The only remaining parameters are $p_{ij}\in\C$ for $(i,j)\in I$
and $\Re(\rho_{ij}^2r_{ji})\in\R$ for $(i,j)\in I\cup I^*$.
We fix the following normalisation:
\begin{equation}
\label{eq-normalisation-rji}
\Re(\rho_{ij}^2 r_{ij})=-\tau_{ij}\quad\mbox{ for $(i,j)\in I^*$}.
\end{equation}
\begin{proposition}
\label{prop-regularity2}
For $t$ small enough, there exists unique values of the parameters
$$\left(p_{ij},\Re(\rho_{ij}^2 r_{ij})\right)_{(i,j)\in I}\in(\C\times\R)^I$$
depending smoothly on $t$ and
$t\log t$, such that Problem \eqref{pb-regularity2} is solved for $(i,j)\in I$.
\end{proposition}
Proof: we first compute the residues in Problem \eqref{pb-regularity2}:
\begin{claim}
\label{claim-regularity2}
 We have for all $(t,\x)$ in a neighborhood of $(0,\x_0)$:
$$\Res_{\infty_{ij}}\whxi_{ij,t,\x;12}^{(-1)}=t\,c_{mij}^0$$
$$\Res_{\infty_{ij}}(z\whxi_{ij,t,\x;12}^{(-1)})=t\,b_{mij}^0.$$
\end{claim}
Proof: as in the proof of Proposition \ref{prop-order}, we write $\Theta_{t,\x}=
\Theta_{0,\x}+\Xi_{t,\x}$ where $\Xi_{t,\x}$ is holomorphic in a neighborhood of $\infty_{ij}$.
Then
$$\xi_{t,\x;12}^{(-1)}=\zeta_{t,\x;12}^{(-1)}-t\left(\chi_{t,\x;12}^{(-1)}+\Theta_{0,\x;12}^{(-1)}+\Xi_{t,\x;12}^{(-1)}\right).$$
By Proposition \ref{prop-gauging} with $G=G_{ij,\x}$:
$$\Theta_{0,\x;12}^{(-1)}=(G_{ij,\x;22}^0)^2\left(\frac{b_{mij}^0\,dz}{(z-m_{ij}^0)^2}+\frac{c_{mij}^0\,dz}{z-m_{ij}^0}\right)$$
By Proposition \ref{prop-gauging} with by $G=G_{ij,\x}^{-1}$:
$$\whxi_{ij,t,\x;12}^{(-1)}=(G_{ij,\x;11}^0)^2\xi_{t,\x;12}^{(-1)}
=\frac{1}{(2g_{ij}^0(z-m_{ij}^0))^2}\left[
\zeta_{t,\x;12}^{(-1)}-t\,\chi_{t,\x;12}^{(-1)}-t\,\Xi_{t,\x;12}^{(-1)}\right]
-\frac{t\,b_{mij}^0\,dz}{(z-m_{ij}^0)^2}-\frac{t\,c_{mij}^0\,dz}{z-m_{ij}^0}.$$
The bracket has at most a simple pole at $\infty_{ij}$ so it does not contribute to the residues (thanks to the factor $(z-m_{ij}^0)^{-2}$ in front). Claim \ref{claim-regularity2} follows.\cqfd

\medskip
As a consequence of Claim \ref{claim-regularity2},
Problem \eqref{pb-regularity2} is equivalent for $t\neq 0$ to
\begin{equation}
\label{pb-bmij0}
\left\{\begin{array}{l}
\Re(b_{mij}^0)=0\\
c_{mij}^0=0
\end{array}\right.
\end{equation}
Assume that the parameters are as in Proposition \ref{prop-monodromy-edges}
and let $\x'$ be the vector of the remaining parameters, namely
$p_{ij}$ and $\Re(\rho_{ij}^2r_{ij})$ for $(i,j)\in I$.
Recall that we computed $b_{mij}^0$ and $c_{mij}^0$ in function of the other
parameters in Proposition \ref{prop-bmij0}.
Since then, some of the parameters involved in this formula have been determined
as functions of $\x'$.
We now compute explicitely $b_{mij}$ and $c_{mij}$ as functions of
$\x'$ at $t=0$.
Assume that $t=0$.
By Equation \eqref{eq-xiC.Gij}, we have:
$$A_{ij}^0=\frac{2g_{ij}^0-h_{ij}^0}{2g_{ij}^0(q_{ij}-m_{ij}^0)},\qquad
A_{ji}^0=\frac{2g_{ij}^0-h_{ij}^0}{2g_{ij}^0(q_{ji}-m_{ij}^0)}$$
\begin{equation}
\label{eq-Cij0}C_{ij}^0=\frac{-1}{4(g_{ij}^0)^2(q_{ij}-m_{ij}^0)^2},\qquad
C_{ji}^0=\frac{-1}{4(g_{ij}^0)^2(q_{ji}-m_{ij}^0)^2}.
\end{equation}
By Proposition \ref{prop-monodromy-nodes-t0}:
$$b_{ij}^0=\frac{r_{ij}(2g_{ij}^0-h_{ij}^0)}{g_{ij}^0(m_{ij}^0-q_{ij})},\qquad
b_{ji}^0=\frac{r_{ji}(2g_{ij}^0-h_{ij}^0)}{g_{ij}^0(m_{ij}^0-q_{ji})}.$$
Using Proposition \ref{prop-bmij0}, we obtain after simplification:
\begin{equation}
\label{eq-bmij00}
b_{mij}^0(0,\x')=\frac{(g_{ij}^0-h_{ij}^0)}{4(g_{ij}^0)^3}\left(\frac{r_{ij}}{(m_{ij}^0-q_{ij})^2}+\frac{r_{ji}}{(m_{ij}^0-q_{ji})^2}\right)
\end{equation}
\begin{equation}
\label{eq-cmij00}
c_{mij}^0(0,\x')=\frac{h_{ij}^0}{4(g_{ij}^0)^3}\left(\frac{r_{ij}}{(m_{ij}^0-q_{ij})^3}+\frac{r_{ji}}{(m_{ij}^0-q_{ji})^3}\right).
\end{equation}
Equations \eqref{eq-bmij00} and \eqref{eq-cmij0} imply that
$b_{0ij}(0,\x'_0)^0=0$ and $c_{0ij}(0,\x'_0)^0=0$ so Problem \eqref{pb-bmij0} is
solved at the central value.
Proposition \ref{prop-regularity2} follows from the following claim and the Implicit
Function Theorem:
\begin{claim}
\label{claim-regularity2-isomorphism}
The partial differential of
$(\Re(b_{mij}^0),c_{mij}^0)$
with respect to $(\Re(\rho_{ij}^2r_{ij}),p_{ij})$ at $(0,\x'_0)$ is an automorphism of $\R\times\C$.
\end{claim}
Proof:
by Proposition \ref{prop-monodromy-nodes-t0} and Equation \eqref{eq-Cij0}, we have at $t=0$:
$$\Im\left(\frac{r_{ij}}{(g_{ij}^0)^2(q_{ij}-m_{ij}^0)^2}\right)=
\Im\left(\frac{r_{ji}}{(g_{ij}^0)^2(q_{ji}-m_{ij}^0)^2}\right)=0.$$
By differentiation with respect to $\x'$ at $(0,\x'_0)$, this gives
\begin{equation}
\label{eq-Imrho2drij}
\Im(\rho_{ij}^2 dr_{ij})=4\Im(\rho_{ij}\tau_{ij}dm_{ij}^0),\qquad
\Im(\rho_{ji}^2 dr_{ji})=4\Im(\rho_{ji}\tau_{ji}dm_{ij}^0).
\end{equation}
By differentiation of Equations \eqref{eq-bmij00} and \eqref{eq-cmij00} with respect to
$\x'$ at $(0,\x'_0)$, we obtain
$$db_{mij}^0=\frac{1}{4}(dg_{ij}^0-dh_{ij}^0)(-4\tau_{ij}-4\tau_{ji})=-2\tau_{ij}(dg_{ij}^0-dh_{ij}^0)$$
$$dc_{mij}^0=\frac{1}{4}\left(
-8\rho_{ij}^3dr_{ij}+12\rho_{ij}^2\tau_{ij}dm_{ij}^0
-8\rho_{ji}^3dr_{ji}+12\rho_{ji}^2\tau_{ji}dm_{ij}^0\right)
=-2\rho_{ij}^3(dr_{ij}-dr_{ji})+24\rho_{ij}^2\tau_{ij}dm_{ij}^0.$$
Recalling Normalisation \eqref{eq-normalisation-rji}:
$$\Re\left(\rho_{ij}^{-1}c_{mij}^0\right)=-2\Re(\rho_{ij}^2dr_{ij})+24\tau_{ij}\Re(\rho_{ij}dm_{ij}^0).$$
Using Equation \eqref{eq-Imrho2drij}:
$$\Im\left(\rho_{ij}^{-1}dc_{mij}^0\right)=
-8\Im(\rho_{ij}\tau_{ij}dm_{ij}^0)
+8\Im(\rho_{ji}\tau_{ji}dm_{ij}^0)
+24\Im(\rho_{ij}\tau_{ij}dm_{ij}^0)=8\tau_{ij}\Im(\rho_{ij}dm_{ij}^0).$$
Using Equations \eqref{eq-dgij}, \eqref{eq-dhij} and \eqref{eq-dmij}, we finally obtain:
$$\Re\left(db_{mij}^0\right)=-2\tau_{ij}\Re\left(\rho_{ij}dp_{ij}+\frac{\rho_{ij}}{\overline{\pi_{ij}}^2}d\overline{p_{ij}}\right)$$
\begin{equation}
\label{eq-regularity2-1}
\Re\left(\rho_{ij}^{-1}dc_{mij}^0\right)=-2\Re(\rho_{ij}^2dr_{ij})+12\tau_{ij}\Re\left(-\rho_{ij}dp_{ij}+\frac{\rho_{ij}}{\overline{\pi_{ij}}^2}d\overline{p_{ij}}\right)
\end{equation}
$$\Im\left(\rho_{ij}^{-1}dc_{mij}^0\right)=4\tau_{ij}\Im\left(-\rho_{ij}dp_{ij}+\frac{\rho_{ij}}{\overline{\pi_{ij}}^2}d\overline{p_{ij}}\right).$$
This implies
\begin{equation}
\label{eq-regularity2-2}
\Re(db_{mij}^0)+\frac{\ii}{2}\Im(\rho_{ij}^{-1}dc_{mij}^0)
=-2\tau_{ij}\left(\rho_{ij}dp_{ij}+\frac{\overline{\rho_{ij}}}{\pi_{ij}^2}dp_{ij}\right)
=-2\frac{\tau_{ij}}{\pi_{ij}}dp_{ij}.
\end{equation}
Claim \eqref{claim-regularity2-isomorphism} follows from Equations
\eqref{eq-regularity2-1} and \eqref{eq-regularity2-2}.
\cqfd
\section{Using the balancing and non-degeneracy hypothesis}
\label{section-balancing}
At this point, all parameters have been determined as smooth functions of 
$t$ and $t\log t$. We write $\x=\x(t)$ and we have $\x(0)=\x_0$.
The central value $\x_0$, as given in Section \ref{section-setup-parameters}, depends on the graph $\Gamma$.
In this section, we use the balancing and non-degeneracy hypothesis to solve the following
problem:
\begin{equation}
\label{pb-balancing}
\left\{\begin{array}{ll}
\forall i\in[1,N],\quad \Res_{\infty_i}\xi_{t,\x(t);21}^0=0&\quad (i)\\
\forall i\in[1,N],\quad\Re\left[\Res_{\infty_i}(z\,\xi_{t,\x(t);21}^0)\right]=0&\quad(ii)\\
\forall (i,j)\in I,\quad\Re\left(\whP_{ij;11}(t,\x(t))(1)\right)=\smallfrac{1}{2}(\ell_{ij}-2)
&\quad(iii).
\end{array}\right.\end{equation}
Item (iii) of Problem \eqref{pb-balancing} is Item (v) of Problem \eqref{pb-monodromy-edges2} which we still have to solve.
We have $\xi_{0,\x}=\xiS$ in $\CC_i$ and the potential $\xi_{t,\x}$ has a pole of multiplicity at most $\minimatrix{1&2\\1&1}$ at $\infty_i$.
Provided the Monodromy Problem is solved and Items (i) and
(ii) of Problem \eqref{pb-balancing} are solved, Corollary \ref{cor-regularity-infty}
in Appendix \ref{appendix-regularity} tells us that $\infty_i$ is a removable singularity.
As in Section \ref{section-regularity2}, Items (i) and (ii) are only three real equations: this is what remains of
the Regularity Problem at $\infty_i$ when the Monodromy Problem is solved.
\begin{proposition}
\label{prop-balancing}
If the graph $\Gamma$ has length-2 edges and is balanced and non-degenerate,
then for $t>0$ small enough, there
exists a deformation $\Gamma(t)$ of $\Gamma$, depending smoothly on $t$ and
$t\log t$, such that Problem \eqref{pb-balancing} is solved.
\end{proposition}
Proof: let $\gamma_{t,\x}=\frac{1}{t}\xi_{t,\x;21}$.
At $t=0$, we have $\xi_{0,\x}=\xiS$ in $\Omega_i$ so the restriction of $\gamma_{t,\x}$
to $\Omega_i$ extends holomorphically at $t=0$.
We define for $i\in[1,N]$:
$$\boE_{i,1}(t)=\Res_{\infty_i}\gamma_{t,\x(t)}^0$$
$$\boE_{i,2}(t)=\Re\left[\Res_{\infty_i}(z\,\gamma_{t,\x(t)}^0)\right]$$
and for $(i,j)\in I$:
$$\boE_{ij,3}(t)=\Re\left[\whP_{ij;11}(t,\x(t))(1)\right]+\smallfrac{1}{2}(2-\ell_{ij}).$$
Problem \eqref{pb-balancing} is equivalent for $t\neq 0$ to the following problem:
\begin{equation}
\label{pb-balancing2}
\left\{\begin{array}{ll}
\forall i\in[1,N],\quad&\boE_{i,1}(t)=0\\
\forall i\in[1,N],\quad&\boE_{i,2}(t)=0\\
\forall (i,j)\in I,\quad&\boE_{ij,3}(t)=0
\end{array}\right.\end{equation}
We have already seen in Section \ref{section-monodromy-edges-implicit}
that
\begin{equation}
\label{eq-balancing-1}
\boE_{ij,3}(0)=\smallfrac{1}{2}(2-\ell_{ij}).
\end{equation}
We have
$$\gamma_{0,\x_0}=\frac{\partial}{\partial t}\xi_{t,\x_0;21}|_{t=0}
=\frac{\partial}{\partial t}\zeta_{t,\x_0;21}|_{t=0}+(\lambda-1)\chi_{0,\x_0;21}.$$
By Theorem \ref{thm-openingnodes-derivee}:
$$\frac{\partial}{\partial t}\zeta_{t,\x;21}(z,\lambda)|_{t=0}
=-\sum_{j\in E_i}\frac{r_{ij}\,dz}{(z-p_{ij})^2}\zeta_{0,\x;21}(q_{ij},\lambda)).$$
Hence at the central value, using Equations \eqref{eq-xiC.Gij0} and
\eqref{eq-2qij-2muij}:
$$\frac{\partial}{\partial t}\zeta_{t,\x_0;21}(z,\lambda)|_{t=0}
=(1-\lambda^2)\sum_{j\in E_i}\frac{\tau_{ij}dz}{(z-\pi_{ij})^2}.$$
By definition of $\chi_{0,\x}$ and using the central value of the parameters as indicated
in Section \ref{section-setup-parameters}:
$$\chi_{0,\x_0;21}=
(\lambda-1)\sum_{j\in E_i}\frac{-2\rho_{ij}\tau_{ij}\,dz}{z-\pi_{ij}}
+(\lambda-1)\sum_{k\in R_i}\left(\frac{\tau_k\,dz}{(z-\pi_k)^2}-\frac{2\rho_k\tau_k\,dz}{z-\pi_k}\right).$$
Collecting all terms and taking $\lambda=0$, we obtain:
$$\gamma_{0,\x_0}^0=\sum_{j\in E_i}\left(\frac{\tau_{ij}dz}{(z-\pi_{ij})^2}
-\frac{2\rho_{ij}\tau_{ij}dz}{z-\pi_{ij}}\right)+\sum_{k\in R_i}\left(\frac{\tau_kdz}{(z-\pi_k)^2}
-\frac{2\rho_k\tau_kdz}{z-\pi_k}\right).$$
This gives
\begin{equation}
\label{eq-balancing-2}
\boE_{i,1}(0)=\sum_{j\in E_i}2\tau_{ij}\frac{\overline{\pi_{ij}}}{1+|\pi_{ij}|^2}
+\sum_{k\in R_i}2\tau_k\frac{\overline{\pi_k}}{1+|p_k|^2}=\overline{\mbox{Horiz}(\force_i)}
\end{equation}
\begin{equation}
\label{eq-balancing-3}
\boE_{i,2}(0)=\sum_{j\in E_i}\tau_{ij}\frac{|\pi_{ij}|^2-1}{|\pi_{ij}|^2+1}+
\sum_{k\in R_i}\tau_k\frac{|\pi_k|^2-1}{|\pi_k|^2+1}=\mbox{Vert}(\force_i)
\end{equation}
where $\mbox{Horiz}(\force_i)$ and $\mbox{Vert}(\force_i)$ are the horizontal and vertical
components of the force defined in Section \ref{section-intro-forces}.
By Equations \eqref{eq-balancing-1}, \eqref{eq-balancing-2} and \eqref{eq-balancing-3}, Problem \eqref{pb-balancing2} at $t=0$ is equivalent to the fact that
the graph $\Gamma$ has length-2 edges and is balanced.
Provided $\Gamma$ is non-degenerate, Problem \eqref{pb-balancing2} can be solved
for $t\neq 0$ using the Implicit Function Theorem.\cqfd
\section{Geometry of the surface}
\label{section-geometry}
At this point, all parameters have been determined as smooth functions of
$t$ and $t\log t$: we write $\x(t)$ for the value of the parameter vector $\x$
and $\xi_t=\xi_{t,\x(t)}$ for the potential at time $t\in[0,\epsilon)$.
Note that by Proposition \ref{prop-balancing}, the graph $\Gamma$ now depends
on $t$, so the numbers $\pi_{ij}(t)$, $\tau_{ij}(t)$ for $(i,j)\in I\cup I^*$ and
$\pi_k(t)$, $\tau_k(t)$ for $k\in[1,n]$ are now functions of $t$.
We adopt the following convention: if the variable $t$ is not written, it means
the value at $t=0$, corresponding to the given balanced and non-degenerate graph
$\Gamma$, so for example, $\pi_k$ means $\pi_k(0)$.
\medskip

We use the notation $C$ for uniform constants, depending only on the graph $\Gamma$. With an argument, $C(r)$ denotes a constant depending only on $r$ and $\Gamma$.
The same letter $C$ may be used to denote different constants.
We use the notation $||\cdot||_{\boW}$ for the functional norm introduced in Section
\ref{section-background-functional} and if $\Omega$ is a compact domain, $||\cdot||_{C^k(\Omega)}$ denotes the standard norm on $C^k(\Omega)$.
\subsection{Differentiability of the potential with respect to $t$}
\label{section-geometry-differentiability}
We define the following domains for $0<r<\frac{\varepsilon}{2}$:
$$\Omega_{i,r}=\{z\in\C_i:\forall j\in E_i,\;|z-\pi_{ij}|>r
\;\mbox{ and }\;
\forall k\in R_i,\;|z-\pi_k|>r
\}\quad\mbox{ for $i\in[1,N]$}$$
$$\Omega_{ij,r}=\{z\in\C_{ij}:
|z-q_{ij}|>r\;\mbox{ and }\;
|z-q_{ji}|>r\} \quad\mbox{ for $(i,j)\in I$}.$$
We denote $\whxi_{ij,t}$ the potential introduced in Section \ref{section-regularity}
at time $t$:
$$\whxi_{ij,t}=\xi_{t,\x(t)}\cdot G_{ij,\x(t)}^{-1}.$$
By Proposition \ref{prop-regularity}, $\whxi_{ij,t}$ extends holomorphically to $\Omega_{ij,r}$.
Note that the map $t\mapsto \x(t)$ is not differentiable at $t=0$. However:
\begin{proposition}
\label{prop-differentiability}
\begin{myenumerate}
\item For $i\in[1,N]$, the restriction of $\xi_t$ to $\Omega_{i,r}$ extends to a $C^1$ function of
$(t,z)\in(-\epsilon,\epsilon)\times\Omega_{i,r}$, with values in $\sl(2,\boW)$.
Moreover, $\xi_0=\xiS$ in $\Omega_{i,r}$.
\item For $(i,j)\in I$, the restriction of $\whxi_{ij,t}$ to $\Omega_{ij,r}$ extends to
a $C^1$ function of $(t,z)\in(-\epsilon,\epsilon)\times\Omega_{ij,r}$, with values
in $\sl(2,\boW)$. Moreover, $\whxi_{ij,0}=\xiC$ in $\Omega_{ij,r}$.
\end{myenumerate}
\end{proposition}
Proof: by Definition \ref{def-tlogt}, there exists a smooth function $\x(t,s)$, defined
for $(t,s)$ in a neighborhood of $(0,0)$, such that $\x(t)=\x(t,t\log t)$.
We denote $\xi_{t,s}=\xi_{t,\x(t,s)}$ and $\whxi_{ij,t,s}=\xi_{t,\x(t,s)}\cdot G_{ij,\x(t,s)}^{-1}$.
We have, at $t=0$ and for all $s$:
$$\xi_{0,s}=\xiS\quad\mbox{ in $\Omega_{i,r}$}$$
$$\whxi_{ij,0,s}=\xiC\cdot G_{ij,\x(0,s)}\cdot G_{ij,\x(0,s)}^{-1}=\xiC\quad\mbox{ in $\Omega_{ij,r}$.}$$
The point is that $\xi_{0,s}$ and $\whxi_{ij,0,s}$ are independent of $s$.
Proposition \ref{prop-differentiability} follows from Proposition \ref{prop-tlogt} in Appendix \ref{appendix-tlogt}.
\cqfd
\subsection{The immersion $f_t$}
\label{section-geometry-immersion}
We introduce the following notations for $0<r<\frac{\varepsilon}{2}$:
\begin{itemize}
\item $\SSigma_t=\SSigma_{t,\x(t)}$, where the compact Riemann surface $\SSigma_{t,\x}$ was defined in Section \ref{section-setup-openingnodes}.
\item $\Omega_{t,r}$ is the domain defined as $\SSigma_t$ minus the points $\infty_i$
for $i\in[1,N]$, $\infty_{ij}$ for $(i,j)\in I$ and the closed disks
$\overline{D}(\pi_k,r)$ for $i\in [1,N]$ and $k\in R_i$.
Note that the domains $\Omega_{i,r}$ and $\Omega_{ij,r}$ defined in Section
\ref{section-geometry-differentiability} are included in $\Omega_{t,r}$.
\item $\wtOmega_{t,r}$ is the universal cover of $\Omega_{t,r}$.
\item $\Omega'_{t,r}=\Omega_{t,r}\setminus\{m_{ij}(t):(i,j)\in I\}$.
The potential $\xi_t$ is holomorphic in $\Omega'_{t,r}$.
\item $\wtOmega'_{t,r}=p^{-1}(\Omega'_{t,r})\subset\wtOmega_{t,r}$.
(This is not the universal cover of $\Omega'_{t,r}$.)
\item $\wt0_1\in\wtOmega_{t,r}$ is an arbitrary point in the fiber of $0_1$.
\item $\Phi_t$ is the solution of the Cauchy Problem
$d\Phi_t=\Phi_t\xi_t$ in $\wtOmega'_{t,r}$
with initial condition $\Phi_t(\wt0_1)=I_2$. Note that $\Phi_t$ is well defined in $\wtOmega'_{t,r}$ because the Regularity Problem at $m_{ij}(t)$ is solved.
\item $f_t=\Sym(\Uni(\Phi_t))$ is the (branched) immersion obtained by the DPW method from $\Phi_t$.
Since the Monodromy Problem for $\Phi_t$ is solved, $f_t$ descends to a well defined (branched) immersion in $\Omega_{t,r}$, still denoted $f_t$.
\item For $i\in[2,N]$, $\wt0_i\in\wtOmega_{t,r}$ is a point in the fiber of $0_i$, depending
continuously on $t$ and defined as follows:
choose a path $c$ from $\v_1$ to $\v_i$ on the graph $\Gamma$ of the
form $c=\prod_{j=1}^k e_{i_j \,i_{j+1}}$ with $i_1=1$ and $i_{k+1}=i$.
Let $\gamma=\prod_{j=1}^k \Gamma_{i_j \,i_{j+1}}$: this is a path from 
$0_1$ to $0_i$ in $\Omega_{t,r}$.
Let $\wtgamma$ be the lift of $\gamma$ to $\wtOmega_{t,r}$ such that $\wtgamma(0)=\wt0_1$.
We take $\wt0_i=\wtgamma(1)$.
\end{itemize}
\begin{proposition}
\label{prop-immersion-extends}
For $t>0$ small enough, $f_t$ extends analytically to
$$\Sigma_t:=\SSigma_t\setminus\{p_1^0(t),\cdots,p_n^0(t)\}.$$
\end{proposition}
Proof:
\begin{itemize}
\item Let $i\in[1,N]$.
Since the Monodromy Problem is solved, Propositions \ref{prop-balancing},
\ref{prop-differentiability} and
Corollary \ref{cor-regularity-infty} in Appendix \ref{appendix-regularity} imply that $\infty_i$ is a removable singularity:
there exists a gauge $G_{\infty i,t}$ such that $\xi_t\cdot G_{\infty i,t}$ extends holomorphically at $\infty_i$. So $f_t$ extends analytically at $\infty_i$.
The gauge $G_{\infty i,t}$ given by Corollary \ref{cor-regularity-infty} has the following form:
\begin{equation}
\label{eq-G8it}
G_{\infty i,t}(z,\lambda)=\matrix{z&0\\\lambda g_{\infty i}(t)(\lambda)&z^{-1}}
\end{equation}
Moreover, $g_{\infty i}$ is a $C^1$ function of $t$ and $g_{\infty i}(0)=-1$.
\item Let $(i,j)\in I$. In the same way, Propositions \ref{prop-regularity2}, \ref{prop-differentiability} and Theorem \ref{thm-regularity-infty} in Appendix \ref{appendix-regularity}
imply that there exists a gauge $G_{\infty ij,t}$ such that $\whxi_{ij,t}\cdot G_{\infty ij,t}$ extends holomorphically at $\infty_{ij}$. So $f_t$ extends analytically at $\infty_{ij}$.
The gauge $G_{\infty ij,t}$ given by Theorem \ref{thm-regularity-infty}
has the following form:
\begin{equation}
\label{eq-G8ijt}
G_{\infty ij,t}(z,\lambda)=\matrix{z^{-1}&g_{\infty ij}(t)(\lambda)\\0&z}.
\end{equation}
Moreover, $g_{\infty ij}$ is a $C^1$ function of $t$ and $g_{\infty ij}(0)=-1$.
\item It remains to prove that $f_t$ extends analytically to the punctured disks
$D^*(p_k^0(t),\varepsilon)$.
Fix $i\in[1,N]$ and $k\in R_i$.
Arguing as in the proof of Proposition 4 of \cite{nnoids}, we consider the change of variable
$$z=\psi_{t,\lambda}(w)=p_k(t)(\lambda)+w$$
and the following domains (independent of $\lambda$):
$$U=\{w\in\C:\smallfrac{3}{4}\varepsilon<|w|<\varepsilon\}$$
$$V_k=\{z\in\C_i:\smallfrac{\varepsilon}{2}<|z-\pi_k|<2\varepsilon\}\subset\Omega_{t,r}.$$
Since $p_k(0)=\pi_k$, we have $\psi_{t,\lambda}(U)\subset V_k$ for $t$ small enough.
Let $\widetilde{U}$ be the universal cover of $U$.
Let $\wtalpha_k$ be the lift of $\alpha_k$ to $\wtOmega_{t,r}$ such that
$\wtalpha_k(0)=\wt0_i$.
Let $\widetilde{V}_k\subset\wtOmega_{t,r}$
be the component of $p^{-1}(V_k)$ containing $\wtalpha_k(1)$.
By Proposition \ref{prop-restriction}, $\widetilde{V}_k$ is a universal cover of $V_k$.
Lift $\psi_{t,\lambda}$ to $\widetilde{\psi}_{t,\lambda}:\widetilde{U}\to\widetilde{V}_k$.
Define for $w\in\widetilde{U}$:
$$\whPhi_t(w,\lambda)=\Phi_t(\widetilde{\psi}_{t,\lambda}(w),\lambda).$$
Let $\whf_t=\Sym(\Uni(\whPhi_t))$. By Corollary 1 in \cite{nnoids},
$\whf_t$ descends to a well defined immersion on $U$ and
\begin{equation}
\label{eq-immersion-1}
\forall w\in U, \quad \whf_t(w)=f_t(\psi_{t,0}(w))=f_t(p_k^0(t)+w).
\end{equation}
Now $\whPhi_t$ solves $d\whPhi_t=\whPhi_t\whxi_t$ with
$\whxi_t=\psi_{t,\lambda}^*\xi_t$. Since the only pole of $\whxi_t$ in $D(0,\varepsilon)$ is at $w=0$, $\whf_t$ extends analytically to $D^*(0,\varepsilon)$. Hence by Equation \eqref{eq-immersion-1}, $f_t$ extends analytically to $D^*(p_k^0(t),\varepsilon)$.\cqfd
\end{itemize}
\begin{proposition} For $t>0$ small enough, $f_t$ is unbranched in $\Sigma_t$ (meaning that it is a regular immersion).
\end{proposition}
Proof: let $\beta_t^0=\xi_{t;12}^{(-1)}$.
\begin{itemize}
\item
By Proposition \ref{prop-gauging}, we have in $\C_{ij}$:
$$\beta_t^0=(G_{ij,\x(t);22}^0)^2\whxi_{ij,t;12}^{\,(-1)}
=(2g_{ij}^0(t)(z-m_{ij}^0))^2\whxi_{ij,t;12}^{\,(-1)}.$$
By Proposition \ref{prop-regularity}, $\whxi_{ij,t}$ is holomorphic at $m_{ij}$,
so $\beta_t^0$ has order at least $2$ at $m_{ij}^0$,
with equality if and only if the potential $\whxi_{ij,t}$ is regular
at $m_{ij}$.
\item Let
$$\whxi_{\infty ij,t}=\whxi_{ij,t}\cdot G_{\infty ij,t}=\xi_t\cdot G_{ij,\x(t)}^{-1}\cdot G_{\infty ij,t}.$$
By Proposition \ref{prop-gauging} and Equation \eqref{eq-G8ijt}, we have in $\C_{ij}$:
$$\beta_t^0=(G_{\infty ij,t;11}^0)^2 (G_{ij,\x(t);22}^0)^2\whxi_{\infty ij,t;12}^{\,(-1)}
=(2g_{ij}^0(t)(z-m_{ij}^0))^2\,z^{-2}\,\whxi_{\infty ij,t;12}^{\,(-1)}.$$
Since $\whxi_{\infty ij,t}$ extends holomorphically at $\infty_{ij}$, $\beta_t^0$ extends holomorphically
at $\infty_{ij}$ and is non-zero if and only if the potential $\whxi_{\infty ij,t}$ is regular
at $\infty_{ij}$.
\item The only remaining poles of $\beta_t^0$ are double poles at $\infty_i$ for $i\in[1,N]$.
Since the genus of $\SSigma_t$ is $g=|I|-N+1$, the number of zeros
of $\beta_t^0$ is
$$2g-2+2N=2(|I|-N+1)-2+2N=2|I|.$$
Hence for $(i,j)\in I$, $\beta_t^0$ has a zero of multiplicity exactly 2 at $m_{ij}^0$,
and has no other zeros in $\SSigma_t$. By the previous points, $f_t$ is unbranched in
$\Sigma_t$.
\end{itemize}
\cqfd
\begin{remark}
This is the only point in the paper where we need to know the genus of
$\SSigma_t$.
\end{remark}
\subsection{Spherical parts}
\label{section-geometry-spherical}
Recall that $\v_1=(-1,0,0)$ and $\fS$ parametrizes the unit sphere centered at
this point.
\begin{proposition}
\label{prop-spherical}
For $i\in[1,N]$, we have for $t$ small enough:
$$|| f_t-f^S-\v_i(t)+\v_1||_{C^1(\Omega_{i,r}\cap D(0,\frac{1}{r}))}\leq C(r)t.$$
\end{proposition}
\medskip

Proof: let $\wtOmega_{i,r}\subset\wtOmega_{t,r}$ be the component of $p^{-1}(\Omega_{i,r})$ 
containing $\wt0_i$.
By Proposition \ref{prop-restriction}, $\wtOmega_{i,r}$ is a universal cover of
$\Omega_{i,r}$.
Let $\Phi_{i,t}$ be the solution of $d\Phi_{i,t}=\Phi_{i,t}\xi_t$ in $\wtOmega_{i,r}$
with initial condition $\Phi_{i,t}(\wt0_i,\lambda)=I_2$ and $f_{i,t}$ be the corresponding
immersion in $\Omega_{i,r}$.
Then
$$\Phi_t(z,\lambda)=\Phi_t(\wt0_i,\lambda)\Phi_{i,t}(z,\lambda)\quad\mbox{ in $\wtOmega_{i,r}$.}$$
By Equation \eqref{eq-principal-morphism} and the definition of $\wt0_i$ in
Section \ref{section-geometry-immersion}:
$$\Phi_t(\wt0_i,\lambda)=\prod_{j=1}^k\boP(\xi_t,\Gamma_{i_j\,i_{j+1}})(\lambda).$$
Since the Monodromy Problem \eqref{pb-monodromy-edges} is solved, we have:
\begin{equation}
\label{eq-spherical-0}
\left\{\begin{array}{l}
\Phi_t(\wt0_i,\cdot)\in\Lambda SU(2)\\
\Phi_t(\wt0_i,1)=I_2\\
\displaystyle\frac{\partial\Phi_t}{\partial\lambda}(\wt0_i,1)=-\ii\sum_{j=1}^k(V_{i_{j+1}}(t)-V_{i_j}(t))
=-\ii (V_i(t)-V_1).\end{array}\right.\end{equation}
By Equation \eqref{eq-dressing}, we obtain:
\begin{equation}
\label{eq-spherical-1}
f_t(z)=f_{i,t}(z)+\v_i(t)-\v_1\quad\mbox{ in $\Omega_{i,r}$.}
\end{equation}
By Point 1 of Proposition \ref{prop-differentiability}:
$$||\xi_t(z,\cdot)-\xiS(z,\cdot)||_{\boW}\leq C(r)t \quad\mbox{ in $\Omega_{i,r}\cap D(0,\frac{1}{r})$}.$$
Let $K$ be a bounded subset of $\wtOmega_{i,r}$ such that $p(K)=\Omega_{i,r}
\cap D(0,\frac{1}{r})$.
Since $\Phi_{i,t}(\wt0_i,\lambda)=\PhiS(0,\lambda)=I_2$, we have for $t$ small enough,
by standard perturbation theory of Ordinary Differential Equations:
\begin{equation}
\label{eq-spherical-3}
||\Phi_{i,t}(z,\cdot)-\PhiS(z,\cdot)||_{\boW}\leq C(r)t\quad\mbox{ in $K$}.
\end{equation}
Let $(F_{i,t},B_{i,t})$ be the Iwasawa decomposition of $\Phi_{i,t}$. Then
$$||F_{i,t}(z,1)-\FS(z,1)||\leq C(r)t\quad\mbox{ in $K$}$$
$$||B_{i,t}(z,0)-\BS(z,0)||\leq C(r)t\quad\mbox{ in $K$}.$$
Using Equations \eqref{eq-dfx}, \eqref{eq-dfy} and \eqref{eq-metric}, we obtain:
$$||df_{i,t}-d\fS||_{C^0(\Omega_{i,r}\cap D(0,\frac{1}{r}))}\leq C(r)t.$$
Hence since $f_{i,t}(0_i)=\fS(0)=0$:
\begin{equation}
\label{eq-spherical-2}
||f_{i,t}-\fS||_{C^1(\Omega_{i,r}\cap D(0,\frac{1}{r}))}\leq C(r)t.
\end{equation}
Proposition \ref{prop-spherical} follows from Equations \eqref{eq-spherical-1}
and \eqref{eq-spherical-2}.
\cqfd
\medskip

To study the immersion $f_t$ in a neighborhood of $\infty_i$ we consider the
change of variable
$$\varphi_i:D(0,r)\to \{z\in\C_i:|z|>\smallfrac{1}{r}\}\cup\{\infty_i\},\quad
\varphi_i(w)=\frac{1}{w}.$$
\begin{proposition}
\label{prop-inftyi}
For $i\in[1,N]$, we have for $t$ small enough
$$|| f_t\circ\varphi_i-f^S\circ\varphi_i-\v_i(t)+\v_1||_{C^1(D(0,r))}\leq C(r)t.$$
\end{proposition}
Proof: let $D_i\subset\CC_i$ be a simply connected domain containing $0_i$, $\infty_i$ and
the domain $|z|>\frac{1}{r}$, and disjoint from the disks $D(\pi_{ij},\varepsilon)$
for $j\in E_i$ and $D(\pi_k,\varepsilon)$ for $k\in R_i$.
The potential $\xi_t$ is holomorphic in $D_i\setminus\{\infty_i\}$.
Since the Regularity Problem at $\infty_i$ is solved,
the solution of $d\Phi_{i,t}=\Phi_{i,t}\xi_t$ in $D_i$ with initial condition 
$\Phi_{i,t}(0_i,\lambda)=I_2$ is well defined in $D_i\setminus\{\infty_i\}$.
By Equation \eqref{eq-spherical-3} (with $r$ replaced by $\frac{r}{4}$):
\begin{equation}
\label{eq-spherical-4}
||\Phi_{i,t}(\smallfrac{2}{r},\cdot)-\PhiS(\smallfrac{2}{r},\cdot)||_{\boW}\leq C(r)t.
\end{equation}
Recalling Equation \eqref{eq-G8it}, we define
$$G^S(z,\lambda)=G_{\infty i,0}(z,\lambda)=\matrix{z&0\\-\lambda&z^{-1}}.$$
Since $g_{\infty i}$ is a $C^1$ function of $t$:
\begin{equation}
\label{eq-spherical-5}
||G_{\infty i,t}(\smallfrac{2}{r},\cdot)-G^S(\smallfrac{2}{r},\cdot)||_{\boW}\leq C(r)t.
\end{equation}
Define
$$\whxi_{i,t}=\varphi_i^*(\xi_t\cdot G_{\infty i,t})$$
$$\whxi^S=\varphi_i^*(\xiS\cdot G^S).$$
$$\whPhi_{i,t}=\varphi_i^*(\Phi_{i,t} G_{\infty i,t})$$
$$\whPhi^S=\varphi_i^*(\PhiS G^S)$$
Since the Regularity Problem at $\infty_i$ is solved, $\whxi_{i,t}$ and
$\whPhi_{i,t}$ extend holomorphically to $D(0,r)$.
By Equations \eqref{eq-spherical-4} and \eqref{eq-spherical-5}:
\begin{equation}
\label{eq-spherical-6}
||\whPhi_{i,t}(\smallfrac{r}{2},\cdot)-\whPhi^S(\smallfrac{r}{2},\cdot)||_{\boW}\leq C(r)t.
\end{equation}
Since $\whxi_{i,t}$ depends $C^1$ on $t$:
\begin{equation}
\label{eq-spherical-7}
||\whxi_{i,t}(z,\cdot)-\whxi^S(z,\cdot)||_{\boW}\leq C(r)t
\quad\mbox{ in $D(0,r)$}.
\end{equation}
Let
$$\whf_{i,t}=\Sym(\Uni(\whPhi_{i,t}))=f_{i,t}\circ\varphi_i$$
$$\whf^S=\Sym(\Uni(\whPhi^S))=\fS\circ\varphi_i.$$
Arguing as in the proof of Proposition \ref{prop-spherical}, Equations \eqref{eq-spherical-6} and \eqref{eq-spherical-7}
imply that
$$||\whf_{i,t}-\whf^S||_{C^1(D(0,r)}\leq C(r)t.$$
Proposition \ref{prop-inftyi} follows from Equation \eqref{eq-spherical-1}.\cqfd
\subsection{Delaunay ends}
\label{section-geometry-delaunay}
\begin{proposition}
\label{prop-delaunay}
Let $i\in[1,N]$ and $k\in R_i$.
\begin{myenumerate}
\item $a_k(t)$ is a real constant (with respect to $\lambda$).
\item $f_t$ has a Delaunay end of weight $8\pi ta_k(t)$ at
$p_k^0(t)$. More precisely:
\item There exists uniform positive numbers $\epsilon$, $c$, $T$ 
and a family of Delaunay immersions
$f_{k,t}^{\boD}:\C^*\to\R^3$ such that for $0<t<T$ and
$0<|w|<\epsilon$
$$||f_t(p_k^0(t)+w)-f_{k,t}^{\boD}(w)||\leq c\,t|w|^{1/2}.$$
\item If $\tau_k>0$, the restriction of $f_t$ to $D^*(p_k^0(t),\epsilon)$
is an embedding.
\item The axis of $f_{k,t}^{\boD}$ converges as $t\to 0$ to
the half-line through $\v_i$ spanned by the vector $\u_k$.
\end{myenumerate}
\end{proposition}
Proof: Points 1 and 2 are proved in Proposition 4 of \cite{nnoids} by gauging the 
potential to a potential with a simple pole and a standard Delaunay residue.
The immersion $f_t$ then has a Delaunay end by Theorem 3.5 in \cite{kilian-rossman-schmitt}.
As in Section \ref{section-monodromy-ends}, the only properties of the potential that are used to prove these
results are Properties \eqref{eq-properties-nnoids}.
\medskip

To prove Points 3 and 4, we use Corollary 2 in \cite{raujouan} as in
Section 10 of \cite{nnoids}. There is a technical issue however, which is that
this result requires the potential to be of class $C^2$ with respect to
$t$ and we do not have that regularity.
As in Section \ref{section-geometry-differentiability}, we denote
$\xi_{t,s}=\xi_{t,\x(t,s)}$, which depends smoothly on $t$ and $s$.
Let $\Phi_{i,t,s}$ be the solution
of $d\Phi_{i,t,s}=\Phi_{i,t,s}\xi_{t,s}$ in $\wtOmega_{i,r}$
with initial condition $\Phi_{i,t,s}(\wt0_i,\lambda)=I_2$.
The crucial point is that $\boM(\Phi_{i,t,s},\gamma_k)=\boP(\xi_{t,s},\gamma_k)$ solves the Monodromy Problem \eqref{pb-monodromy-ends}
for all $(t,s)$ in a neighborhood of $(0,0)$, because that problem was solved
before Section \ref{section-monodromy-edges}.
Applying Corollary 2 in \cite{raujouan} as in Section 10 of \cite{nnoids}, for fixed value
of $s$,
there exist uniform positive numbers $\epsilon$, $c$, $T$ 
and a family of Delaunay immersions
$f_{k,t,s}^{\boD}:\C^*\to\R^3$ such that for $0<t<T$
and $0<|w|<\epsilon$
$$||f_{i,t,s}(p_k^0(t,s)+w)-f_{k,t,s}^{\boD}(w)||\leq c\,t|w|^{1/2}.$$
Moreover if $\tau_k>0$, the restriction of $f_{i,t,s}$ to $D^*(p_k^0(t,s),\epsilon)$
is an embedding. The numbers $\epsilon$, $c$ and $T$ can be chosen independent of $s$
by continuity.
Specializing to $s=t\log t$ we obtain
\begin{equation}
\label{eq-delaunay-1}
||f_{i,t}(p_k^0(t)+w)-f_{k,t,t\log t}^{\boD}(w)||\leq c\,t|w|^{1/2}.
\end{equation}
We define
$$f_{k,t}^{\boD}=f_{k,t,t\log t}^{\boD}+\v_i(t)-\v_1.$$
Equations \eqref{eq-spherical-1} and \eqref{eq-delaunay-1} give Point 3 of Proposition \ref{prop-delaunay}.
By Proposition 5 of \cite{nnoids}, the axis of $f_{k,t,s}^{\boD}$ converges to the line
through $(0,0,-1)$ directed by $\u_k$, which gives Point 5.
\cqfd
\subsection{Catenoidal parts}
\label{section-geometry-catenoidal}
\begin{proposition}
\label{prop-catenoidal}
For $(i,j)\in I$,
there exists
a continuous family of rigid motions $h_{ij,t}$, $t\in[0,\epsilon)$
and a complete minimal immersion
$$\Psi_{ij}:\C\cup\{\infty\}\setminus\{q_{ij},q_{ji}\}\to\R^3$$
such that:
\begin{myenumerate}
\item $h_{ij,0}$ is the translation of vector $\v_i+\u_{ij}$.
\item The restriction of $f_t$ to $\Omega_{ij,r}\cap D(0,\frac{1}{r})$ satisfies
$$\lim_{t\to 0} ||\smallfrac{1}{t}h_{ij,t}^{-1}\circ f_t-\Psi_{ij}||_{C^1(\Omega_{ij,r}\cap D(0,\frac{1}{r}))}=0.$$
\item $\Psi_{ij}$ parametrizes a catenoid with necksize $4|\tau_{ij}|$
and axis directed by the vector $\u_{ij}$.
(The axis is oriented from the end at $q_{ij}$ to
the end at $q_{ji}$).
\item The Gauss map of $\Psi_{ij}$ points towards the ``inside'' if $\tau_{ij}>0$
and the ``outside'' if $\tau_{ij}<0$.
(By ``inside'', we mean the component of the complement of the catenoid which contains
its axis.)
\end{myenumerate}
\end{proposition}
Proof: 
recall that $\delta_{ij}$ is a path from $0_i$ to $O_{ij}$ (see Section \ref{section-monodromy-paths}).
Let $\gamma$ be the composition of $\delta_{ij}$ with a fixed path from $O_{ij}$ to
$0_{ij}$ in $\Omega_{ij,r}$. Let $\wtgamma$ be the lift of $\gamma$ to $\wtOmega_{t,r}$ such that $\wtgamma(0)=\wt0_i$.
We define $\wt0_{ij}=\wtgamma(1)$.
Let $\wtOmega_{ij,r}\subset\wtOmega_{t,r}$ be the component of $p^{-1}(\Omega_{ij,r})$
containing $\wt0_{ij}$. By Proposition \ref{prop-restriction}, $\wtOmega_{ij,r}$ is a universal cover of $\Omega_{ij,r}$.
Let
$$\whPhi_{ij,t}=\Phi_t G_{ij,\x(t)}^{-1}.$$
\begin{claim}
\label{claim-catenoidal}
We have
$$\whPhi_{ij,0}(\wt0_{ij},\lambda)=\Phi_0(\wt0_i,\lambda)U(q_{ij},\lambda)\in\Lambda SU(2)$$
where $U(q_{ij},\lambda)$ is given by Equation \eqref{eq-Uqij}.
\end{claim}
Proof: recall that $G_{ij,\x(0)}=G_{ij}$.
We have, by the computations in the proof of Proposition
\ref{prop-Pij} (omitting the variable $\lambda$):
\begin{eqnarray*}
\whPhi_{ij,0}(\wt0_{ij})&=&
\Phi_0(\wt0_{ij})G_{ij}(0)^{-1}\\
&=&\Phi_0(\wt0_i)\boP(\xi_0,\delta_{ij})\boP(\xi_0,O_{ij},0_{ij})\,G_{ij}(0)^{-1}\\
&=&\Phi_0(\wt0_i)\PhiS(\pi_{ij})[\PhiC(q_{ij})G_{ij}(q_{ij})]^{-1}\,[\PhiC(0)G_{ij}(0)] \,G_{ij}(0)^{-1}\\
&=&\Phi_0(\wt0_i)U_{ij}(q_{ij}).\end{eqnarray*}
.\cqfd

Let
$$H_{ij,t}=\Uni(\whPhi_{ij,t}(\wt0_{ij},\cdot))\matrix{0&-1\\1&0}\in\Lambda SU(2)$$
and 
$h_{ij,t}$ be the rigid motion given by the action \eqref{eq-actionSU2} of
$H_{ij,t}$ on $\su(2)$.
At $t=0$, we have by claim \ref{claim-catenoidal}:
$$H_{ij,0}(\lambda)=\whPhi_{ij,0}(\wt0_{ij},\lambda)\matrix{0&-1\\1&0}
=\Phi_0(\wt0_i,\lambda)U(q_{ij},\lambda)\matrix{0&-1\\1&0}.$$
By Equations \eqref{eq-Uqij} and \eqref{eq-spherical-0}:
$$H_{ij,0}(1)=I_2\matrix{0&1\\-1&0}\matrix{0&-1\\1&0}=I_2$$
\begin{eqnarray*}
\ii\frac{\partial H_{ij,0}}{\partial\lambda}(1)&=&
\ii\frac{\partial\Phi_0}{\partial\lambda}(\wt0_i,1)+
\frac{\ii}{1+|\pi_{ij}|^2}\matrix{\pi_{ij}&-|\pi_{ij}|^2\\-|\pi_{ij}|^2&-\overline{\pi}_{ij}}
\matrix{0&-1\\1&0}\\
&=&V_i-V_1-\frac{\ii}{2}
\frac{1}{(1+|\pi_{ij}|^2)}\matrix{2|\pi_{ij}|^2&2\pi_{ij}\\2\overline{\pi}_{ij}&-2|\pi_{ij}|^2}
\end{eqnarray*}
Hence $h_{ij,0}$ is the translation of vector:
$$\v_i-\v_1+\frac{1}{1+|\pi_{ij}|^2}\left(2\Re(\pi_{ij}),2\Im(\pi_{ij}),-2|\pi_{ij}|^2\right)
=\v_i+\pi_S^{-1}(\pi_{ij})=\v_i+\u_{ij}$$
which proves Point 1 of Proposition \ref{prop-catenoidal}.
We consider the dressing by $H_{ij,t}^{-1}$ and define in $\wtOmega_{ij,r}$:
$$\wtPhi_{ij,t}=H_{ij,t}^{-1}\whPhi_{ij,t}$$
$$\wtf_{ij,t}=\Sym(\Uni(\wtPhi_{ij,t}))=h_{ij,t}^{-1}\circ f_t.$$
At $t=0$, we have in $\Omega_{ij,r}$, since $\whxi_{ij,0}=\xiC$:
\begin{equation}
\label{eq-catenoidal-1}\wtPhi_{ij,0}(z,\lambda)=
\matrix{0&1\\-1&0}\whPhi_{ij,0}(\wt0_{ij},\lambda)^{-1}\whPhi_{ij,0}(z,\lambda)
=\matrix{0&1\\-1&0}
\matrix{1&0\\z&1}=\matrix{z&1\\-1&0}.
\end{equation}
Consequently, the unitary part of $\wtPhi_{ij,0}$ is in
$\Lambda SU(2)\cap \Lambda_+ SL(2,\C)$ so is constant with respect to $\lambda$.
By the Sym-Bobenko formula, $\wtf_{ij,0}\equiv 0$ in $\Omega_{ij,r}$.
To compute the limit of $\frac{1}{t}\wtf_{ij,t}$ as $t\to 0$, we use
Theorem \ref{thm-minimal} in Appendix \ref{appendix-minimal}.
By Proposition \ref{prop-differentiability}, $\whxi_{ij,t}$ is of class $C^1$ in $\Omega_{ij,r}$,
and
\begin{eqnarray*}
\frac{d}{dt}\whxi_{ij,t}|_{t=0}&=&\frac{d}{dt}\whxi_{ij,t,\x(t,t\log t)}|_{t=0}\\
&=&\frac{\partial}{\partial t}\whxi_{ij,t,\x(t,s)}|_{(t,s)=(0,0)}\quad\mbox{by the proof of Proposition \ref{prop-tlogt}}\\
&=&\frac{d}{dt}\whxi_{ij,t,\x(0)}|_{t=0}\quad\mbox{ by the chain rule and using $\whxi_{ij,0,\x}=\xiC$}.
\end{eqnarray*}
By Proposition \ref{prop-bmij0}, we have $b_{mij}^0(0,\x(0))=c_{mij}^0(0,\x(0))=0$, so
$$\Theta_{0,\x(0);12}^{(-1)}=0.$$
By definition of $\xi_{t,\x}$ and Theorem \ref{thm-openingnodes-derivee}, we obtain
in $\C_{ij}$:
$$\frac{d}{dt}\xi_{t,\x(0);12}^{(-1)}|_{t=0}=-\frac{r_{ij}\,dz}{(z-q_{ij})^2}
-\frac{r_{ji}\,dz}{(z-q_{ji})^2}+\frac{b_{ij}^0dz}{z-q_{ij}}
+\frac{b_{ji}^0dz}{z-q_{ji}}=\frac{4\tau_{ij}(z-\mu_{ij})^2\,dz}{\rho_{ij}^2(z-q_{ij})^2(z-q_{ji})^2}.$$
By Proposition \ref{prop-gauging} with $G=G_{ij,\x(0)}^{(-1)}$:
$$\whxi_{ij,t,\x(0);12}^{(-1)}=(G_{ij,\x(0);11}^{(0)})^2\xi_{t,\x(0);12}^{(-1)}
=\frac{1}{4(z-\mu_{ij})^2}\xi_{t,\x(0);12}^{(-1)}.$$
Hence by differentiation:
\begin{equation}
\label{eq-catenoidal-2}
\frac{d}{dt}\whxi_{ij,t;12}^{(-1)}|_{t=0}
=\frac{d}{dt}\whxi_{ij,t,\x(0);12}^{(-1)}|_{t=0}
=\frac{1}{4(z-\mu_{ij})^2}\frac{d}{dt}\xi_{t,\x(0);12}^{(-1)}|_{t=0}
=\frac{\tau_{ij}\,dz}{\rho_{ij}^2(z-q_{ij})^2(z-q_{ji})^2}.
\end{equation}
By Theorem \ref{thm-minimal} and Equations \eqref{eq-catenoidal-1} and
\eqref{eq-catenoidal-2},
we have in $\Omega_{ij,r}\cap D(0,\frac{1}{r})$:
$$\lim_{t\to 0}\frac{1}{t}\wtf_{ij,t}=\sigma\circ\psi_{ij}=:\Psi_{ij}$$
where $\sigma(x_1,x_2,x_3)=(x_1,-x_2,x_3)$ and
$\psi_{ij}$ is a minimal
immersion with the following Weierstrass data:
$$g(z)=\frac{-1}{z},\quad
\omega(z)=\frac{4\tau_{ij}\,z^2\,dz}{\rho_{ij}^2(z-q_{ij})^2(z-q_{ji})^2}.$$
The immersion $\psi_{ij}$ is regular at $z=0$ and $z=\infty$.
A computation gives
\begin{equation}
\label{eq-catenoidal-3}
\Res_{q_{ij}}\left(\frac{1}{2}(1-g^2)\omega,\frac{\ii}{2}(1+g^2)\omega,g\omega\right)
=\frac{4\tau_{ij}}{1+|\pi_{ij}|^2}\left(
2\Re(\pi_{ij}),-2\Im(\pi_{ij}),1-|\pi_{ij}|^2\right).
\end{equation}
Since this residue is real, $\psi_{ij}$ is well defined in $\CC\setminus\{q_{ij},q_{ji}\}$
and has two catenoidal ends at $z=q_{ij}$ and
$z=q_{ji}$ so is a catenoid. The necksize and the direction of the axis are determined
by Equation \eqref{eq-catenoidal-3} in a standard way.
\cqfd
\begin{remark}
Point 2 of Proposition \ref{prop-catenoidal} can be extended to a neighbhorhood of
$\infty_{ij}$ in the way that Proposition \ref{prop-spherical} was extended to
a neighborhood of $\infty_i$ in Proposition \ref{prop-inftyi}.
\end{remark}
\subsection{Transition regions}
For $(i,j)\in I\cup I^*$ and $0<r\leq\frac{\varepsilon}{2}$, let $\boA_{ij,t,r}$ be the
annulus
$$\boA_{ij,t,r}=\{z\in\C_i:\smallfrac{|t_{ij}(t)|}{r}\leq |v_{ij}|<r\}$$
 which is identified
with the annulus $\frac{|t_{ij}(t)|}{r}\leq |w_{ij}|<r$ in $\C_{ij}$.
Let $N_t$ be the Gauss map of $f_t$ in $\Sigma_t$.
The goal of this section is to prove:
\begin{proposition}
\label{prop-transition}
For all $\alpha>0$, there exists $r>0$ and $T>0$ such that
for all $(i,j)\in I\cup I^*$ and $t<T$:
\begin{myenumerate}
\item $|N_t+\u_{ij}|\leq \alpha$ in $\boA_{ij,t,r}$.
\item $f_t(\boA_{ij,t,r})$ is a graph over a domain in the plane orthogonal to $\u_{ij}$.
\item If moreover $\tau_{ij}>0$, then $f_t(\boA_{ij,t,r})$ and $f_t(\boA_{ji,t,r})$ are disjoint.
\end{myenumerate}
\end{proposition}
Proof:
\begin{myenumerate}
\item let $\Phi_{i,t}$ be the solution of $d\Phi_{i,t}=\Phi_{i,t}\xi_t$ with initial
condition $\Phi_{i,t}(\wt0_i)=I_2$.
The idea is to prove that $\Phi_{i,t}$ is close to $\PhiS(\pi_{ij},\cdot)$ on the two boundary
components of $\boA_{ij,t,r}$ and extend this to the interior of $\boA_{ij,t,r}$ by the
maximum principle. The problem is that $\Phi_{i,t}$ is not well-defined in $\boA_{ij,t,r}$,
so we multiply it by a suitable factor $\boG_{ij,t}$ to obtain a well-defined holomorphic
function in $\boA_{ij,t,r}$.
\medskip

Let $\alpha_{ij}$ and $\beta_{ij}$ be the paths defined in Section \ref{section-monodromy-paths} with $\varepsilon$ replaced by $r$, so $\alpha_{ij}$
goes from $0_i$ to $p_{ij}(t)+r$ in $\C_i$, and $\beta_{ij}$ goes from
$p_{ij}(t)+r$ in $\C_i$ to $q_{ij}(t)+r$ in $\C_{ij}$.
Let $\wtalpha_{ij}$ be the lift of $\alpha_{ij}$ to $\wtOmega_{t,r}$ such that
$\wtalpha_{ij}(0)=\wt0_i$ and $\wtboA_{ij,t,r}\subset\wtOmega_{t,r}$
be the component of $p^{-1}(\boA_{ij,t,r})$ containing $\wtalpha_{ij}(1)$.
By Proposition \ref{prop-restriction}, $\wtboA_{ij,t,r}$ is a universal cover
of $\boA_{ij,t,r}$.
By uniqueness of the universal cover up to isomorphism, we may identify
$\wtboA_{ij,t,r}$ with the domain
$$\{\wtz\in\C:\log|t_{ij}|-\log r\leq\Re(\wtz)\leq\log r\}$$
so that $\wtalpha_{ij}(1)$ is identified with $\wtz=\log r$ and
the universal covering map is $\wtz\mapsto p_{ij}(t)+e^{\wtz}$.
Then the lift $\wtbeta_{ij}$ of $\beta_{ij}$ such that $\wtbeta_{ij}(0)=\wtalpha_{ij}(1)$ is the path parametrized by $s\mapsto\wtbeta_{ij}(s)=(1-2s)\log r+s\log t_{ij}$.
Under this identification, the translation $\sigma:\wtz\mapsto\wtz+2\pi\ii$ is a generator of
$\mbox{Deck}(\wtboA_{ij,t,r}/\boA_{ij,t,r})$.
By Equations  \eqref{eq-principal-morphism} and \eqref{eq-boM-boP}:
$$\boM(\Phi_{i,t},\sigma)=\Phi_{i,t}(\wtalpha_{ij}(1),\cdot)
\boP(\xi_t,C_{ij})\Phi_{i,t}(\wtalpha_{ij}(1),\cdot)^{-1}
=\boP(\xi_t,\alpha_{ij})\boP(\xi_t,C_{ij})\boP(\xi_t,\alpha_{ij})^{-1}
=\boP(\xi_t,\gamma_{ij}).$$
We define for $\wtz\in\wtboA_{ij,t,r}$:
$$\boG_{ij,t}(\wtz,\lambda)=\exp\left[\frac{\wtz}{2\pi\ii}\log\boP(\xi_t,\gamma_{ij})(\lambda)\right]
=\exp\left[\frac{\wtz}{2\pi\ii}t(\lambda-1)\wtM_{ij}(t,\x(t))(\lambda)\right].$$
Then
$$\boG_{ij,t}(\sigma(\wtz),\lambda)=\boP(\xi_t,\gamma_{ij})(\lambda)\boG_{ij,t}(\wtz,\lambda).$$
Consequently, the function $\boG_{ij,t}^{-1}\Phi_{i,t}$ descends to a well-defined holomorphic function in $\boA_{ij,t,r}$ which we denote $\boH_{ij,t}$:
$$\boH_{ij,t}(p_{ij}(t)+e^{\wtz},\lambda)=\boG_{ij,t}^{-1}(\wtz,\lambda)\Phi_{i,t}(\wtz,\lambda).$$
\begin{claim}
\label{claim-transition}
There exists uniform constants $C_1$ and $C_2(r)$ such that
for all $z\in\boA_{ij,t,r}$ and $t$ small enough:
\begin{equation}
\label{eq-transition-1}
||\boH_{ij,t}(z,\cdot)-\PhiS(\pi_{ij},\cdot)||_{\boW}\leq C_1r+C_2(r)|t\log t|.
\end{equation}
\end{claim}
Proof: we prove that Estimate \eqref{eq-transition-1} holds on the two boundary
components of $\boA_{ij,t,r}$ and we conclude by the maximum principle.
\begin{enumerate}[a)]
\item
We have in the subset $0\leq \Im(\wtz)\leq 2\pi$ of $\wtboA_{ij,t,r}$:
\begin{equation}
\label{eq-transition-2}
||\boG_{ij,t}(\wtz,\cdot)-I_2||_{\boW}\leq C(r)|t\log t|.
\end{equation}
Consider the segment
$$J_1=[\log r,\log r+2\pi\ii]\subset\wtboA_{ij,t,r}$$
which projects onto the circle $|v_{ij}|=r$.
By Equation \eqref{eq-spherical-3}, we have on $J_1$:
$$||\Phi_{i,t}(\wtz,\cdot)-\PhiS(p_{ij}(t)+e^{\wtz},\cdot)||_{\boW}\leq C(r)t.$$
We have on the circle $|v_{ij}|=r$:
$$||\PhiS(z,\cdot)-\PhiS(\pi_{ij},\cdot)||_{\boW}\leq Cr.$$
Hence, on the segment $J_1$:
\begin{equation}
\label{eq-transition-4}
||\boG_{ij,t}(\wtz,\cdot)^{-1}\Phi_{i,t}(\wtz,\cdot)-\PhiS(\pi_{ij},\cdot)||_{\boW}
\leq Cr+C(r)|t\log t|.
\end{equation}
So Estimate \eqref{eq-transition-1} holds on the circle $|v_{ij}|=r$.
\item By Points 1 and 3 of the Proof of Proposition \ref{prop-Pij},
$\boP(\xi_t,\alpha_{ij}\beta_{ij})$ extends at $t=0$ as a smooth function of $t$
and $t\log t$, with value at $t=0$:
\begin{eqnarray*}
\boP(\xi_0,\alpha_{ij}\beta_{ij})&=&\boP(\xiS,0_i,\pi_{ij})\boP(\xiC\cdot G_{ij},q_{ij},q_{ij}+r)\\
&=&\PhiS(\pi_{ij},\cdot)[\PhiC(q_{ij})G_{ij}(q_{ij},\cdot)]^{-1}\PhiC(q_{ij}+r)G_{ij}(q_{ij}+r,\cdot).
\end{eqnarray*}
Hence
\begin{equation}
\label{eq-transition-3}
||\boP(\xi_t,\alpha_{ij}\beta_{ij})-\PhiS(\pi_{ij},\cdot)||_{\boW}\leq Cr+C(r)|t\log t|.
\end{equation}
Consider the segment
$$J_2=[\log|t_{ij}|-\log r,\log|t_{ij}|-\log r+2\pi\ii]\subset\wtboA_{ij,t,r}$$
which projects onto the circle
$|v_{ij}|=\frac{|t_{ij}|}{r}$, identified with the circle $|w_{ij}|=r$.
On this circle, $\xi_t$ depends smoothly on $t$ and $t\log t$ so
$$||\Phi_{i,t}(\wtz,\cdot)-\boP(\xi_t,\alpha_{ij}\beta_{ij})||_{\boW}=||\Phi_{i,t}(\wtz,\cdot)-\Phi_{i,t}(\wtbeta_{ij}(1),\cdot)||_{\boW}\leq C(r)|t\log t|.$$
Hence by Estimate \eqref{eq-transition-3}:
$$||\Phi_{i,t}(\wtz,\cdot)-\PhiS(\pi_{ij},\cdot)||_{\boW}\leq Cr+C(r)|t\log t|
\quad\mbox{ in $J_2$.}$$
Using Estimate \eqref{eq-transition-2}, we conclude that
Estimate \eqref{eq-transition-4} holds on $J_2$, so
Estimate \eqref{eq-transition-1} holds on the circle $|w_{ij}|=r$.
By the maximum principle, Estimate \eqref{eq-transition-1} is true for all $z\in\boA_{ij,t,\varepsilon}$.\cqfd
\medskip

Returning to the proof of Proposition \ref{prop-transition}, let $(F_{i,t},B_{i,t})$ be the Iwasawa decomposition of $\Phi_{i,t}$.
By Claim \eqref{claim-transition} and Estimate \eqref{eq-transition-1}, we have in the
domain $0\leq \Im(\wtz)\leq 2\pi$ of $\wtboA_{ij,t,r}$:
$$||\Phi_{i,t}(\wtz,\cdot)-\PhiS(\pi_{ij},\cdot)||_{\boW}\leq Cr+C(r)|t\log t|.$$
This implies
$$||F_{i,t}(\wtz,1)-\FS(\pi_{ij},1)||\leq Cr+C(r)|t\log t|.$$
By Equation \eqref{eq-spherical-0}, we have
$$N_t(p_{ij}(t)+e^{\wtz})=\Nor(F_{i,t}(\wtz,1))$$
so
$$||N_t(z)-\NS(\pi_{ij})||\leq Cr+C(r)|t\log t|\quad\mbox{ in $\boA_{ij,t,r}$}.$$
Since $\NS(\pi_{ij})=-\u_{ij}$, Point 1 of Proposition \ref{prop-transition} follows.
\end{enumerate}
\item
Fix $t>0$ and let $\pi$ be the projection on the plane orthogonal to $\u_{ij}(0)$.
By Point 1, $\pi\circ f_t$ is a local diffeomorphism on $\boA_{ij,t,r}$.
To prove that it is a global diffeomorphism onto its image, we use a topological 
argument.
Let $D_1$ and $D'_1$ be the disks in $\C_i$ with center $p_{ij}(t)$ and respective
radius $r$ and $\frac{r}{2}$.
Let $A_1$ be the annulus $D_1\setminus\overline{D'_1}$.
By Proposition \ref{prop-spherical}, $\pi\circ f_t$ is a global diffeomorphism from $A_1$ onto its image. Moreover, it maps the outside circle $\partial D_1$
to the outside boundary component of the image.
Therefore, we may extend $\pi\circ f_t$ to a homeomorphism $h_1$ from the Riemann sphere minus $\overline{D'_1}$ to the Riemann sphere minus the disk bounded
by $\pi\circ f_t(\partial D'_1)$.
\medskip

Let $D_2$ and $D'_2$ be the disks in $\C_i$ with center $p_{ij}(t)$ and respective
radius $\frac{2|t_{ij}|}{r}$ and $\frac{|t_{ij}|}{r}$.
Then $\boA_{ij,t,r}=D_1\setminus D_2'$.
Let $A_2$ be the annulus $D_2\setminus\overline{D'_2}$,
which is identified with the
annulus $\frac{r}{2}<|w_{ij}|<r$ in $\C_{ij}$.
By Proposition \ref{prop-catenoidal}, $\pi\circ f_t$ is a global diffeomorphism from $A_2$ onto its image.
Moreover, it maps the inside circle $\partial D'_2$ onto the inside boundary component
of the image.
Therefore, we may extend $\pi\circ f_t$ to a homeomorphism $h_2$ from
the disk $D_2$ to the disk bounded by $\pi\circ f_t(\partial D_2)$.
We define a local homeomorphism $h:\CC\to\CC$ by
$$h=\left\{\begin{array}{l}
h_1\mbox{ in }\CC\setminus D'_1\\
\pi\circ f_t\mbox{ in } D_1\setminus D_2'\\
h_2\mbox{ in } D_2\end{array}\right.$$
Since the Riemann sphere is compact, $h$ is a covering map, and since it is simply connected, $h$ is a homeomorphism. Hence the restriction of $h$ to
$D_1\setminus D_2'$ is a homeomorphism onto its image. In other words,
$\pi\circ f_t$ is a diffeomorphism from $\boA_{ij,t,r}$ onto its image, which
proves Point 2.
\item Assume that $\tau_{ij}>0$.
We use barrier arguments to prove that the images of $\boA_{ij,t,r}$ and $\boA_{ji,t,r}$ are disjoint.
Let $\gamma\subset\C_{ij}$ be the circle whose image by
the catenoidal immersion $\Psi_{ij}$ is its waist circle.
Let $\omega_t$ be the center of mass of $f_t(\gamma)$.
Let $\Delta_t$ be the half-line issued from $\omega_t$ and containing $\v_i(t)$.
By Proposition \ref{prop-catenoidal}, for small $t$, $f_t(\gamma)$ is at distance $o(t)$
from the circle of center $\omega_t$ and radius $4t\tau_{ij}$ contained in the plane
orthogonal to $\Delta_t$.
Let $(x,y,z)$ be an orthonormal coordinate system in $\R^3$ with origin at $\omega_t$
and such that $\Delta_t$ is the positive $x$-axis.
Let $\delta_t$ be the distance between $\omega_t$ and $\v_i(t)$, so
$\v_i(t)=(\delta_t,0,0)$ in the $(x,y,z)$ coordinate system.
\begin{claim}
\label{claim-barrier}
There exists a uniform constant $C$ such that for $t$ small enough,
$$\delta_t\geq 1+2\tau_{ij}|t\log t|-Ct.$$
\end{claim}
Proof: let $D_t$ be a half-period of the Delaunay surface of neck-radius $2t\tau_{ij}$
and axis $\Delta_t$, bounded by a circle of minimum radius on the left (a neck)
and a circle of maximum radius on the right (a bulge). (Here the words ``left'' and
``right'' refer to the $x$-axis.)
Slide $D_t$ from the right until a first contact point $p$ between $D_t$ and $A_t$
occurs.
Let $d_t^-$ and $d_t^+$ be the $x$-coordinate of the centers of the left and right
boundary circles of $D_t$.
Observe that since $\tau_{ij}>0$, the Gauss maps of $D_t$ and $A_t$ both point
to the ``inside''. By the maximum principle for CMC-1 surfaces, the first contact point $p$ must
be on the boundary of $A_t$ or $D_t$. It cannot be on the left boundary circle
of $D_t$ (which is too small) nor on the right boundary circle (which is too big).
So it has to be on the right boundary of $A_t$, namely on the image of the circle
$|v_{ij}|=r$.
By Proposition \ref{prop-spherical}, the point $p$ is at distance less than $Ct$ from
the unit sphere centered at $(\delta_t,0,0)$, and is outside the cylinder of axis
$\Delta_t$ and radius $\frac{r}{4}$.
Outside of this cylinder, for $t$ small enough, the half-Delaunay surface $D_t$ is at distance less than $Ct$ from the sphere
centered at $(d_t^+,0,0)$. Hence
\begin{equation}
\label{eq-barrier-1}
|\delta_t-d_t^+|\leq Ct.
\end{equation}
Let $(d_t',0,0)$ be the center of the circle of radius $8t\tau_{ij}$ on $D_t$.
Because $p$ is the first contact point, we have $d_t'>0$.
Since $\frac{1}{t}(D_t-(d_t^-,0,0))$ converges to a catenoid, we have
\begin{equation}
\label{eq-barrier-2}
|d_t'-d_t^-|\leq Ct.
\end{equation}
Finally, it is known that the half-period of the Delaunay surface of necksize 
$2t\tau_{ij}$ satisfies:
\begin{equation}
\label{eq-barrier-3}
d_t^+-d_t^-=1+2\tau_{ij}|t\log t|+o(t\log t).
\end{equation}
Claim \ref{claim-barrier} follows from Estimates \eqref{eq-barrier-1},
\eqref{eq-barrier-2}, \eqref{eq-barrier-3} and $d_t'>0$.
\cqfd
\medskip

Returning to the proof of Point 3 of Proposition \ref{prop-transition},
let $S$ be the hemisphere with center at $(0,0,1)$ and contained in the half-space
$x\leq 1$.
By Claim \ref{claim-barrier}, the right boundary component of $f_t(\boA_{ij,t,r})$
(namely the image of the circle $|v_{ij}|=r$) lies on the right of $S$.
By Proposition \ref{prop-catenoidal}, the left boundary component of
$f_t(\boA_{ij,t,r})$ (namely the image of the circle $|w_{ij}|=r$) also lies
on the right of $S$. By the maximum principle, $f_t(\boA_{ij,t,r})$ lies on the
right of $S$. In particular, $f_t(\boA_{ij,t,r})$ is in the half-space $x>0$.
By the same argument, $f_t(\boA_{ji,t,r})$ is in the half-space $x<0$,
so $f_t(\boA_{ij,t,r})$ and $f_t(\boA_{ji,t,r})$ are disjoint.\cqfd
\end{myenumerate}
\subsection{Embeddedness}
Assume that the weighted graph $\Gamma$ is pre-embedded (see Definition
\ref{def-preembedded}).
Fix a small positive $r$ such that $2r\leq\epsilon$, where $\epsilon$ is the number given by Proposition \ref{prop-delaunay}.
By Propositions \ref{prop-spherical}, \ref{prop-delaunay} and
\ref{prop-transition}, the $f_t$-images of the following domains in $\Sigma_t$
are embedded for $t>0$ small enough:
\begin{itemize}
\item the domains $\Omega_{i,r}\cup\{\infty_i\}$ for $i\in[1,N]$,
\item the domains $\Omega_{ij,r}\cup\{\infty_{ij}\}$ for $(i,j)\in I$,
\item the punctured disks $D^*(p_k^0(t),2r)$
for $i\in[1,N]$ and $k\in R_i$,
\item the annuli $\boA_{ij,t,2r}$ for $(i,j)\in I\cup I^*$.
\end{itemize}
These domains cover $\Sigma_t$ and we have good control on the position of their images.
By Claim \ref{claim-barrier}, if $\v_i$ and $\v_j$ are adjacent, we have
$$||\v_i(t)-\v_j(t)||\geq 2+4\tau_{ij}|t\log t|-Ct$$
so by Proposition \ref{prop-spherical}, the $f_t$-images of
$\Omega_{i,r}$ and $\Omega_{j,r}$ are disjoint.
If $\v_i$ and $\v_j$ are not adjacent, then 
By Point 2 of Definition \ref{def-preembedded} and Proposition \ref{prop-spherical}, the $f_t$-images of $\Omega_{i,r}$ and $\Omega_{j,r}$ are disjoint.
Points 2 and 3 of Definition \ref{def-preembedded} and Proposition \ref{prop-delaunay}
ensure that for $i\in[1,N]$ and $k\in R_i$, $f_t(D^*(p_k^0(t),2r))$ only intersects
$f_t(\Omega_{i,r})$.
It is then rather clear that $M_t=f_t(\Sigma_t)$ is embedded. A formal proof can be written
by covering $M_t$ by a finite number of open sets in $\R^3$ whose pre-image by
$f_t$ is included in one of the domains in the above list. See Section 10.2 of \cite{nnoids}
where the complete argument is given in the case $N=1$.
This concludes our proof of Theorem \ref{thm-main}.
\appendix
\section{On removable singularities}
\label{appendix-regularity}
In this appendix, we are interested in the following question:
let $\xi$ be a DPW potential in the punctured disk $D^*(p,r)=\{z\in\C:0<|z|<r\}$
with a pole at $p$.
Assume that we have a solution $\Phi$ in the universal cover of $D^*(p,r)$
which solves the Monodromy Problem \eqref{pb-monodromy}.
How does that help us to prove that $p$ is a removable singularity of
$\xi$ in the sense of Definition \ref{def-removable}?
The question is certainly very vague, but as we shall see, in certain cases,
provided the Monodromy Problem is solved,
the Regularity Problem boils down to a finite number of real equations.
The results that we present here are tailored to our needs in this paper,
but the question deserves further investigation.
\medskip

In this appendix, we consider families of potentials $\xi_t$ in a domain in the complex plane
and solutions $\Phi_t$ of $d\Phi_t=\Phi_t\xi_t$ in its universal cover, both depending
on $t\in(-\epsilon,\epsilon)$. We assume that
the map $(t,z)\mapsto \xi_t(z,\cdot)$ is $C^1$ with values in $\sl(2,\boW)$ and
the map $(t,z)\mapsto \Phi_t(z,\cdot)$ is $C^1$ with values in $SL(2,\boW)$
(this is the meaning of $C^1$ in the statements).
Our first result concerns potentials which are meromorphic perturbations of the
catenoid potential $\xiC$ (see Section \ref{section-background-DPW-catenoid}).
\begin{theorem}
\label{thm-regularity-xiC}
Let $\xi_t$ be a $C^1$ family of DPW potentials on $D^*(0,R)$ and $\Phi_t$ a $C^1$ family of
solutions of $d\Phi_t=\Phi_t\xi_t$ in its universal cover $\widetilde{D^*}(0,R)$.
Let $\wtz_0\in\widetilde{D^*}(0,R)$ be an arbitrary point in $\widetilde{D^*}(0,R)$.
Assume that:
\begin{myenumerate}
\item $\xi_0=\xiC$.
\item For $t\neq 0$, $\xi_t$ has at most a triple pole at $z=0$, with principal part
$$\xi_t(z,\lambda)=\matrix{0&\lambda^{-1}\\0&0}\left(\frac{a_t}{z^3}+\frac{b_t}{z^2}+\frac{c_t}{z}\right)dz+\Xi_t(z,\lambda)$$
where $\Xi_t$ is holomorphic in $D(0,R)$.
\item $\Phi_0(\wtz_0,\cdot)$ is holomorphic in the annulus $\A_{\rho'}$ for some
$\rho'>\rho$.
\label{hyp-thm-regularity-4}
\item For all $t\neq 0$, $\Phi_t$ solves the Monodromy Problem \eqref{pb-monodromy}.
\item For all $t\neq 0$, $a_t^0=0$ and $\Re(b_t^0)=0$.
\end{myenumerate}
Then for $t$ small enough, $a_t=b_t=c_t=0$, so $\xi_t$ is holomorphic at $z=0$.
\end{theorem}
Proof: let $z_0$ be the projection of $\wtz_0$ in $D^*(0,R)$.
Let $\gamma$ be a generator of $\pi_1(D^*(0,R),z_0)$ and
$\sigma$ be the corresponding Deck transformation of $\widetilde{D^*}(0,R)$.
Let $(F_0,B_0)$ be the Iwasawa decomposition of $\Phi_0(\wtz_0,\cdot)\PhiC(z_0)^{-1}$.
By hypothesis \ref{hyp-thm-regularity-4}, $B_0$ extends holomorphically to $\D_{\rho'}$
with $\rho'>\rho$ so is in $SL(2,\boW^{\pos})$.
Replacing $\Phi_t$ by $F_0^{-1}\Phi_t$ for all $t$, we can assume without loss of generality that $F_0=I_2$.  (This does not change the hypothesis that the Monodromy Problem is solved since $F_0\in\Lambda SU(2)$). Then
\begin{equation}
\label{eq-regularity-xiC-1}
\Phi_0(\wtz_0,\cdot)\PhiC(z_0)^{-1}=B_0\in\Lambda_+^{\R}SL(2,\C).
\end{equation}
For $\x=(a,b,c)\in(\boW^\pos)^3$, let $\xi_{t,\x}$ be the potential in $D^*(0,R)$ defined
by
$$\xi_{t,\x}(z,\lambda)=\matrix{0&\lambda^{-1}\\0&0}\left(\frac{a}{z^3}+\frac{b}{z^2}+\frac{c}{z}\right)dz+\Xi_t(z,\lambda).$$
Let $\Phi_{t,\x}$ be the solution of $d\Phi_{t,\x}=\Phi_{t,\x}\xi_{t,\x}$
in $\widetilde{D^*}(0,R)$ with initial condition
$$\Phi_{t,\x}(\wtz_0,\lambda)=\Phi_t(\wtz_0,\lambda).$$
We consider the following Problem:
\begin{equation}
\label{pb-regularity-xiC}
\left\{\begin{array}{l}
\mbox{$\Phi_{t,\x}$ solves the Monodromy Problem \eqref{pb-monodromy}}\\
a^0=0\\
\Re(b^0)=0\end{array}\right.\end{equation}
Writing $\x_t=(a_t,b_t,c_t)$,
we have $\xi_t=\xi_{t,\x_t}$ and $\Phi_t=\Phi_{t,\x_t}$.
So Theorem \ref{thm-regularity-xiC} follows from the following
\begin{lemma}
For $(t,\x)$ in a neighborhood of $(0,0)$, Problem \eqref{pb-regularity-xiC} is equivalent to $\x=0$.
\end{lemma}
Proof: we use an Implicit Function argument.
Let
$M(t,\x)=\boM(\Phi_{t,\x},\sigma)$ be the monodromy of $\Phi_{t,\x}$ with respect
to $\sigma$.
\begin{claim}\begin{myenumerate}
\item $(t,\x)\mapsto M(t,\x)$ is a $C^1$ map from $(-\epsilon,\epsilon)\times(\boW^\pos)^3$ to
$SL(2,\boW)$.
\item For all $t\in (-\epsilon,\epsilon)$, $M(t,0)=I_2$.
\item The partial differential of $M$ with respect to $\x$ at $(0,0)$ is given by
$$d_{\x}M(0,0)=\frac{2\pi\ii}{\lambda}B_0\matrix{-db&dc\\-da&db}B_0^{-1}.$$
\end{myenumerate}
\end{claim}
Proof:
\begin{myenumerate}
\item $(t,\x)\mapsto \xi_{t,\x}$ is a $C^1$ map with value in $\su(2,\boW)$ so Point 1 follows
from standard ODE theory.
\item Point 2 follows from the fact that $\xi_{t,0}=\Xi_t$ is holomorphic in $D(0,R)$.
\item Let $\wtgamma$ be the lift of $\gamma$ to $\widetilde{D^*}(0,R)$ such that
$\wtgamma(0)=\wtz_0$. By Equation \ref{eq-boM-boP}, we have
$$M(t,\x)=\Phi_t(\wtz_0,\cdot)\boP(\xi_{t,\x},\gamma)\Phi_t(\wtz_0,\cdot)^{-1}.$$
Since $\boP(\xi_{0,0},\gamma)=I_2$,
$$d_{\x}M(0,0)=\Phi_0(\wtz_0,\cdot)d_{\x}\boP(\xi_{t,\x},\gamma)|_{(t,\x)=(0,0)}\Phi_0(\wtz_0,\cdot)^{-1}.$$
By Proposition \ref{prop-principal-derivee}, since $\xi_{0,0}=\xiC$:
$$d_{\x}\boP(\xi_{t,\x},\gamma)|_{(t,\x)=(0,0)}=
\int_{\wtgamma}\left[\PhiC(z_0)^{-1}\PhiC\right]d_{\x}\xi_{t,\x}|_{(t,\x)=(0,0)}
\left[\PhiC(z_0)^{-1}\PhiC\right]^{-1}.$$
Hence using Equation \eqref{eq-regularity-xiC-1} and the Residue Theorem:
\begin{eqnarray*}
d_{\x}M(0,0)
&=&\int_{\gamma} B_0\PhiC d_{\x}\boP(\xi_{t,\x},\gamma)|_{(t,\x)=(0,0)}(\PhiC)^{-1}
B_0^{-1}\\
&=&2\pi\ii\, B_0\Res_0\left[\matrix{1&0\\z&1}\matrix{0&\lambda^{-1}\\0&0}
\matrix{1&0\\-z&1}\left(\frac{da}{z^3}+\frac{db}{z^2}+\frac{dc}{z}\right)\right]B_0^{-1}\\
&=&\frac{2\pi\ii}{\lambda}B_0\matrix{-db&dc\\-da&db}B_0^{-1}.
\end{eqnarray*}
\cqfd\end{myenumerate}

We define for $(t,\x)$ in a neighborhood of $(0,0)$:
$$\wtM(t,\x)=\boL_1(\lambda\log M(t,\x))$$
$$\whM(t,\x)=\boL_1(\wtM(t,\x))$$
where $\boL_1$ is the operator introduced in Proposition \ref{prop-boL}.
Explicitly:
$$\wtM(t,\x)(\lambda)=\frac{1}{\lambda-1}\left(\lambda\log M(t,\x)(\lambda)
-\log M(t,\x)(1)\right)$$
$$\whM(t,\x)(\lambda)=\frac{1}{\lambda-1}\left(\wtM(t,\x)(\lambda)-\wtM(t,\x)(1)\right).$$
\begin{claim}
\label{claim-regularity-xiC-2}
\begin{myenumerate}
\item $\wtM$ and $\whM$ are $C^1$ maps from
a neighborhood of $(0,0)$ in $\R\times(\boW^\pos)^3$ to $\sl(2,\boW)$.
\item The Monodromy Problem \eqref{pb-monodromy} for $\Phi_{t,\x}$ is equivalent to:
$$\left\{\begin{array}{l}
\whM(t,\x)\in\Lambda\su(2)\\
M(t,\x)(1)=I_2\\
\wtM(t,\x)(1)=0.\end{array}\right.$$
\end{myenumerate}
\end{claim}
Proof:
\begin{myenumerate}
\item Point 1 follows from Proposition \ref{prop-boL} and the fact that $\boW$ is a Banach algebra.
\item Using $d\log(I_2)=\mbox{id}$, we have
$$M(t,\x)(1)=I_2\Rightarrow \wtM(t,\x)(1)=\frac{\partial}{\partial\lambda} M(t,x)(1)$$
$$M(t,\x)(1)=I_2\mbox{ and } \wtM(t,\x)(1)=0\Rightarrow
\whM(t,\x)(\lambda)=\frac{\lambda}{(\lambda-1)^2}\log M(t,\x)(\lambda).$$
Point 2 follows from the fact that $(\lambda-1)^2\lambda^{-1}$ is real on the unit circle.
\cqfd
\end{myenumerate}
We introduce the auxiliary variables $(p,q,r)$ defined in function of $(a,b,c)$ by
$$\matrix{-q&r\\-p&q}=B_0\matrix{-b&c\\-a&b}B_0^{-1}.$$
This change of variable is an automorphism of $(\boW^{\pos})^3$ because
$B_0\in SL(2,\boW^{\pos})$.
Then
\begin{equation}
\label{eq-dxM00}
d_{\x}M(0,0)=\frac{2\pi\ii}{\lambda}\matrix{-dq&dr\\-dp&dq}.
\end{equation}
Writing $B_0(0)=\minimatrix{\rho&\mu\\0&\rho^{-1}}$ with $\rho>0$ and $\mu\in\C$, we have
$$\matrix{-q^0&r^0\\-p^0&q^0}=
\matrix{-b^0-\frac{\mu}{\rho}a^0&\mu^2a^0+2\mu\rho b+\rho^2c\\-\frac{1}{\rho^2}a^0&b^0+\frac{\mu}{\rho}a^0}.$$
Hence
\begin{equation}
\label{eq-aequivp}
\left\{\begin{array}{l}
a^0=0\\
\Re(b^0)=0\end{array}\right.\Leftrightarrow
\left\{\begin{array}{l}
p^0=0\\
\Re(q^0)=0\end{array}\right.\end{equation}
Using Proposition \ref{prop-boL} twice, we decompose an arbitrary parameter $x\in\boW^\pos$ as
$$x(\lambda)=x(1)+(\lambda-1)x'(1)+(\lambda-1)^2\widetilde{x}(\lambda)\quad
\mbox{with $\widetilde{x}\in\boW^\pos$.}$$
Then by Equation \eqref{eq-dxM00}:
\begin{equation}
\label{eq-dxwtM00}
d_{\x}\wtM(0,0)(\lambda)=2\pi\ii\matrix{-dq'(1)-(\lambda-1)d\wtq(\lambda)&
dr'(1)+(\lambda-1)d\wtr(\lambda)\\-dp'(1)-(\lambda-1)d\wtp(\lambda)&dq'(1)+
(\lambda-1)d\wtq(\lambda)}
\end{equation}
\begin{equation}
\label{eq-dxwhM00}
d_{\x}\whM(0,0)=2\pi\ii\matrix{-d\wtq&d\wtr\\-d\wtp&d\wtq}.
\end{equation}
We define
$$\boE_1(t,\x)=\whM_{11}(t,\x)+\whM_{11}(t,\x)^*\in\boW$$
$$\boE_2(t,\x)=\whM_{12}(t,\x)+\whM_{21}(t,\x)^*\in\boW$$
$$\boE_3(t,\x)=\left(M_{11}(t,\x)(1)-1,M_{12}(t,\x)(1),M_{21}(t,\x)(1)\right)\in\C^3$$
$$\boE_4(t,\x)=\left(\wtM_{11}(t,\x)(1),\wtM_{12}(t,\x)(1),\wtM_{21}(t,\x)(1)\right)\in\C^3$$
and finally
$$\boF(t,\x)=\left[\boE_1^+,\boE_2^+,(\boE_2^-)^*,\boE_2^0,\boE_3,\boE_4,p^0,4\pi\Re(q^0)+\ii \Re(\boE_1^0)\right](t,\x)
\quad\in(\boW^+)\times\C^9.$$
By definition, $\boE_1=\boE_1^*$ so $\boE_1=0$ is equivalent to
$\boE_1^+=0$ and $\Re(\boE_1^0)=0$.
By Claim \ref{claim-regularity-xiC-2} and Equation \eqref{eq-aequivp}
Problem \eqref{pb-regularity-xiC} is equivalent to $\boF(t,\x)=0$.
We have $\boF(0,0)=0$. Using Equation \eqref{eq-dxwhM00}:
\begin{equation}
\label{eq-xic-1}
d\boE_1(0,0)^+=-2\pi\ii \,d\wtq^{\,+}
\end{equation}
\begin{equation}
\label{eq-xic-2}
d\boE_2(0,0)^+=2\pi\ii \,d\wtr^{\,+}
\end{equation}
\begin{equation}
\label{eq-xic-3}
(d\boE_2(0,0)^-)^*=-2\pi\ii \,d\wtp^{\,+}
\end{equation}
\begin{equation}
\label{eq-xic-4}
dp^0=dp(1)-dp'(1)+d\wtp^{\,0}
\end{equation}
$$dq^0=dq(1)-dq'(1)+d\wtq^{\,0}$$
\begin{equation}
\label{eq-xic-6}
4\pi \Re(dq^0)+\ii \Re(d\boE_1(0,0)^0)=4\pi\left(\Re(dq(1))-\Re(dq'(1))+d\wtq^{\,0}\right)
\end{equation}
\begin{equation}
\label{eq-xic-7}
d\boE_2(0,0)^0=2\pi\ii(d\wtr^{\,0}+d\overline{\wtp{\,^0}}).
\end{equation}
Using Equations \eqref{eq-dxM00} and \eqref{eq-dxwtM00}:
\begin{equation}
\label{eq-xic-8}
d\boE_3(0,0)=2\pi\ii(-dq(1),dr(1),-dp(1))
\end{equation}
\begin{equation}
\label{eq-xic-9}
d\boE_4(0,0)=2\pi\ii(-dq'(1),dr'(1),-dp'(1))
\end{equation}
\begin{claim}
Let $L$ be the partial differential of $\boF$ at $(0,0)$ with respect to
$$(\wtp^+,\wtq^+,\wtr^+,p(1),q(1),r(1),p'(1),q'(1),r'(1),\wtp^0,\wtq^0,\wtr^0).$$
Then $L$ is an ($\R$-linear) automorphism of $(\boW^+)^3\times\C^9$.
\end{claim}
Proof: the partial differential of $(\boE_1^+,\boE_2^+,(\boE_2^-)^*)$ with respect to
$(\wtp^+,\wtq^+,\wtr^+)$ is clearly an automorphism of $(\boW^+)^3$. By Proposition \ref{prop-fredholm}, it suffices to
prove that $L$ is injective, so let us formally solve $dL=0$.
\begin{itemize}
\item Equation \eqref{eq-xic-8} gives $dp(1)=dq(1)=dr(1)=0$.
\item Equation \eqref{eq-xic-9} gives $dp'(1)=dq'(1)=dr'(1)=0$.
\item Equations \eqref{eq-xic-4}, \eqref{eq-xic-6} and \eqref{eq-xic-7} give $d\wtp^{\,0}=d\wtq^{\,0}=d\wtr^{\,0}=0$.
\item Equations \eqref{eq-xic-1}, \eqref{eq-xic-2} and \eqref{eq-xic-3} give
$d\wtp^{\,+}=d\wtq^{\,+}=d\wtr^{\,+}=0$.\cqfd
\end{itemize}
By the Implicit Function Theorem, for $(t,\x)$ in a neighborhood of $(0,0)$, Problem
\eqref{pb-regularity-xiC} uniquely determines $\x$ as a function of $t$. Since we know that $\x=0$
is a solution, it has to be the only one, so Problem \eqref{pb-regularity-xiC} is equivalent to
$\x=0$.
\cqfd
\medskip

Using Theorem \ref{thm-regularity-xiC}, we now prove a result for perturbations of the standard catenoid
potential in a neighborhood of $z=\infty$. Let $D^*(\infty,R)=\{z\in\C,|z|>R^{-1}\}$,
$\widetilde{D^*}(\infty,R)$ its universal cover and $\wtz_0\in\widetilde{D^*}(\infty,R)$
a base point.
\begin{theorem}
\label{thm-regularity-infty}
Let
$\xi_t=\minimatrix{\alpha_t& \lambda^{-1}\beta_t\\\gamma_t&-\alpha_t}$ be a $C^1$
family of DPW potentials in $D^*(\infty,R)$ with a pole of multiplicity at most $\minimatrix{1&1\\2&1}$ at $\infty$ and such that $\xi_0=\minimatrix{0&0\\1&0}dz$.
Let $\Phi_t$ be a $C^1$ family of solutions of $d\Phi_t=\Phi_t\xi_t$.
Assume that $\Phi_0(\wtz_0,\cdot)$ is holomorphic in the annulus $\A_{\rho'}$ for some
$\rho'>\rho$ and that for all $t$:
$$\left\{\begin{array}{l}
\mbox{$\Phi_t$ solves the Monodromy Problem with respect to $\sigma$}\\
\Res_{\infty}\beta_t^0=0\\
\Re\left(\Res_{\infty}(z\beta_t^0)\right)=0
\end{array}\right.$$
Then for $t$ small enough, there exists a gauge $G_t$ such that $\xi_t\cdot G_t$
extends holomorphically at $\infty$. In other words, $\infty$ is a removable singularity.
\end{theorem}
Proof: 
we consider the gauged potential
\begin{equation}
\label{eq-regularity-infty-whxit}
\whxi_t=\matrix{\whalpha_t&\lambda^{-1}\whbeta_t\\\whgamma_t&-\whalpha_t}=\xi_t\cdot G_t\quad\mbox{ with }\quad G_t(z,\lambda)=\matrix{\frac{1}{z}&g_t(\lambda)\\0&z}.
\end{equation}
Here $g_t\in\boW^{\pos}$ is a parameter to be determined.
A computation gives
$$\whxi_t=
\matrix{\alpha_t-g_tz^{-1}\gamma_t-\frac{dz}{z}&
2g_tz\alpha_t+\lambda^{-1}z^2\beta_t-g_t^2\gamma_t-g_tdz\\
z^{-2}\gamma_t&
-\alpha_t+g_tz^{-1}\gamma_t+\frac{dz}{z}}.$$
This implies that $\whxi_t$ has a pole of multiplicity at most $\minimatrix{1&3\\0&1}$
at $\infty$.
At $t=0$, we have
$$\Res_{\infty}(z^{-1}\gamma_0)=\Res_{\infty}\frac{dz}{z}=-1.$$
Hence $\Res_{\infty}(z^{-1}\gamma_t)$ does not vanish for $t$ small enough.
We take
$$g_t=\frac{1+\Res_{\infty}\alpha_t}{\Res_{\infty}(z^{-1}\gamma_t)}$$
so that $\Res_{\infty}\whalpha_t=0$ and $\whalpha_t$ is holomorphic at $\infty$.
We introduce the local coordinate $w=\frac{-1}{z}$ in a neighborhood of $\infty$. At $t=0$,
we have $g_0=-1$ and by Equation \eqref{eq-regularity-infty-whxit}:
$$\whxi_0=\matrix{0&0\\\frac{dz}{z^2}&0}=\matrix{0&0\\dw&0}.$$
We write the principal part of $\whbeta_t$ at $\infty$ as follows:
$$\whbeta_t=\left(\frac{a_t}{w^3}+\frac{b_t}{w^2}+\frac{c_t}{w}+O(w^0)\right)dw.$$
By Proposition \ref{prop-gauging}, we have
$\beta_t^0=w^2\whbeta_t^0$ so
$$a_t^0=\Res_{\infty}\beta_t^0=0$$
$$\Re(b_t^0)=\Re(\Res_{\infty}(w^{-1}\beta_t^0))=-\Re(\Res_{\infty}(z\beta_t^0))=0.$$
Let $\whPhi_t=\Phi_tG_t$. Then $\whPhi_t$ solves the Monodromy Problem
and since $g_0=-1$, $\whPhi_0(\wtz_0,\cdot)$ is holomorphic in $\boA_{\rho'}$.
Applying Theorem \ref{thm-regularity-xiC} to $(\whxi_t,\whPhi_t)$, we obtain $(a_t,b_t,c_t)=(0,0,0)$.
Hence $\whxi_t$ is holomorphic at $\infty$.\cqfd
\medskip

By duality (see Section \ref{section-background-DPW-duality}), we obtain the following result for perturbations of the standard spherical potential:
\begin{corollary}
\label{cor-regularity-infty}
Let $\xi_t=\minimatrix{\alpha_t& \lambda^{-1}\beta_t\\\gamma_t&-\alpha_t}$ be a $C^1$ family of DPW potentials in $D(\infty,R)$ with a pole of multiplicity at most $\minimatrix{1&2\\1&1}$ at $\infty$ and such that $\xi_0=\minimatrix{0&\lambda^{-1}\\0&0}dz$.
Let $\Phi_t$ be a $C^1$ family of solutions of $d\Phi_t=\Phi_t\xi_t$.
Assume that $\Phi_0(\wtz_0,\cdot)$ is holomorphic in the annulus $\A_{\rho'}$ for some
$\rho'>\rho$ and that for all $t$:
$$\left\{\begin{array}{l}
\mbox{$\Phi_t$ solves the Monodromy Problem with respect to $\sigma$}\\
\Res_{\infty}\gamma_t^0=0\\
\Re\left(\Res_{\infty}(z\gamma_t^0)\right)=0
\end{array}\right.$$
Then for $t$ small enough, there exists a gauge $G_t$ such that $\xi_t\cdot G_t$
extends holomorphically at $\infty$.
\end{corollary}
By duality, the gauge $G_t$ has the following form:
$$G_t=H^{-1}\matrix{\frac{1}{z}&g_t\\0&z}H=\matrix{z&0\\\lambda g_t&\frac{1}{z}}.$$
\section{Principal solution in an annulus}
\label{appendix-principal}
Fix some numbers $0<\varepsilon<\eta$.
For $t\in\C$ such that $0<|t|<\eta^2$, we define the annulus
$$\boA_t=\{z\in\C:\smallfrac{|t|}{\eta}< |z|<\eta\}$$
and let $\psi_t:\boA_t\to\boA_t$ be the diffeomorphism defined by $\psi_t(z)=\frac{t}{z}$.
What we have in mind here is the case of opening nodes, where $\boA_t$ would be the annulus
$\frac{|t_i|}{\varepsilon}<|v_i|<\varepsilon$ and $\psi_t$ would be the change of coordinate $w_i\circ v_i^{-1}$ (see Section \ref{section-background-openingnodes}).
\medskip

Consider a family of $\sl(n,\C)$-valued 1-forms $\xi_t(z)$ depending holomorphically
on $t\in D^*(0,\eta^2)$ and $z\in\boA_{t}$.
In this section, we are interested in the behavior as $t\to 0$
of the principal solution of $\xi_t$
on the path $s\mapsto\varepsilon^{1-2s}t^s$, which connects $z=\varepsilon$
to $z=\frac{t}{\varepsilon}$. This spiral depends on the choice of the argument of
$t$, so we consider the universal cover $\exp:\C\to\C^*$ and we write
$t=e^w$. For $\Re(w)<2\log\eta$,
let $\beta_w:[0,1]\to\boA_t$ be the spiral parametrized by
$$\beta_w(s)=\varepsilon^{1-2s}e^{sw}.$$
Let $\gamma:[0,1]\to\boA_t$ be the circle parametrized by $\gamma(s)=\varepsilon e^{2\pi\ii s}$.
For $w\in\C$ such that $\Re(w)<2\log\eta$, we set $t=e^w$ and we define
$$\wtF(w)=\boP(\xi_t,\gamma)^{-\frac{w}{2\pi\ii}} \boP(\xi_t,\beta_w)
=\exp\left(-\frac{w}{2\pi\ii}\log\boP(\xi_t,\gamma)\right)\boP(\xi_t,\beta_w).$$
\begin{theorem}
\label{thm-principal}
Let $\whxi_t=\psi_t^*\xi_t$.
Assume that there exists $\sl(n,\C)$-valued 1-forms $\xi_0$ and $\whxi_0$, holomorphic
in $D(0,\eta)$, such that on compact subsets of $D^*(0,\eta)$, we have:
$$
\lim_{t\to 0}\xi_t=\xi_0\quad\mbox{ and }\quad
\lim_{t\to 0}\whxi_t=\whxi_0.$$
Then for $|t|$ small enough:
\begin{myenumerate}
\item The function $\wtF(w)$ satisfies $\wtF(w+2\pi\ii)=\wtF(w)$ so
descends to a holomorphic function $F$ in $D^*(0,\eta)$, such that $\wtF(w)=F(e^w)$.
\item The function $F$ extends holomorphically at $t=0$ with
$$F(0)=\boP(\xi_0,\varepsilon,0)\boP(\whxi_0,0,\varepsilon).$$
\end{myenumerate}
\end{theorem} 
Here $\boP(\xi_0,\varepsilon,0)$ denotes the principal solution of $\xi_0$ on an arbitrary
path from $\varepsilon$ to $0$ in $D(0,\eta)$ (see Remark \ref{remark-principal}).
\medskip

Proof: we make the following change of variables:
$$\wtz=\varphi(z)=\varepsilon^{-1}z,\quad \wtt=\varepsilon^{-2}t\quad
\mbox{ and } \quad \wtw=w-2\log\varepsilon.$$
Let $\widetilde{\eta}=\frac{\eta}{\varepsilon}$.
Then
$$\frac{|t|}{\eta}<|z|<\eta\;\Leftrightarrow\;\frac{|\wtt|}{\widetilde{\eta}}<|\wtz|<\widetilde{\eta}$$
$$\varphi(\psi_t(z))=\varepsilon^{-1}\frac{\varepsilon^2\wtt}{\varepsilon\wtz}
=\frac{\wtt}{\wtz}$$
$$\varphi(\beta_w(s))=\varepsilon^{-1}\varepsilon^{1-2s}
e^{s(2\log\varepsilon+\wtw)}=e^{s\wtw}.$$
Thanks to this change of variables, it suffices to prove Theorem \ref{thm-principal}
in the case $\varepsilon=1$, which we assume from now on.
\medskip
We use the letter $C$ to denote uniform constants (independent of $z$ and $t$).
On the circle $|z|=1$, the map $t\mapsto \xi_t(z)$ is holomorphic in a punctured neighborhood of $t=0$ and extends continuously at $t=0$ so
is holomorphic in a neighborhood of $0$. Hence:
\begin{equation}
\label{eq-principal-1}
||\xi_t(z)-\xi_0(z)||\leq C|t|.
\end{equation}
Since $\boP(\xi_0,\gamma)=I_2$, Inequality \eqref{eq-principal-1} gives
\begin{equation}
\label{eq-principal-7}
||\boP(\xi_t,\gamma)-I_2||\leq C|t|.
\end{equation}
Consequently, for $t\neq 0$ small enough, $\log\boP(\xi_t,\gamma)$ is well defined
so $\wtF(w)$ is well defined.
The path $\beta_{w+2\pi\ii}$ is homotopic to the product $\gamma\beta_w$.
By Equation \eqref{eq-principal-morphism}, we have
$$
\wtF(w+2\pi\ii)=\boP(\xi_t,\gamma)^{-\frac{w+2\pi\ii}{2\pi\ii}}
\boP(\xi_t,\beta_{w+2\pi\ii})=
\boP(\xi_t,\gamma)^{-\frac{w}{2\pi\ii}}\boP(\xi_t,\gamma)^{-1}\boP(\xi_t,\gamma)\boP(\xi_t,\beta_w)=\wtF(w).$$
To prove Point 2, fix $t$ and write $t=e^w$ where $w$ is chosen so that
$|\Im(w)|\leq \pi$.
We split the path $\beta_w$ into $\beta_w=\alpha_w\whalpha_w^{-1}$
where
$$\alpha_w(s)=\beta_w(\smallfrac{s}{2})=e^{\frac{s}{2}w}\quad\mbox{ and }\quad
\whalpha_w(s)=\beta_w(1-\smallfrac{s}{2})=\psi_t(\alpha_w(s))\quad\mbox{ for $s\in[0,1]$}.$$
\begin{claim}
\label{claim-principal}
There exists a uniform constant $C$ such that for $t\neq 0$ small enough:
$$\int_0^1\left\|\xi_0(\alpha_w(s))\alpha'_w(s)\right\|ds\leq C$$
\begin{equation}
\label{eq-principal-claim-2}
\int_0^1\left\|\,\left[\xi_t(\alpha_w(s))-\xi_0(\alpha_w(s))\right]\alpha'_w(s)\right\|ds\leq C|t|^{1/2}.\end{equation}
\end{claim}
Proof: we have
$$\alpha_w'(s)=\frac{w}{2}e^{\frac{s}{2}w}.$$
Provided $|t|\leq e^{-\pi}$, we have:
$$-\Re(w)=-\log|t|\geq\pi$$
$$|w|\leq |\Re(w)|+|\Im(w)|\leq -\Re(w)+\pi\leq -2\Re(w)$$
\begin{equation}
\label{eq-principal-4}
|\alpha_w'(s)|\leq -\Re(w)e^{\frac{s}{2}\Re(w)}.
\end{equation}
Let
$$C_1=\max_{z\in D(0,1)}\left\|\frac{\xi_0(z)}{dz}\right\|.$$
Then
$$\int_0^1\left\|\xi_0(\alpha_w(s))\alpha'_w(s)\right\|ds\leq
C_1\int_0^1-\Re(w)e^{\frac{s}{2}\Re(w)}ds=2C_1\left(1-e^{\frac{1}{2}\Re(w)}\right)\leq 2C_1.$$
To prove Inequality \eqref{eq-principal-claim-2},
fix some $\rho$ such that $1<\rho<\eta$.
By Estimate \eqref{eq-principal-1} (which also holds for $|z|=\rho$):
\begin{equation}
\label{eq-principal-2}
\int_{C(0,\rho)}||\xi_t-\xi_0||\leq C |t|.
\end{equation}
By the change of variable formula and the convergence of $\whxi_t$ to $\whxi_0$
on compact subsets of $D^*(0,\eta)$:
$$\int_{C(0,\rho^{-1}|t|)}||\xi_t||=\int_{C(0,\rho)}||\psi^*\xi_t||=\int_{C(0,\rho)}||\whxi_t||\leq C.$$
Since $\xi_0$ is holomorphic in $D(0,\eta)$,
\begin{equation}
\label{eq-principal-3}
\int_{C(0,\rho^{-1}|t|)}||\xi_t-\xi_0||
\leq\int_{C(0,\rho^{-1}|t|)}||\xi_t||+||\xi_0||
\leq C.
\end{equation}
We expand $\xi_t-\xi_0$ in Laurent series in the annulus $\rho^{-1}|t|\leq|z|\leq\rho$ as
$$\xi_t(z)-\xi_0(z)=\sum_{k\in\Z}A_k(t)z^kdz$$
where the matrices $A_k(t)$ are given by
$$A_k(t)=\frac{1}{2\pi\ii}\int_{C(0,\rho)}\frac{\xi_t(z)-\xi_0(z)}{z^{k+1}}
=\frac{1}{2\pi\ii}\int_{C(0,\rho^{-1}|t|)}\frac{\xi_t(z)-\xi_0(z)}{z^{k+1}}$$
Estimates \eqref{eq-principal-2} and \eqref{eq-principal-3} give us respectively:
$$||A_k(t)||\leq\frac{1}{2\pi\rho^{k+1}}\int_{C(0,\rho)}||\xi_t-\xi_0||\leq C\frac{|t|}{\rho^{k+1}}$$
$$||A_k(t)||\leq\frac{\rho^{k+1}}{2\pi |t|^{k+1}}\int_{C(0,\rho^{-1}|t|)}||\xi_t-\xi_0||\leq
C\frac{\rho^{k+1}}{|t|^{k+1}}$$
Hence
\begin{eqnarray*}
\lefteqn{\int_0^1\left\|[\xi_t(\alpha_w(s))-\xi_0(\alpha_w(s))]\alpha'_w(s)\right\|ds
\leq\sum_{k\in\Z}\int_0^1\|A_k(t)\|\cdot|\alpha_w(s)|^k\cdot|\alpha'_w(s)|ds}\\
&\leq&\sum_{k\in\Z}\int_0^1-\|A_k(t)\|e^{\frac{s}{2}(k+1)\Re(w)}\Re(w)ds
\quad\mbox{ using \eqref{eq-principal-4}}\\
&=&-||A_{-1}(t)||\Re(w)+\sum_{k\neq -1}\frac{2}{k+1}||A_k(t)||\left(
1-e^{\frac{1}{2}(k+1)\Re(w)}\right)\\
&\leq &|A_{-1}(t)|\cdot|\log|t|\,|+\sum_{k\geq 0}\frac{2}{k+1}||A_k(t)||+
\sum_{k\leq -2}\frac{-2}{k+1}||A_k(t)||\cdot |t|^{\frac{k+1}{2}}\\
&\leq& C|t\log|t|\,|+\sum_{k\geq 0}\frac{2C|t|}{\rho^{k+1}}
+\sum_{k\leq -2}2C\left|\frac{\rho}{t}\right|^{k+1}|t|^{\frac{k+1}{2}}
\quad\mbox{ using \eqref{eq-principal-2} and \eqref{eq-principal-3}}\\
&\leq& C|t\log|t|\,|+C|t|+C|t|^{\frac{1}{2}}.
\end{eqnarray*}
\cqfd
\medskip

Returning to the proof of Theorem \ref{thm-principal}, let $\Phi_0$ be the solution of
$d\Phi_0=\Phi_0\xi_0$ in $D(0,\eta)$ with initial condition
$\Phi_0(1)=I_n$.
We first estimate the principal solution of $\xi_t$ on the path $\alpha_w$.
Let $Y_w(s)$ be the solution on $[0,1]$ of the Cauchy Problem
$$\left\{\begin{array}{l}
Y'_w(s)=Y_w(s)\xi_t(\alpha_w(s))\alpha'_w(s)\\
Y_w(0)=I_2.\end{array}\right.$$
By definition, $\boP(\xi_t,\alpha_w)=Y_w(1)$.
Define
$$Z_w(s)=Y_w(s)-\Phi_0(\alpha_w(s)).$$
Then
\begin{eqnarray*}
Z'_w(s)&=&Y_w(s)\xi_t(\alpha_w(s))\alpha'_w(s)-\Phi_0(\alpha_w(s))\xi_0(\alpha_w(s))\alpha'_w(s)\\
&=&Z_w(s)\xi_t(\alpha_w(s))\alpha'_w(s)+\Phi_0(\alpha_w(s))[\xi_t(\alpha_w(s))-\xi_0(\alpha_w(s))]\alpha'_w(s).
\end{eqnarray*}
Since $Z_w(0)=0$, we have for all $s\in[0,1]$:
\begin{eqnarray*}
||Z_w(s)||&=&\left\|\int_0^s Z'_w(x)dx\right\|\\
&\leq&\int_0^s||Z_w(x)||\cdot||\xi_t(\alpha_w(x))\alpha'_w(x)||dx
+\int_0^s||\Phi_0(\alpha_w(x))||\cdot\left\|\,[\xi_t(\alpha_w(x))-\xi_0(\alpha_w(x))]\alpha'_w(x)\right\|dx
\end{eqnarray*}
By Gr\"onwall inequality:
$$||Z_w(1)||\leq \int_0^1||\Phi_0(\alpha_w(s))||\cdot\left\|\,[\xi_t(\alpha_w(s))-\xi_0(\alpha_w(s))]\alpha'_w(s)\right\|ds\times\exp\left(\int_0^1||\xi_t(\alpha_w(s))\alpha'_w(s)||ds\right).$$
Since $\Phi_0$ is bounded in $D(0,1)$, Claim \ref{claim-principal} gives
$$||\boP(\xi_t,\alpha_w)-\Phi_0(\alpha_w(1))||=||Z_w(1)||\leq C|t|^{1/2}.$$
Since $\Phi_0$ is holomorphic in $D(0,\eta)$,
$$||\Phi_0(\alpha_w(1))-\Phi_0(0)||\leq C|\alpha_w(1)|=C|t|^{1/2}$$
Hence
\begin{equation}
\label{eq-principal-5}
||\boP(\xi_t,\alpha_w)-\Phi_0(0)||\leq C|t|^{1/2}.
\end{equation}
Let $\whPhi_0$ be the solution of $d\whPhi_0=\whPhi_0\xi_0$ with initial condition
$\Phi_0(1)=I_n$.
Replacing $\xi_t$ by $\whxi_t$, we have:
$$||\boP(\whxi_t,\alpha_w)-\whPhi_0(0)||\leq C|t|^{1/2}.$$
By Equation \eqref{eq-principal-pullback}:
$$\boP(\whxi_t,\alpha_w)=\boP(\psi_t^*\xi_t,\alpha_w)=\boP(\xi_t,\psi_t\circ\alpha_w)
=\boP(\xi_t,\whalpha_w).$$
Hence
\begin{equation}
\label{eq-principal-6}
||\boP(\xi_t,\whalpha_w)-\whPhi_0(0)||\leq C|t|^{1/2}.
\end{equation}
By Equation \eqref{eq-principal-morphism} and Estimates \eqref{eq-principal-5} and \eqref{eq-principal-6}:
$$||\boP(\xi_t,\beta_w)-\Phi_0(0)\whPhi_0(0)^{-1}||=
||\boP(\xi_t,\alpha_w)\boP(\xi_t,\whalpha_w)^{-1}-\Phi_0(0)\whPhi_0(0)^{-1}||
\leq C|t|^{1/2}.$$
Using Estimate \eqref{eq-principal-7}, we finally obtain
$$\left\|F(t)-\Phi_0(0)\whPhi_0(0)^{-1}\right\|\leq C|t|^{1/2}.$$
Hence $F$ extends continuously at $t=0$, and holomorphically by Riemann Extension
Theorem.\cqfd
\section{Convergence to a minimal surface}
\label{appendix-minimal}
The following theorem is proven in \cite{minoids}.
\begin{theorem}
\label{thm-minimal}
Let $\Sigma$ be a Riemann surface and $\wtSigma$ its universal covering. Let $\xi_t$ be a $C^1$ family of holomorphic DPW potentials on
$\Sigma$ and $\Phi_t$ be a continuous family of solutions of $d\Phi_t=\Phi_t\xi_t$
in $\wtSigma$, where $t$ is a real parameter in a neighborhood of $0$.
Assume that $\Phi_t$ solves the Monodromy Problem and
let $f_t=\Sym(\Uni(\Phi_t)):\Sigma\to\R^3$ be the immersion obtained by the DPW method.
Assume that $\Phi_0(z,\lambda)$ is independent of $\lambda$:
$$\Phi_0(z,\lambda)=\matrix{a(z)&b(z)\\c(z)&d(z)}$$
and that
$$\frac{\partial \xi_{t;12}^{(-1)}}{\partial t}|_{t=0}\not\equiv 0.$$
Then 
$$\lim_{t\to 0}\frac{1}{t}f_t(z)=\sigma\circ\psi(z)$$ where
$\sigma$ is the symmetry with respect to the plane $x_2=0$ and
$\psi:\Sigma\to\R^3$ is a minimal (branched) immersion with the following Weierstrass data:
$$g(z)=\frac{c(z)}{a(z)}\qquad
\omega=4a(z)^2\frac{\partial \xi_{t;12}^{(-1)}}{\partial t}|_{t=0}$$
The limit is the uniform $C^1$ convergence on compact subsets of $\Sigma$
minus the zeros of $a$.
\end{theorem}
Observe that since the Monodromy of $\Phi_t$ at $\lambda=1$ is $I_2$, $\Phi_0(z)$
descends to a well defined map in $\Sigma$, so $g(z)$ is a well defined 
meromorphic function on $\Sigma$.
\section{Differentiability of smooth functions of $t$ and $t\log t$}
\label{appendix-tlogt}
\begin{proposition}
\label{prop-tlogt}
Let $g(t,s,z)$ be a smooth function of $(t,s)$ in a neighborhood of $(0,0)$ in $\R^2$
and $z\in\Omega\subset E$ where $E$ is a finite dimensional normed space,
with values in a Banach space $\boB$.
Define
$$f(t,z)=\left\{\begin{array}{ll}
g(t,t\log|t|,z)&\mbox{ if $t\neq 0$}\\
g(0,0,z)&\mbox{ if $t=0$.}\end{array}\right.$$
Assume that $g(0,s,z)$ does not depend on $s$.
Then $f$ is of class $C^1$.
\end{proposition}
Proof: let $z_0\in\Omega$. We have:
$$g(t,s,z)=g(0,0,z_0)+\frac{\partial g}{\partial t}(0,0,z_0)t+d_zg(0,0,z_0)(z-z_0)
+O(t^2+s^2+||z-z_0||^2)$$
$$f(t,z)=f(0,z_0)+\frac{\partial g}{\partial t}(0,0,z_0)t+d_zg(0,0,z_0)(z-z_0)
+O( (t\log |t|)^2+||z-z_0||^2).$$
Since $O( (t\log |t|)^2)=o(t)$, $f$ is differentiable at $(0,z_0)$ with
$$df(0,z_0)=\frac{\partial g}{\partial t}(0,0,z_0)dt+d_zg(0,0,z_0).$$
For $t\neq 0$, we have by the chain rule:
$$df(t,z)=\frac{\partial g}{\partial t}(t,t\log |t|,z)dt+
\frac{\partial g}{\partial s}(t,t\log|t|,z)(1+\log|t|)dt+
d_zg(t,t\log|t|,z).$$
By the mean value inequality:
$$||\frac{\partial g}{\partial s}(t,s,z)||=
||\frac{\partial g}{\partial s}(t,s,z)-\frac{\partial g}{\partial s}(0,s,z)||
\leq C|t|.$$
Hence
$$\lim_{(t,z)\to(0,z_0)}df(t,z)=df(0,z_0).$$
\cqfd

\bigskip

\noindent
Martin Traizet\\
Institut Denis Poisson\\
Universit\'e de Tours, 37200 Tours, France\\
\verb$martin.traizet@univ-tours.fr$

\begin{thebibliography}{9}
\bibitem{chae}
S. B. Chae. {\em Holomorphy and calculus in normed spaces.}
Monographs and textbooks in pure and applied mathematics, vol. 92 (1985).
\bibitem{dorfmeister-haak}
J. Dorfmeister, G. Haak:
{\em Constant mean curvature surfaces with periodic metric.}
Pacific Journal of Mathematics 182 (1998), 229--287.
\bibitem{dorfmeister-pedit-wu}
J. Dorfmeister, F. Pedit, H. Wu:
{\em Weierstrass type representation of
harmonic maps into symmetric spaces.}
Communications in Analysis and Geometry 6 (1998), 633-668.
\bibitem{dorfmeister-wu}
J. Dorfmeister, H. Wu:
{\em Construction of constant mean curvature n-noids
from holomorphic potentials.}
Mathematische Zeitschrift 258 (2008), 773--803.
\bibitem{fay} J.D. Fay: {\em Theta Functions on Riemann Surfaces}. Lecture Notes 
in Mathematics 352 (1973).
\bibitem{forster}
O. Forster:
{\em Lectures on Riemann surfaces.}
Graduate texts in Mathematics, Springer Verlag.
\bibitem{fujimori-kobayashi-rossman}
S. Fujimori, S. Kobayashi, W. Rossman:
{\em Loop group methods for constant mean curvature surfaces.}
arXiv:math/0602570.
\bibitem{gerding-pedit-schmitt}
A. Gerding, F. Pedit, N. Schmitt:
{\em Constant mean curvature surfaces: an integrable systems perspective.}
Harmonic maps and differential geometry,
Contemp. Math. 542 (2011), Amer. Math. Soc., 7--39.
\bibitem{karsten-kusner-sullivan}
K. Gro\ss e-Brauckmann, R. Kusner, J. Sullivan:
{\em Triunduloids: embedded constant mean curvature surfaces with three ends and genus zero.}
J. Reine Angew. Math. 564 (2003), 35--61.
\bibitem{karsten-kusner-sullivan2}
K. Gro\ss e-Brauckmann, R. Kusner, J. Sullivan:
{\em Coplanar constant mean curvature surfaces.}
Comm. Anal. Geom. 15 (2007), no. 5, 985--1023.
\bibitem{heller-heller-schmitt}
L. Heller, S. Heller, N. Schmitt:
{\em Navigating the Space of Symmetric CMC Surfaces.}
arXiv:1501.01929.
\bibitem{heller1}
S. Heller:
{\em Higher genus minimal surfaces in $\S^3$ and stable bundles.}
J. Reine Angew. Math. 685 (2013), 105--122. 
\bibitem{heller2}
S. Heller:
{\em Lawson's genus two surface and meromorphic connections.}
Mathematische Zeitschrift 274 (2013), 745--760.
\bibitem{heller3}
S. Heller:
{\em A spectral curve approach to Lawson symmetric CMC surfaces of genus 2.}
Math. Annalen 360, Issue 3 (2014), 607--652. 
\bibitem{kapouleas}
N. Kapouleas:
{\em Complete constant mean curvature surfaces in euclidean three-space.}
Annals of Mathematics 131 (1990), 239--330.
\bibitem{kilian-kobayashi-rossman-schmitt}
M. Kilian, S. Kobayashi, W. Rossman, N. Schmitt:
{\em Constant mean curvature surfaces of any positive genus.}
J. London Math. Soc. 72 (2005), 258--272.
\bibitem{kilian-mcintosh-schmitt}
M. Kilian, I. McIntosh, N. Schmitt:
{\em New constant mean curvature surfaces.}
Experiment. Math. 9 (2000), 595--611.
\bibitem{kilian-rossman-schmitt}
M. Kilian, W. Rossman, N. Schmitt:
{\em Delaunay ends of constant mean curvature surfaces.}
Compositio Mathematica 144 (2008), 186--220.
\bibitem{korevaar-kusner-solomon}
N. Korevaar, R. Kusner, B. Solomon:
{\em The structure of complete embedded surfaces with constant mean curvature.}
J. Diff. Geom. 30 (1989), 465--503.
\bibitem{kusner-mazzeo-pollack}
R. Kusner, R. Mazzeo, D. Pollack:
{\em The moduli space of complete embedded constant mean curvature surfaces.}
GAFA 6, Issue 1 (1996), 120--137.
\bibitem{masur} H. Masur: {\em The extension of the Weil-Petersson metric to the 
boundary of Teichmuller space}. Duke Math. J. 43  (1976), 623--635.
\bibitem{raujouan}
T. Raujouan:
{\em On Delaunay ends in the DPW method.}
arXiv:1710.00768.
\bibitem{schmitt}
N. Schmitt:
{\em Constant mean curvature $n$-noids with platonic symmetries.}
arXiv:math/0702469.
\bibitem{schmitt-kilian-kobayashi-rossman}
N. Schmitt, M. Kilian, S. Kobayashi, W. Rossman:
{\em Unitarization of monodromy
representations and constant mean curvature trinoids in 3-dimensional space
forms.}
Journal of the London Mathematical Society 75 (2007), 563--581.
\bibitem{taylor} M. Taylor: {\em Introduction to Differential Equations.}
Pure and Applied Undergraduate Texts 14, American Math. Soc. (2011).
\bibitem{teschl} G. Teschl: Ordinary differential equations and dynamical systems.
{\em Graduate Studies in Mathematics} 140, American Math. Soc. (2010).
\bibitem{nosym}
M. Traizet: An embedded minimal surface with no symmetries.
{\em J. Differential Geometry}, 60(1) (2002), 103--153.
\bibitem{triply} M. Traizet:
{\em On the genus of triply periodic minimal surfaces.}
J. Diff. Geom. 79 (2008), 243--275.
\bibitem{crelle} M. Traizet:
{\em Opening infinitely many nodes.}
J. reine angew. Math. 684 (2013), 165--186.
\bibitem{bryant}
M. Traizet: 
{\em Opening nodes on horosphere packings.}
Trans. Amer. Math. Soc. 368 (2016), 5701--5725.
\bibitem{nnoids}
M. Traizet:
{\em Construction of constant mean curvature $n$-noids using
the DPW method.}
arXiv 1709.00924 (2017).
\bibitem{minoids}
M. Traizet:
{\em Gluing Delaunay ends to minimal n-noids using the DPW method.}
arXiv 1710.09261 (2017).
\end{thebibliography}
\end{document}